\documentclass[12pt]{amsart}
\usepackage{amsfonts,amsthm,amssymb,amscd}
\usepackage{afterpage} 
\usepackage{lscape}

\advance\textheight\topskip
\usepackage[OT2,T1]{fontenc}

\usepackage[mathscr]{eucal}
\usepackage{mathabx}

\usepackage{times}
\usepackage{mathptmx}
\newcommand{\myboldsymbol}{\pmb}
\allowdisplaybreaks
\usepackage{t-angles}

\usepackage[matrix,arrow,curve,dvips]{xy} 
\UsePSspecials{dvips}
\usepackage[dvips, bookmarks=false, hyperfigures=false,
setpagesize=false]{hyperref}

\DeclareMathAlphabet{\mathpzc}{OT1}{pzc}{m}{it}

\allowdisplaybreaks

\numberwithin{equation}{section}
\numberwithin{figure}{section}
\makeatletter
\@addtoreset{equation}{section}
\@addtoreset{figure}{section}
\@addtoreset{subsubsection}{section}

\def\@secnumfont{\bfseries}
\def\subsubsection{\@startsection{subsubsection}{3}%
  \z@{.5\linespacing\@plus.7\linespacing}{-.5em}%
  {\normalfont\bfseries}}
\def\paragraph{\@startsection{paragraph}{4}%
  \z@\z@{-\fontdimen2\font}%
  \normalfont\bfseries}
\def\subparagraph{\@startsection{subparagraph}{5}%
  \z@\z@{-\fontdimen2\font}%
  \normalfont\bfseries}
\makeatother

\newcommand{\bRi}{\bar{R}\up{1}}
\newcommand{\bRii}{\bar{R}\up{2}}
\newcommand{\Ri}{R\up{1}}
\newcommand{\Rii}{R\up{2}}
\newcommand{\Pone}{P\up{1}}
\newcommand{\Pii}{P\up{2}}
\newcommand{\Uminus}{U_{<}}

\newcommand{\pii}{{\pi^{(1)}}}
\newcommand{\piii}{{\pi^{(2)}}}
\newcommand{\rhoi}{{\rho^{(1)}}}
\newcommand{\rhoii}{{\rho^{(2)}}}
\newcommand{\sigmai}{{\sigma^{(1)}}}
\newcommand{\sigmaii}{{\sigma^{(2)}}}
\newcommand{\taui}{{\tau^{(1)}}}
\newcommand{\tauii}{{\tau^{(2)}}}

\newcommand{\tDelta}{\widetilde{\Delta}}

\newcommand{\drmap}{\boldsymbol{\chi}}
\newcommand{\Zcenter}{\mathfrak{Z}}
\newcommand{\Ch}{\mathfrak{Ch}}

\newcommand{\TOP}{{}^{{\scriptscriptstyle\bigtriangleup}}}
\newcommand{\BOT}{{}_{{\scriptscriptstyle\bigtriangledown}}}

\newcommand{\est}{\mathbf{e}} 
\newcommand{\etyp}{\mathbf{e}} 
\newcommand{\wtip}{\mathbf{w}} 
\newcommand{\watypb}{\mathbf{w}} 
\newcommand{\watypm}{\mathbf{w}} 
\newcommand{\Watypm}{\mathbf{W}} 

\newcommand{\ZC}[1]{\ket{{\textstyle#1}}^{\!\leftarrow}}
\newcommand{\ZB}[1]{\ket{{\textstyle#1}}^{\!\rightarrow}}
\newcommand{\ZD}[1]{\ket{{\textstyle#1}}^{\downarrow}}
\newcommand{\ZU}[1]{\ket{{\textstyle#1}}^\uparrow}

\newcommand{\something}[1]{*{#1}\ar@{};[0,0];}
\newcommand{\myatop}[2]{{#1\atop\mbox{\small$(#2)$}}}

\newcommand{\zpo}{\boxminus}
\newcommand{\zpp}{\boxslash}
\newcommand{\zso}{\boxplus}
\newcommand{\zss}{\boxtimes}
\newcommand{\zoz}{\boxasterisk}

\newcommand{\Size}{\Large}
\newcommand{\Zpo}{\mbox{\Size$\zpo$}}
\newcommand{\Zpp}{\mbox{\Size$\zpp$}}
\newcommand{\Zso}{\mbox{\Size$\zso$}}
\newcommand{\Zss}{\mbox{\Size$\zss$}}
\newcommand{\Zoz}{\mbox{\Size$\zoz$}}
\newcommand{\Zpr}{\mbox{\Size${\boxdot}$}}

\newcommand{\cpo}{\varpi}
\newcommand{\cpp}{\pi}
\newcommand{\cso}{\omega}
\newcommand{\css}{\sigma}
\newcommand{\coz}{\zeta}

\newcommand{\mylabel}[3]{{#2_{#3}}} 
\newcommand{\vmylabel}[3]{{({#3})}} 
\newcommand{\arlabel}[1]{\mbox{\rule[-5pt]{5pt}{0pt}\rule[-5pt]{0pt}{12pt}\kern-8pt\larger[-6]$#1$}}

\newcommand{\myxy}[1]{\xymatrix@R=15pt@C=0pt{#1}}

\newcommand{\mysmallmatrix}[1]{{\left(\!\begin{smallmatrix}#1\end{smallmatrix}\!\right)}}

\newcommand{\poiApBZ}{{\left(\!\begin{smallmatrix}\alpha &\beta \\p&0\end{smallmatrix}\!\right)_{\!\!\mylabel{}{1}{0}}\kern-4pt}\ar@{};[0,0];}
\newcommand{\poiieemAIBZ}{{\left(\!\begin{smallmatrix}-\alpha &\beta \\1&0\end{smallmatrix}\!\right)_{\!\!\mylabel{ee}{2}{0}}\kern-4pt}\ar@{};[0,0];}
\newcommand{\poiiefIApmimBZ}{{\left(\!\begin{smallmatrix}\alpha &-\beta \\p-1&0\end{smallmatrix}\!\right)_{\!\!\mylabel{ef1}{2}{1}}\kern-4pt}\ar@{};[0,0];}
\newcommand{\poiiffmAImBZ}{{\left(\!\begin{smallmatrix}-\alpha &-\beta \\1&0\end{smallmatrix}\!\right)_{\!\!\mylabel{ff}{2}{2}}\kern-4pt}\ar@{};[0,0];}
\newcommand{\poiiiCcApBZ}{{\left(\!\begin{smallmatrix}\alpha &\beta \\p&0\end{smallmatrix}\!\right)_{\!\!\mylabel{Cc}{3}{0}}\kern-4pt}\ar@{};[0,0];}
\newcommand{\poiiiEImAIBI}{{\left(\!\begin{smallmatrix}-\alpha &\beta \\1&1\end{smallmatrix}\!\right)_{\!\!\mylabel{E1}{3}{1}}\kern-4pt}\ar@{};[0,0];}
\newcommand{\poiiiEfcApmiBpmi}{{\left(\!\begin{smallmatrix}\alpha &\beta \\p-1&p-1\end{smallmatrix}\!\right)_{\!\!\mylabel{Efc}{3}{2}}\kern-4pt}\ar@{};[0,0];}
\newcommand{\poiiiFImAImBI}{{\left(\!\begin{smallmatrix}-\alpha &-\beta \\1&1\end{smallmatrix}\!\right)_{\!\!\mylabel{F1}{3}{3}}\kern-4pt}\ar@{};[0,0];}
\newcommand{\poivBBApmimBZ}{{\left(\!\begin{smallmatrix}\alpha &-\beta \\p-1&0\end{smallmatrix}\!\right)_{\!\!\mylabel{BB}{4}{1}}\kern-4pt}\ar@{};[0,0];}
\newcommand{\poivEEmAIBZ}{{\left(\!\begin{smallmatrix}-\alpha &\beta \\1&0\end{smallmatrix}\!\right)_{\!\!\mylabel{EE}{4}{0}}\kern-4pt}\ar@{};[0,0];}
\newcommand{\poivFFmAImBZ}{{\left(\!\begin{smallmatrix}-\alpha &-\beta \\1&0\end{smallmatrix}\!\right)_{\!\!\mylabel{FF}{4}{2}}\kern-4pt}\ar@{};[0,0];}
\newcommand{\povbotApBZ}{{\left(\!\begin{smallmatrix}\alpha &\beta \\p&0\end{smallmatrix}\!\right)_{\!\!\mylabel{bot}{5}{0}}\kern-4pt}\ar@{};[0,0];}

\newcommand{\poiieemAIBZm}{{\left(\!\begin{smallmatrix}\alpha &\beta \\1&0\end{smallmatrix}\!\right)_{\!\!\mylabel{ee}{2}{0}}\kern-4pt}\ar@{};[0,0];}
\newcommand{\poivEEmAIBZm}{{\left(\!\begin{smallmatrix}\alpha &\beta \\1&0\end{smallmatrix}\!\right)_{\!\!\mylabel{EE}{4}{0}}\kern-4pt}\ar@{};[0,0];}

\newcommand{\ppiApmiBpmi}{{\left(\!\begin{smallmatrix}\alpha &\beta \\p-1&p-1\end{smallmatrix}\!\right)_{\!\!\mylabel{}{1}{0}}\kern-4pt}\ar@{};[0,0];}
\newcommand{\ppiiccApmiimBpmii}{{\left(\!\begin{smallmatrix}\alpha &-\beta \\p-2&p-2\end{smallmatrix}\!\right)_{\!\!\mylabel{cc}{2}{0}}\kern-4pt}\ar@{};[0,0];}
\newcommand{\ppiieemAIBZ}{{\left(\!\begin{smallmatrix}-\alpha &\beta \\1&0\end{smallmatrix}\!\right)_{\!\!\mylabel{ee}{2}{1}}\kern-4pt}\ar@{};[0,0];}
\newcommand{\ppiiffmAImBZ}{{\left(\!\begin{smallmatrix}-\alpha &-\beta \\1&0\end{smallmatrix}\!\right)_{\!\!\mylabel{ff}{2}{2}}\kern-4pt}\ar@{};[0,0];}
\newcommand{\ppiiiEemAIImBZ}{{\left(\!\begin{smallmatrix}-\alpha &-\beta \\2&0\end{smallmatrix}\!\right)_{\!\!\mylabel{Ee}{3}{0}}\kern-4pt}\ar@{};[0,0];}
\newcommand{\ppiiiEfbApmiBpmi}{{\left(\!\begin{smallmatrix}\alpha &\beta \\p-1&p-1\end{smallmatrix}\!\right)_{\!\!\mylabel{Efb}{3}{1}}\kern-4pt}\ar@{};[0,0];}
\newcommand{\ppiiiFfmAIIBZ}{{\left(\!\begin{smallmatrix}-\alpha &\beta \\2&0\end{smallmatrix}\!\right)_{\!\!\mylabel{Ff}{3}{2}}\kern-4pt}\ar@{};[0,0];}
\newcommand{\ppiiiBbApBZ}{{\left(\!\begin{smallmatrix}\alpha &\beta \\p&0\end{smallmatrix}\!\right)_{\!\!\mylabel{Bb}{3}{3}}\kern-4pt}\ar@{};[0,0];}
\newcommand{\ppivEEmAIBZ}{{\left(\!\begin{smallmatrix}-\alpha &\beta \\1&0\end{smallmatrix}\!\right)_{\!\!\mylabel{EE}{4}{1}}\kern-4pt}\ar@{};[0,0];}
\newcommand{\ppivFFmAImBZ}{{\left(\!\begin{smallmatrix}-\alpha &-\beta \\1&0\end{smallmatrix}\!\right)_{\!\!\mylabel{FF}{4}{2}}\kern-4pt}\ar@{};[0,0];}
\newcommand{\ppivCCApmiimBpmii}{{\left(\!\begin{smallmatrix}\alpha &-\beta \\p-2&p-2\end{smallmatrix}\!\right)_{\!\!\mylabel{CC}{4}{0}}\kern-4pt}\ar@{};[0,0];}
\newcommand{\ppvbotApmiBpmi}{{\left(\!\begin{smallmatrix}\alpha &\beta \\p-1&p-1\end{smallmatrix}\!\right)_{\!\!\mylabel{bot}{5}{0}}\kern-4pt}\ar@{};[0,0];}

\newcommand{\soiAsBZ}{{\left(\!\begin{smallmatrix}\alpha &\beta \\s&0\end{smallmatrix}\!\right)_{\!\!\mylabel{}{1}{0}}\kern-4pt}\ar@{};[0,0];}
\newcommand{\soiietmApmsmBpms}{{\left(\!\begin{smallmatrix}-\alpha &-\beta \\p-s&p-s\end{smallmatrix}\!\right)_{\!\!\mylabel{et}{2}{0}}\kern-4pt}\ar@{};[0,0];}
\newcommand{\soiiftmApmsBpms}{{\left(\!\begin{smallmatrix}-\alpha &\beta \\p-s&p-s\end{smallmatrix}\!\right)_{\!\!\mylabel{ft}{2}{1}}\kern-4pt}\ar@{};[0,0];}
\newcommand{\soiitlAsmimBZ}{{\left(\!\begin{smallmatrix}\alpha &-\beta \\s-1&0\end{smallmatrix}\!\right)_{\!\!\mylabel{tl}{2}{2}}\kern-4pt}\ar@{};[0,0];}
\newcommand{\soiitrAsimBZ}{{\left(\!\begin{smallmatrix}\alpha &-\beta \\s+1&0\end{smallmatrix}\!\right)_{\!\!\mylabel{tr}{2}{3}}\kern-4pt}\ar@{};[0,0];}
\newcommand{\soiiiBtAsBZ}{{\left(\!\begin{smallmatrix}\alpha &\beta \\s&0\end{smallmatrix}\!\right)_{\!\!\mylabel{Bt}{3}{0}}\kern-4pt}\ar@{};[0,0];}
\newcommand{\soiiiElmApimsBpims}{{\left(\!\begin{smallmatrix}-\alpha &\beta \\p-s+1&p-s+1\end{smallmatrix}\!\right)_{\!\!\mylabel{El}{3}{1}}\kern-4pt}\ar@{};[0,0];}
\newcommand{\soiiiErmApmsmiBpmsmi}{{\left(\!\begin{smallmatrix}-\alpha &\beta \\p-s-1&p-s-1\end{smallmatrix}\!\right)_{\!\!\mylabel{Er}{3}{2}}\kern-4pt}\ar@{};[0,0];}
\newcommand{\soiiiFlmApimsmBpims}{{\left(\!\begin{smallmatrix}-\alpha &-\beta \\p-s+1&p-s+1\end{smallmatrix}\!\right)_{\!\!\mylabel{Fl}{3}{3}}\kern-4pt}\ar@{};[0,0];}
\newcommand{\soiiiFrmApmsmimBpmsmi}{{\left(\!\begin{smallmatrix}-\alpha &-\beta \\p-s-1&p-s-1\end{smallmatrix}\!\right)_{\!\!\mylabel{Fr}{3}{4}}\kern-4pt}\ar@{};[0,0];}
\newcommand{\soiiiTbAsBZ}{{\left(\!\begin{smallmatrix}\alpha &\beta \\s&0\end{smallmatrix}\!\right)_{\!\!\mylabel{Tb}{3}{5}}\kern-4pt}\ar@{};[0,0];}
\newcommand{\soivBLAsmimBZ}{{\left(\!\begin{smallmatrix}\alpha &-\beta \\s-1&0\end{smallmatrix}\!\right)_{\!\!\mylabel{BL}{4}{0}}\kern-4pt}\ar@{};[0,0];}
\newcommand{\soivBRAsimBZ}{{\left(\!\begin{smallmatrix}\alpha &-\beta \\s+1&0\end{smallmatrix}\!\right)_{\!\!\mylabel{BR}{4}{1}}\kern-4pt}\ar@{};[0,0];}
\newcommand{\soivEBmApmsmBpms}{{\left(\!\begin{smallmatrix}-\alpha &-\beta \\p-s&p-s\end{smallmatrix}\!\right)_{\!\!\mylabel{EB}{4}{2}}\kern-4pt}\ar@{};[0,0];}
\newcommand{\soivFBmApmsBpms}{{\left(\!\begin{smallmatrix}-\alpha &\beta \\p-s&p-s\end{smallmatrix}\!\right)_{\!\!\mylabel{FB}{4}{3}}\kern-4pt}\ar@{};[0,0];}
\newcommand{\sovbotAsBZ}{{\left(\!\begin{smallmatrix}\alpha &\beta \\s&0\end{smallmatrix}\!\right)_{\!\!\mylabel{bot}{5}{0}}\kern-4pt}\ar@{};[0,0];}

\newcommand{\ssiAsBs}{{\left(\!\begin{smallmatrix}\alpha &\beta \\s&s\end{smallmatrix}\!\right)_{\!\!\mylabel{}{1}{0}}\kern-4pt}\ar@{};[0,0];}
\newcommand{\ssiibcemApmsBZ}{{\left(\!\begin{smallmatrix}-\alpha &\beta \\p-s&0\end{smallmatrix}\!\right)_{\!\!\mylabel{bce}{2}{0}}\kern-4pt}\ar@{};[0,0];}
\newcommand{\ssiibcefAsmimBsmi}{{\left(\!\begin{smallmatrix}\alpha &-\beta \\s-1&s-1\end{smallmatrix}\!\right)_{\!\!\mylabel{bcef}{2}{2}}\kern-4pt}\ar@{};[0,0];}
\newcommand{\ssiibcfmApmsmBZ}{{\left(\!\begin{smallmatrix}-\alpha &-\beta \\p-s&0\end{smallmatrix}\!\right)_{\!\!\mylabel{bcf}{2}{1}}\kern-4pt}\ar@{};[0,0];}
\newcommand{\ssiibefAsimBsi}{{\left(\!\begin{smallmatrix}\alpha &-\beta \\s+1&s+1\end{smallmatrix}\!\right)_{\!\!\mylabel{bef}{2}{3}}\kern-4pt}\ar@{};[0,0];}
\newcommand{\ssiiiBcAsBs}{{\left(\!\begin{smallmatrix}\alpha &\beta \\s&s\end{smallmatrix}\!\right)_{\!\!\mylabel{Bc}{3}{0}}\kern-4pt}\ar@{};[0,0];}
\newcommand{\ssiiiBemApmsmimBZ}{{\left(\!\begin{smallmatrix}-\alpha &-\beta \\p-s-1&0\end{smallmatrix}\!\right)_{\!\!\mylabel{Be}{3}{2}}\kern-4pt}\ar@{};[0,0];}
\newcommand{\ssiiiCemApimsmBZ}{{\left(\!\begin{smallmatrix}-\alpha &-\beta \\p-s+1&0\end{smallmatrix}\!\right)_{\!\!\mylabel{Ce}{3}{1}}\kern-4pt}\ar@{};[0,0];}
\newcommand{\ssiiiCefAsBs}{{\left(\!\begin{smallmatrix}\alpha &\beta \\s&s\end{smallmatrix}\!\right)_{\!\!\mylabel{Cef}{3}{5}}\kern-4pt}\ar@{};[0,0];}
\newcommand{\ssiiiCfmApimsBZ}{{\left(\!\begin{smallmatrix}-\alpha &\beta \\p-s+1&0\end{smallmatrix}\!\right)_{\!\!\mylabel{Cf}{3}{3}}\kern-4pt}\ar@{};[0,0];}
\newcommand{\ssiiiBfmApmsmiBZ}{{\left(\!\begin{smallmatrix}-\alpha &\beta \\p-s-1&0\end{smallmatrix}\!\right)_{\!\!\mylabel{Bf}{3}{4}}\kern-4pt}\ar@{};[0,0];}
\newcommand{\ssivBBAsimBsi}{{\left(\!\begin{smallmatrix}\alpha &-\beta \\s+1&s+1\end{smallmatrix}\!\right)_{\!\!\mylabel{BB}{4}{1}}\kern-4pt}\ar@{};[0,0];}
\newcommand{\ssivCCAsmimBsmi}{{\left(\!\begin{smallmatrix}\alpha &-\beta \\s-1&s-1\end{smallmatrix}\!\right)_{\!\!\mylabel{CC}{4}{0}}\kern-4pt}\ar@{};[0,0];}
\newcommand{\ssivEEmApmsBZ}{{\left(\!\begin{smallmatrix}-\alpha &\beta \\p-s&0\end{smallmatrix}\!\right)_{\!\!\mylabel{EE}{4}{2}}\kern-4pt}\ar@{};[0,0];}
\newcommand{\ssivFFmApmsmBZ}{{\left(\!\begin{smallmatrix}-\alpha &-\beta \\p-s&0\end{smallmatrix}\!\right)_{\!\!\mylabel{FF}{4}{3}}\kern-4pt}\ar@{};[0,0];}
\newcommand{\ssvbotAsBs}{{\left(\!\begin{smallmatrix}\alpha &\beta \\s&s\end{smallmatrix}\!\right)_{\!\!\mylabel{bot}{5}{0}}\kern-4pt}\ar@{};[0,0];}
\newcommand{\oziAIBZ}{{\kern-6pt\left(\!\begin{smallmatrix}\alpha &\beta \\1&0\end{smallmatrix}\!\right)_{\!\!\mylabel{}{1}{0}}\kern-8pt}\ar@{};[0,0];}
\newcommand{\oziibIIIAIBI}{{\left(\!\begin{smallmatrix}\alpha &\beta \\1&1\end{smallmatrix}\!\right)_{\!\!\mylabel{b3}{2}{0}}\kern-4pt}\ar@{};[0,0];}
\newcommand{\oziicIVAIImBZ}{{\left(\!\begin{smallmatrix}\alpha &-\beta \\2&0\end{smallmatrix}\!\right)_{\!\!\mylabel{c4}{2}{5}}\kern-4pt}\ar@{};[0,0];}
\newcommand{\oziieIImApmBZ}{{\left(\!\begin{smallmatrix}-\alpha &-\beta \\p&0\end{smallmatrix}\!\right)_{\!\!\mylabel{e2}{2}{1}}\kern-4pt}\ar@{};[0,0];}
\newcommand{\oziieVImApmimBpmi}{{\left(\!\begin{smallmatrix}-\alpha &-\beta \\p-1&p-1\end{smallmatrix}\!\right)_{\!\!\mylabel{e6}{2}{3}}\kern-4pt}\ar@{};[0,0];}
\newcommand{\oziifImApBZ}{{\left(\!\begin{smallmatrix}-\alpha &\beta \\p&0\end{smallmatrix}\!\right)_{\!\!\mylabel{f1}{2}{2}}\kern-4pt}\ar@{};[0,0];}
\newcommand{\oziifVmApmiBpmi}{{\left(\!\begin{smallmatrix}-\alpha &\beta \\p-1&p-1\end{smallmatrix}\!\right)_{\!\!\mylabel{f5}{2}{4}}\kern-4pt}\ar@{};[0,0];}
\newcommand{\oziiibImApmimBZ}{{\left(\!\begin{smallmatrix}-\alpha &-\beta \\p-1&0\end{smallmatrix}\!\right)_{\!\!\mylabel{b1}{3}{1}}\kern-4pt}\ar@{};[0,0];}
\newcommand{\oziiibIImApmiBZ}{{\left(\!\begin{smallmatrix}-\alpha &\beta \\p-1&0\end{smallmatrix}\!\right)_{\!\!\mylabel{b2}{3}{0}}\kern-4pt}\ar@{};[0,0];}
\newcommand{\oziiisVAIBZ}{{\left(\!\begin{smallmatrix}\alpha &\beta \\1&0\end{smallmatrix}\!\right)_{\!\!\mylabel{s5}{3}{2}}\kern-4pt}\ar@{};[0,0];}
\newcommand{\oziiicIZmApmiiBpmii}{{\left(\!\begin{smallmatrix}-\alpha &\beta \\p-2&p-2\end{smallmatrix}\!\right)_{\!\!\mylabel{c10}{3}{8}}\kern-4pt}\ar@{};[0,0];}
\newcommand{\oziiicIXmApmiimBpmii}{{\left(\!\begin{smallmatrix}-\alpha &-\beta \\p-2&p-2\end{smallmatrix}\!\right)_{\!\!\mylabel{c9}{3}{9}}\kern-4pt}\ar@{};[0,0];}
\newcommand{\oziiisVIAIBZ}{{\left(\!\begin{smallmatrix}\alpha &\beta \\1&0\end{smallmatrix}\!\right)_{\!\!\mylabel{s6}{3}{5}}\kern-4pt}\ar@{};[0,0];}
\newcommand{\oziiisVIIAIBZ}{{\left(\!\begin{smallmatrix}\alpha &\beta \\1&0\end{smallmatrix}\!\right)_{\!\!\mylabel{s7}{3}{7}}\kern-4pt}\ar@{};[0,0];}
\newcommand{\oziiiltVIIIAImBZ}{{\left(\!\begin{smallmatrix}\alpha &-\beta \\1&0\end{smallmatrix}\!\right)_{\!\!\mylabel{lt8}{3}{3}}\kern-4pt}\ar@{};[0,0];}
\newcommand{\oziiimtIVAIBZ}{{\left(\!\begin{smallmatrix}\alpha &\beta \\1&0\end{smallmatrix}\!\right)_{\!\!\mylabel{mt4}{3}{4}}\kern-4pt}\ar@{};[0,0];}
\newcommand{\oziiihtIIIAImBZ}{{\left(\!\begin{smallmatrix}\alpha &-\beta \\1&0\end{smallmatrix}\!\right)_{\!\!\mylabel{ht3}{3}{6}}\kern-4pt}\ar@{};[0,0];}
\newcommand{\ozivBIIIAIBI}{{\left(\!\begin{smallmatrix}\alpha &\beta \\1&1\end{smallmatrix}\!\right)_{\!\!\mylabel{B3}{4}{1}}\kern-4pt}\ar@{};[0,0];}
\newcommand{\ozivFImApBZ}{{\left(\!\begin{smallmatrix}-\alpha &\beta \\p&0\end{smallmatrix}\!\right)_{\!\!\mylabel{F1}{4}{2}}\kern-4pt}\ar@{};[0,0];}
\newcommand{\ozivEIImApmBZ}{{\left(\!\begin{smallmatrix}-\alpha &-\beta \\p&0\end{smallmatrix}\!\right)_{\!\!\mylabel{E2}{4}{0}}\kern-4pt}\ar@{};[0,0];}
\newcommand{\ozivCIVAIImBZ}{{\left(\!\begin{smallmatrix}\alpha &-\beta \\2&0\end{smallmatrix}\!\right)_{\!\!\mylabel{C4}{4}{5}}\kern-4pt}\ar@{};[0,0];}
\newcommand{\ozivEVImApmimBpmi}{{\left(\!\begin{smallmatrix}-\alpha &-\beta \\p-1&p-1\end{smallmatrix}\!\right)_{\!\!\mylabel{E6}{4}{3}}\kern-4pt}\ar@{};[0,0];}
\newcommand{\ozivFVmApmiBpmi}{{\left(\!\begin{smallmatrix}-\alpha &\beta \\p-1&p-1\end{smallmatrix}\!\right)_{\!\!\mylabel{F5}{4}{4}}\kern-4pt}\ar@{};[0,0];}
\newcommand{\ozvbotAIBZ}{{\kern-6pt\left(\!\begin{smallmatrix}\alpha &\beta \\1&0\end{smallmatrix}\!\right)_{\!\!\mylabel{bot}{5}{0}}\kern-8pt}\ar@{};[0,0];}

\newcommand{\Exti}{\mathop{\mathrm{Ext}^1}}

\newcommand{\KK}{\mathcal{K}\kern-5.7pt\raisebox{-3.9pt}{\footnotesize\textit{2}}\,}

\newcommand{\ccirc}{\mathbin{\raisebox{1pt}{\,$\scriptscriptstyle\circ$\,}}}

\newcommand{\bExti}{\boldsymbol{\textup{Ext}}^{\myboldsymbol{1}}}

\newcommand{\Ures}{\overline{\mathscr{U}}_{\q} s\ell(2)}
\newcommand{\Uresk}{\overline{\mathscr{U}}^{*}_{\q} s\ell(2)}

\newcommand{\algU}{{\mathbf{U}(X)}}

\newcommand{\balgU}{{\mathbf{U}\myboldsymbol{(X)}}}

\newcommand{\brepQ}{\bar{\repQ}}

\newcommand{\repP}{\mathscr{P}}
\newcommand{\repQ}{\mathscr{Q}}

\newcommand{\repX}{\mathscr{X}}
\newcommand{\repY}{\mathscr{Y}}
\newcommand{\repZ}{\mathscr{Z}}

\newcommand{\dF}{\widehat{F}}

\newcommand{\bNich}{\boldsymbol{\mathfrak{B}}}
\newcommand{\Nich}{\mathfrak{B}}

\newcommand{\symm}{\mathrm{s}}
\newcommand{\asymm}{\mathrm{a}}

\newcommand{\Fi}{B}
\newcommand{\Fii}{F}

\newcommand{\one}{\boldsymbol{1}}

\newcommand{\zDelta}{\boldsymbol{\Delta}}
\newcommand{\zS}{\boldsymbol{S}}

\newcommand{\DD}{\mathscr{D}}

\newcommand{\fobject}[1]{\object{\scriptstyle #1}}

\newcommand{\dotact}{\mathbin{\pmb{.}}}

\newcommand{\category}{\mathscr}
\newcommand{\catC}{\category{C}}
\newcommand{\catCrho}{\catC_{\rho}}
\newcommand{\bcatCrho}{\boldsymbol{\catC}_{\myboldsymbol{\rho}}}
\newcommand{\catCU}{\algU\mbox{-\textsc{mod}}}
\newcommand{\catCD}{\DD(U_{\leq})\mbox{-\textsc{mod}}}

\newcommand{\functor}{\mathcal}
\newcommand{\functorDB}{\functor{F}_{\!\text{\textsc{db}}}}
\newcommand{\bfunctorDB}{\myboldsymbol{\functor{F}}_{\!\text{\textsc{db}}}}
\newcommand{\functorB}{\functor{F}_{\!\text{\textsc{b}}}}
\newcommand{\bfunctorB}{\myboldsymbol{\functor{F}}_{\!\text{\textsc{b}}}}
\newcommand{\functorT}{\functor{F}_{\text{\textsc{t}}}}
\newcommand{\functorDD}{\functor{F}_{\text{\textsc{dd}}}}

\newcommand{\leftact}{\kern1pt{\rightharpoonup}\kern1pt}
\newcommand{\rightact}{\kern1pt{\leftharpoonup}\kern1pt}

\newcommand{\cY}{\mathscr{Y}}
\newcommand{\cZ}{\mathscr{Z}}

\newcommand{\cR}{\mathscr{R}}

\newcommand{\eval}[2]{\langle#1,\,#2\rangle\,}

\newcommand{\id}{\mathrm{id}}

\newcommand{\bdelta}{\underline{\delta}}
\newcommand{\bDelta}{\underline{\Delta}}

\newcommand{\YDname}{\mathscr{Y\kern-1ptD}}

\newcommand{\BByd}{{}\mbox{\small${}^{\Nich(X)}_{\Nich(X)}$}\YDname}
\newcommand{\bBByd}{{}\mbox{\small${}^{\bNich\myboldsymbol{(X)}}_{\bNich\myboldsymbol{(X)}}$}\boldsymbol{\YDname}}
\newcommand{\HHyd}{{}\mbox{\small${}^{H}_{H}$}\YDname}

\newcommand{\BBHHYD}{{}\mbox{\small${}^{\Nich(X)\Smash H}_{\Nich(X)\Smash H}$}\YDname}

\newcommand{\RRYD}{{}\mbox{\small${}^{\cR}_{\cR}$}\YDname}
\newcommand{\RRyd}{{}{}^{\cR}_{\cR}\YDname}
\newcommand{\RHRHYD}{{}\mbox{\small${}^{\cR\Smash H}_{\cR\Smash
      H}$}\YDname} \newcommand{\RHRHyd}{{}^{\cR\Smash H}_{\cR\Smash
    H}\YDname}

\newcommand{\hA}{\mbox{\large$\mathpzc{s}\kern-.8pt$}}

\newcommand{\Smash}{\mathbin{\hash}}

\newcommand{\up}[1]{^{{\scriptscriptstyle(#1)}}}
\newcommand{\bup}[1]{^{{\scriptscriptstyle\!\underline{\,#1_{\vphantom{.}}\,}\!}}}
\newcommand{\zup}[1]{^{{\scriptscriptstyle[#1]}}}

\newcommand{\zero}{_{_{(0)}}}
\newcommand{\mone}{_{_{(-1)}}}

\newcommand{\bzero}{_{_{\,\underline{\,0\,}}}}
\newcommand{\bmone}{_{_{\underline{\!{-1}\!}}}}
\newcommand{\zzero}{_{_{[0]}}}
\newcommand{\zmone}{_{_{[-1]}}}

\def\verma{\bullet\rule[-2pt]{0pt}{9pt}
  \ar@{-}@*{[|(2.5)]}[0,0]+<0pt,0pt>;[0,0]-<10pt,0pt>
  \ar@{-}@*{[|(2.5)]}[0,0]-<1pt,-1pt>;[0,0]+<7pt,-7pt>
  \ar@{}[0,0]+<0pt,0pt>;}

\def\vermaiv{\bullet\rule[-2pt]{0pt}{9pt}
  \ar@{-}@*{[|(4)]}[0,0]+<0pt,0pt>;[0,0]-<10pt,0pt>
  \ar@{-}@*{[|(4)]}[0,0]-<1pt,-1pt>;[0,0]+<7pt,-7pt>
  \ar@{}[0,0]+<0pt,0pt>;}

\def\Cverma{\circ\rule[-2pt]{0pt}{9pt}
  \ar@{-}@*{[|(2.5)]}[0,0]+<-3pt,0pt>;[0,0]-<10pt,0pt>
  \ar@{-}@*{[|(2.5)]}[0,0]+<2pt,-2pt>;[0,0]+<6pt,-7pt>
  \ar@{}[0,0]+<0pt,0pt>;}

\def\iverma{\bullet\rule[-2pt]{0pt}{9pt}
  \ar@{-}@*{[|(2.5)]}[]+<0pt,0pt>;[0,0]+<10pt,0pt>
  \ar@{-}@*{[|(2.5)]}[]+<1pt,1pt>;[0,0]+<-7pt,-7pt>
  \ar@{}[0,0]+<0pt,0pt>;}

\def\ivermaiv{\bullet\rule[-2pt]{0pt}{9pt}
  \ar@{-}@*{[|(4)]}[]+<0pt,0pt>;[0,0]+<10pt,0pt>
  \ar@{-}@*{[|(4)]}[]+<1pt,1pt>;[0,0]+<-7pt,-7pt>
  \ar@{}[0,0]+<0pt,0pt>;}

\def\Civerma{\circ\rule[-2pt]{0pt}{9pt}
  \ar@{-}@*{[|(2.5)]}[]+<3pt,0pt>;[0,0]+<10pt,0pt>
  \ar@{-}@*{[|(2.5)]}[]-<2pt,2pt>;[0,0]+<-6pt,-7pt>
  \ar@{}[0,0]+<0pt,0pt>;}

\def\nothing{*{\mbox{}}\ar@{};[0,0];}

\def\Bullet{*{\rule[-1pt]{0pt}{7pt}\bullet}\ar@{};[0,0];}
\def\Circ{*{\rule[-1pt]{0pt}{7pt}\circ}\ar@{};[0,0];}

\newcommand{\ket}[1]{\mathchoice{%
    {\bigl|\smash[t]{#1}\bigr\rangle}}{|{#1}\rangle}{|{#1}\rangle}{|{#1}\rangle}}
\newcommand{\charSL}[2]{\chi_{%
    {\phantom{h}\kern-3pt #2}}^{\phantom{y}\kern-3pt #1}}

\newcommand{\bref}[1]{\textbf{\ref{#1}}}
\newcommand{\leftii}{\mathbin{{\smalltriangleright}\kern1pt}}

\newcommand{\thisfont}{\textsc}
\newcommand{\ab}[1]{\thisfont{a\kern-1pt b}_{#1}}
\newcommand{\ba}[1]{\thisfont{b\kern-1pt a}_{#1}}
\newcommand{\aba}[1]{\thisfont{a\kern-1pt b\kern-1pt a}_{#1}}
\newcommand{\bab}[1]{\thisfont{b\kern-1pt a\kern-1pt b}_{#1}}

\newcommand{\F}[1]{\thisfont{f\kern-1pt}_{#1}}
\newcommand{\Fb}[1]{\thisfont{f\kern-1pt b\kern-.5pt}_{#1}}
\newcommand{\bF}[1]{\thisfont{b\kern-1pt f\kern-1pt}_{#1}}
\newcommand{\bX}[1]{\thisfont{x\kern-1pt b\kern-1pt}_{#1}}
\newcommand{\bFb}[1]{\thisfont{b\kern-1pt f\kern-1pt b\kern-.5pt}_{#1}}

\newcommand{\qint}[1]{{\textstyle[#1]}}
\newcommand{\qfac}[1]{{\textstyle[#1]!\,}}

\newcommand{\tensor}{\mathbin{\otimes}}

\newcommand{\Dynkin}[3]{\xymatrix@C28pt@1{\circ\ar@{-}[r]^(.55){#2}\ar@{}^{#1}[]&\circ\ar@{}^{#3}[]}}

\newcommand{\mfrac}[2]{\raisebox{.8pt}{\mbox{\small$\displaystyle\frac{#1}{#2}$}}}
\newcommand{\ffrac}[2]{\raisebox{.5pt}{\mbox{\footnotesize$\displaystyle\frac{#1}{#2}$}}}
\newcommand{\fffrac}[2]{\raisebox{.9pt}{\mbox{\scriptsize$\displaystyle\frac{#1}{#2}$}}}
\newcommand{\half}{%
  \mathchoice{\ffrac{1}{2}}{\frac{1}{2}}{\frac{1}{2}}{\frac{1}{2}}}

\newcommand{\q}{\mathfrak{q}}

\newcommand{\hSL}[1]{\widehat{s\ell}(#1)}

\newcommand{\oC}{\mathbb{C}}
\newcommand{\oZ}{\mathbb{Z}}

\swapnumbers
\newtheorem{Thm}[subsection]{Theorem}
\newtheorem{thm}[subsubsection]{Theorem}

\newtheorem{Prop}[subsection]{Proposition}
\newtheorem{prop}[subsubsection]{Proposition}

\theoremstyle{definition}

\newtheorem{dfn}[subsubsection]{Definition}
\newtheorem{rem}[subsubsection]{Remark}

\begin{document}

\title{Representations of $\overline{\mathscr{U}}_{\q} s\ell(2|1)$ at
  even roots of unity}\thanks{This paper was written as
  ``\textbf{Logarithmic $\hSL2$ CFT{} models from Nichols
    algebras. 2.}'' But because it has no direct bearing on
  logarithmic models (which is left for the future) and in the hope
  that the algebraic structures discussed here may be interesting in
  their own right, we do not mention logarithmic CFT{} in the title.}

\author{A.\,M.\;Semikhatov} \author{I.\,Yu.\;Tipunin}

\address{Lebedev Physics Institute\hfill\mbox{}\linebreak
  \texttt{asemikha@gmail.com}, \ \texttt{tipunin@pli.ru}}

\begin{abstract}
  We\, construct\, all\, projective\, modules\, of\, the\,
  restricted\, quantum\, group $\overline{\mathscr{U}}_{\q}
  s\ell(2|1)$ at an even, $2p$th, root of unity.  This
  $64p^4$-dimensional Hopf algebra is a common double bosonization,
  $\Nich(X^*)\tensor\Nich(X)\tensor H$, of two rank-2 Nichols algebras
  $\Nich(X)$ with fermionic generator(s), with
  $H=\oZ_{2p}\tensor\oZ_{2p}$.  The category of
  $\overline{\mathscr{U}}_{\q} s\ell(2|1)$-modules is equivalent to
  the category of Yetter--Drinfeld $\Nich(X)$-modules in
  $\catCrho=\HHyd$, where coaction is defined by a universal
  $R$-matrix $\rho$.  As an application of the projective module
  construction, we find the associative algebra structure and the
  dimension, $5p^2-p+4$, of the $\overline{\mathscr{U}}_{\q}
  s\ell(2|1)$ center.
\end{abstract}

\maketitle
\thispagestyle{empty}

\section{\textbf{Introduction}}
We study the representation theory of a particular version of the
$s\ell(2|1)$ quantum group at even roots of unity.  For an integer
$p\geq 2$, our $64p^4$-dimensional $\overline{\mathscr{U}}_{\q}
s\ell(2|1)$ at $ \q = e^{i\pi/p}$ is a double bosonization of any of
the rank-$2$ Nichols algebras defined by the braiding matrices
\begin{equation}\label{qij-first}
  Q_{\asymm} = 
  \begin{pmatrix}
    -1 & \q^{-1} \\
    \q^{-1} & \q^2
  \end{pmatrix}
  \qquad\text{and}
  \qquad
  Q_{\symm} = 
  \begin{pmatrix}
    -1 & -\q \\
    -\q & -1
  \end{pmatrix},\qquad \q = e^{\frac{i\pi}{p}}.
\end{equation}

The interest in these Nichols algebras is motivated by the fact that
they centralize extended chiral algebras (vertex-operator algebras) of
logarithmic models of two-dimen\-sional conformal field theory with
$\hSL2_k$ symmetry, at the level $k=\frac{1}{p} - 2$.\footnote{We
  refer to \cite{[Gurarie],[Sa],[Kausch],[GK+]} for the origin of
  logarithmic conformal field theory (LCFT) and to \cite{[FHST]} for
  the idea to define LCFT models as centralizers of screening
  operators.  The screenings themselves turn out to generate a Nichols
  algebra~\cite{[STbr]}, hence the importance of Nichols algebras in
  this context.}  We recall that Nichols algebras are graded braided
Hopf algebras ``universally'' associated with a braided vector
space~$X$
(see~\cite{[Nich],[Wor],[Lu-intro],[Rosso-inv],[AG],[AS-pointed],
  [AS-onthe],[Heck-Weyl],[Heck-class],[AHS],[ARS],
  [Ag-0804-standard],[Ag-1008-presentation]} and the (co)references
therein).  The ``centralizer'' relation between vertex-operator
algebras and (braided or nonbraided) Hopf algebras is a key ingredient
underlying the Kazhdan--Lusztig correspondence~\cite{[KL],[Fink]} and
the quest for nice relations between suitable categories defined in
terms of Hopf algebras and representation categories of extended
symmetry algebras of the corresponding CFT{} models.

A ``logarithmic'' Kazhdan--Lusztig correspondence can be particularly
nice (an equivalence of categories)~\cite{[FGST2],[NT]} (also
see~\cite{[AM-lattice]} and, among more recent papers, \cite{[TW1]}).
However, the case studied in~\cite{[FGST2],[NT]} is nearly trivial in
the language of Nichols algebras, the case of rank
one~\cite{[b-fusion]}; true, the braided world can be quite rich, and
already the corresponding nonbraided Hopf algebra is $\Ures$, a not
altogether trivial quantum $s\ell(2)$ at an even ($2p$th) root of
unity.  It is this $\Ures$ that features in the categorial equivalence
with the $(1,p)$ LCFT models~\cite{[FGST],[FGST2]}; the correspondence
between Hopf-algebraic and LCFT realms, moreover, extends to modular
group representations: those on the quantum group center and on the
torus amplitudes of the logarithmic model turn out to
coincide.\footnote{A feature that remarkably survives in the case
  where the representation category on the LCFT side is not that
  nice~\cite{[FGST3],[FGST-q],[GRW],[RGW]}.}

The ideas of the logarithmic Kazhdan--Lusztig correspondence need to
be extended to higher-rank Nichols algebras; this would show, among
other things, how much of what we know in the rank-1 case is
``accidental,'' and which features are indeed generic.
Moving to higher-rank Nichols algebras was initiated
in~\cite{[c-charge]} and, with precisely the two Nichols algebras
defined by braiding matrices~\eqref{qij-first},
in~\cite{[nich-sl2-1]}.

First and foremost, with the structural theory of finite-dimensional
Nichols algebras with diagonal braiding
completed~\cite{[Heck-class],[Ag-0804-standard],[Ag-1008-presentation]},
knowledge about their appropriate representation categories (of
Yetter--Drinfeld modules) is highly desirable.\footnote{Once again on
  the subject of LCFT models, we note that their \textit{fusion} has a
  good chance to be described just by the ring structure of a suitable
  category of finite-dimensional modules, finding which, difficult
  though it may be, is ``infinitely'' easier than deriving the fusion
  algebra directly.  The examples where the fusion algebra known or
  reasonably conjectured in other approaches coincides with the one
  taken over from the Hopf-algebra side are quite encouraging (see,
  e.g., \cite{[FHST],[FGST],[BFGT],[BGT]} and the references
  therein).}  In this paper, we address the representation theory of
the common double bosonization of the two chosen Nichols algebras,
which is the $64 p^4$-dimensional $\overline{\mathscr{U}}_{\q}
s\ell(2|1)$.  We first show that the category of
$\overline{\mathscr{U}}_{\q} s\ell(2|1)$ modules is equivalent to the
category $\BByd_{\catCrho}$ of Yetter--Drinfeld $\Nich(X)$ modules in
$\catCrho=\HHyd$, where $\Nich(X)$ is the corresponding Nichols
algebra and $H=\oZ_{2p}\tensor\oZ_{2p}$ is the Cartan subalgebra of
$\overline{\mathscr{U}}_{\q} s\ell(2|1)$, with its coaction on $\HHyd$
defined by the universal $R$-matrix $\rho\in H\tensor H$.  Our main
result is then the construction of projective
$\overline{\mathscr{U}}_{\q} s\ell(2|1)$ modules.  The ensuing picture
is rather involved, and may be a good illustration of the intricacies
occurring at roots of unity.

An essential part of the structure of projective modules can be
conveniently expressed in terms of directed graphs whose vertices are
simple subquotients and the edges are associated with elements of
$\Exti$ groups, weighted with some coefficients (finding which is a
major part of the existence proof for a given projective module).  The
paper therefore contains a number of pictures showing graphs of
projective modules.

As an application of using the data presented in these graphs, we find
the associative algebra structure of the center of
$\overline{\mathscr{U}}_{\q} s\ell(2|1)$:
\begin{equation*}
  \Zcenter = Z_{\mathrm{at}}\oplus\bigoplus_{j=1}^{(p-1)(p-2)}Z_{\mathrm{t}}^j
  \oplus\bigoplus_{j=1}^{4(p-1)}Z_{\mathrm{st}}^j
\end{equation*}
where $Z_{\mathrm{at}}$, $Z_{\mathrm{t}}^j$, and $Z_{\mathrm{st}}^j$
are commutative associative unital algebras; all $Z_{\mathrm{st}}^j$
are $1$-dimensional, each $Z_{\mathrm{t}}^j$ is $5$-dimensional
(contains 4 nilpotents in addition to the unit), and $Z_{\mathrm{at}}$
is $10p-2$ dimensional ($10p-3$ nilpotents in addition to the unit).

In Sec.~\ref{algUa-def}, we define our Hopf algebra
$\algU\equiv\overline{\mathscr{U}}_{\q} s\ell(2|1)$ and discuss some
of its simple properties (including Casimir elements and the universal
$R$-matrix).  In Sec.~\ref{sec:equivalence}, we prove that the
category of its modules is equivalent to $\BByd$.  \ In
Sec.~\ref{sec:simples}, we classify its simple modules, describing
them quite explicitly.  We continue in Sec.~\ref{sec:simples} with
listing the $\Exti$ spaces for the simple modules.  We then construct
the projective $\algU$ modules in Sec.~\ref{sec:proj}. \ As an
application of our treatment of projective modules, we find the center
of $\algU$ in Sec.~\ref{sec:center}.  In an attempt to improve the
readability of this inevitably technical paper, we isolate some
computational details in the appendices.

We use $q$-integers and factorials
defined as
\begin{equation*}
  \qint{n}=\ffrac{q^{n}-q^{-n}}{q-q^{-1}},
  \qquad
  \qfac{n}=\qint{1}\dots\qint{n},
\end{equation*}
all of which are assumed specialized to $q=\q$ in~\eqref{qij-first}.

All (co)modules in this paper are finite-dimensional.

\section{\textbf{The Hopf algebra $\balgU$}}\label{algUa-def}
\subsection{}
The notation $\algU$ for our $64p^4$-dimensional quantum group
$\overline{\mathscr{U}}_{\q} s\ell(2|1)$ is a legacy of the
Nichols-algebra setting, where $X$ in $\Nich(X)$ is a two-dimensional
braided vector space---in our case, specifically, the one with
diagonal braiding defined by any of the two braiding matrices
in~\eqref{qij-first}.  The Hopf algebra $\algU$ was derived from each
of these two Nichols algebras in~\cite{[nich-sl2-1]}.

\subsubsection{The Hopf algebra
  $\balgU\equiv\boldsymbol{\overline{\mathscr{U}}_{\q}}
  \myboldsymbol{s\ell(2|1)}$}\label{U(X)-gens}
Our $\algU$ is the algebra on $ B $, $F$, $k$, $K$, $ C $, $ E $ with
the relations listed in~\eqref{Ures}--\eqref{the-other}. \ We first
identify an important Hopf subalgebra in~$\algU$, the restricted
quantum group $\Ures$.  It is generated by $ E $, $K$, and $F$, with
the relations
\begin{equation}\label{Ures}
  \begin{gathered}
    K F=\q^{-2} F K,\quad E F - F E
    =\mfrac{K-K^{-1}}{\q-\q^{-1}},\quad K E =\q^{2} E K,
    \\
    F^p=0,\qquad E^p=0,\qquad K^{2p}=\one.
  \end{gathered}
\end{equation}
Next, $\Ures$ and $k$ generate a larger algebra $\Uresk$ with further
relations
\begin{gather}\label{Ustar}
  k F = \q F k, \quad k E = \q^{-1} E k, \quad k^{2p}=\one,\quad kK =
  Kk.
\end{gather}
The other relations in $\algU$\,---\,those involving fermions $B$ and
$C$\,---\,are
\begin{equation}\label{the-other}
  \begin{gathered}
    K B = \q B K,\quad k B =- B k, \quad K C = \q^{-1} C K,\quad k C
    =- C k,
    \\[-2pt]
    B^2 = 0, \quad B C - C B =\mfrac{k-k^{-1}}{\q-\q^{-1}}, \quad C^2
    = 0,
    \\
    F C - C F =0, \qquad B E - E B =0,
    \\
    F F B - \qint{2} F B F + B F F = 0, \quad E E C - \qint{2} E C E +
    C E E = 0.
  \end{gathered}
\end{equation}

The Hopf-algebra structure of $\algU$ is furnished by the coproduct,
antipode, and counit given by
\begin{align*}
  &\begin{alignedat}{2} \zDelta( F )&= F \otimes1 + K^{-1}\otimes F,\qquad &\Delta( E )&= E \tensor K+\one\tensor E,
    \\
    \zDelta( B )&= B \otimes1 + k^{-1}\otimes B, &\Delta( C )&= C
    \tensor k+\one\tensor C,
  \end{alignedat}
  \\
  &\begin{alignedat}{2} \zS( B )&=-k B,\quad \zS( F )=-K F,
    \quad&S(C)&=- C k^{-1}, \quad S( E )=- E K^{-1},
  \end{alignedat}
  \\
  &\begin{alignedat}{2} \epsilon( B )&=0,\quad\epsilon( F )=0, &\quad
    \epsilon( C )&=0,\quad \epsilon( E )=0,
  \end{alignedat}
\end{align*}
with $k$ and $K$ group-like.

By the $\algU$ generators in what follows, we always mean $ B $, $ F
$, $ C $, $ E $, $k$, $K$.

\subsubsection{The second Hopf algebra structure~\cite{[nich-sl2-1]}}
\label{2ndHopf}
The Hopf algebra $\algU$ admits a nontrivial ``twist''---an invertible
normalized 2-cocycle (see~\bref{app:MM})
\begin{equation*}
  \Phi = 1\tensor 1  + (\q - \q^{-1}) B k\tensor C k^{-1}
  \in\algU\tensor\algU
\end{equation*}
twisting by which gives rise to the second coalgebra structure
\begin{equation*}
  \tDelta(x)=\Phi^{-1}\Delta(x)\Phi.
\end{equation*}
(Two coalgebra structures naturally come from the two underlying
Nichols algebras.)

For the $\algU$ generators chosen above, the coproducts are
\begin{align*}
  \tDelta F &=F\tensor 1 + K^{-1}\tensor F + (\q - \q^{-1}) F B
  k\tensor C k^{-1} + (1 - \q^2) B F k\tensor C k^{-1},
  \\
  \tDelta B&= B\tensor k^{-2} + k^{-1}\tensor B,
  \\
  \tDelta E&=1\tensor E + E\tensor K + (\q - \q^{-1}) B k\tensor E C
  k^{-1} + (1 - \q^2) B k\tensor C E k^{-1},
  \\
  \tDelta C&= C\tensor k + k^{2}\tensor C.
\end{align*}
We note that the new coproduct has the following skew-primitive
elements:
\begin{align*}
  \tDelta(B k) &= 1\tensor B k + B k\tensor k^{-1},
  \\
  \tDelta(F B - \q B F) &= (F B - \q B F)\tensor 1 + K^{-1}
  k^{-1}\tensor (F B - \q B F),
  \\
  \tDelta(C k^{-1}) &=C k^{-1}\tensor 1 + k\tensor C k^{-1},
  \\
  \tDelta(E C - \q C E) &=1\tensor (E C - \q C E) + (E C - \q C E)\tensor K k.
\end{align*}

\subsection{PBW basis in $\balgU$} We have a linear-space isomorphism
\begin{equation}\label{UU-decomp}
  U_{>}\tensor U_{\leq} \to \algU,
\end{equation}
where $U_{>}$ is the subalgebra in $\algU$ generated by $E$ and $C$,
and $U_{\leq}$ is the subalgebra generated by $F$, $B$, $K$, and $k$.
A PBW basis in $U_{\leq}$ (which we refer to as \textit{the} PBW basis
in $U_{\leq}$) can be chosen as $F^{n} K^i k^{j}$, $F^{n}BF K^i
k^{j}$, $F^{n} B F K^i k^{j}$, $F^{n} B F B K^i k^{j}$, where $0\leq
n\leq p-1$ and $0\leq i,j\leq 2p-1$; the PBW basis in $U_{>}$ can be
chosen similarly, as $E^{n}$, $E^{n + 1} C$, $E^{n} C E$, $E^{n} C E
C$, where $0\leq n\leq p-1$.  The PBW basis in $\algU$ is then given
by the product of these two PBW bases.

\subsection{Casimir elements}\label{sec:casimirs}
There are two elements of a particular form in the center of $\algU$,
which we call Casimir elements:
\begin{align*}
  C_1&=
  E F K k^{2} + C B K^{2} k^{3} - \q (E C - \q C E) (F B - \q^{-1} B F) K k^{3}
  \\[-2pt]
  &\quad{} + \ffrac{\q}{(\q - \q^{-1})^2}\, k^{2} +
  \ffrac{\q^{-1}}{(\q - \q^{-1})^2}\, K^{2} k^{2} -
  \ffrac{\q^{-1}}{(\q - \q^{-1})^2}\, K^{2} k^{4},
  \\[4pt]
  C_2&=
  E F K^{-1} k^{-2} + C B K^{-2} k^{-3} - \q^{-1} (E C - \q^{-1} C E)
  (F B - \q B F) K^{-1} k^{-3}
  \\[-2pt]
  &\quad{} + \ffrac{\q^{-1}}{(\q - \q^{-1})^2}\, k^{-2} +
  \ffrac{\q}{(\q - \q^{-1})^2}\, K^{-2} k^{-2} - \ffrac{\q}{(\q -
    \q^{-1})^2}\,K^{-2} k^{-4}.
\end{align*}
That they are central is verified directly.

\begin{prop}
  Each of the Casimir elements satisfies a minimal polynomial relation
  of degree $p^2 - 2 p + 4$:
  \begin{equation*}
    \Bigl(C_1 - \ffrac{\q}{(\q - \q^{-1})^2}\Bigr)^{\!3}\;
    \prod_{s=1}^{p - 1}\prod_{r=1}^{s - 1}
    \Bigl(C_1 - \ffrac{
      \q^{1 - 4 r} (\q^{2 r} - \q^{2 s} + \q^{2 r + 2 s})}{(\q - \q^{-1})^2}
    \Bigl)^{\!2}\;
    \prod_{r=1}^{p - 1}
    \Bigl(C_1 - \ffrac{\q^{1 - 4 r} (2 \q^{2 r} - 1)}{(\q - \q^{-1})^2}\Bigr)
    = 0,
  \end{equation*}
  and the minimal relation for $C_2$ is obtained from here by
  replacing $\q\to\q^{-1}$.
\end{prop}
This is an immediate corollary of our construction of projective
$\algU$-modules in Sec.~\ref{sec:proj}.
We simply take the product of $(C_i - \lambda)$ factors over all
different eigenvalues $\lambda$ of a given $C_i$ on the (linkage
classes of) projective modules, with each factor raised to the power
given by the corresponding Jordan cell size (three for the atypical
linkage class, two for each of the typical linkage classes, and one
for each Steinberg module$/$class; see~\bref{sec:link-cl}).

A somewhat more involved derivation shows that the two Casimir
elements satisfy a degree-$p$ polynomial relation,
\begin{equation*}
  \sum_{i=1}^{p} \ffrac{(-1)^{i}}{i}
  \mbox{\small$\displaystyle\binom{p - 1}{i - 1}$}
  (\q - \q^{-1})^{2 i} \bigl(\q^{-i} C_1^i -
  \q^{i} C_2^i\bigr) = 0.
\end{equation*}
The full list of (``mixed'') relations satisfied by $C_1$ and $C_2$ is
currently unknown.

\subsection{Quasitriangular and related structures}
\begin{thm}\label{thm:R}
  The algebra $\algU$ is quasitriangular, with the universal
  $R$-matrix given by
  \begin{align}
    R &=\rho\bar{R},\nonumber
    \\
    \label{the-rho}
    \rho &= \ffrac{1}{(2 p)^2} \sum_{i=0}^{2 p - 1}\sum_{j=0}^{2 p -
      1} \sum_{m=0}^{2 p - 1}\sum_{n=0}^{2 p - 1} (-1)^{j n} \q^{-2 i
      m + j m + i n} K^{i} k^{j}\tensor K^{m} k^{n},
    \\
    \label{barR-factored}
    \bar{R} &=\sum_{a=0}^{p - 1} \ffrac{\q^{\half a (a - 1)} (\q -
      \q^{-1})^a}{\qfac{a}} E^{a}\tensor F^{a}\; \Bigl(1\tensor 1 -
    (\q - \q^{-1}) C\tensor B\Bigr)
    \\*
    &\quad\times \Bigl(1\tensor 1 + \q (\q - \q^{-1}) \Bar{C}\tensor
    \Bar{B}\Bigr)
    \Bigl(1\tensor 1 + (\q - \q^{-1})^3 \Bar{\Bar{C}}\tensor
    \Bar{\Bar{B}}\Bigr), \notag
  \end{align}
  where
  \begin{alignat*}{2}
    \Bar{C} &= E C - \q^{-1} CE,&\qquad \Bar{B} &= F B - \q^{-1} B F,\\
    \Bar{\Bar{C}} &= C \Bar{C} = C E C,& \Bar{\Bar{B}} &= B \Bar{B} =
    B F B.
  \end{alignat*}
\end{thm}
This must be possible to extract from~\cite{[AY]}; we give the proof
by direct calculation in~\bref{app:R}.

The category of $\algU$ modules is therefore braided, with the
braiding $\repY\tensor\repZ\to\repZ\tensor\repY$ given by $y\tensor
z\mapsto \Rii z\tensor\Ri y$ (where we standardly write
$R=\Ri\tensor\Rii$).

We let $\catCU$ denote the braided monoidal category of
$\algU$-modules.

\subsubsection{}\label{R-properties}We recall the general properties
of the universal $R$-matrix for further reference:
$R\Delta(x)=\Delta^{\text{op}}(x)R$, i.e., $\Ri x'\tensor\Rii
x''=x''\Ri\tensor x'\Rii$ for all $x$ in the algebra, and
\begin{align*}
  {\Ri}'\tensor{\Ri}''\tensor\Rii&= \Ri\tensor\Pone\tensor\Rii\Pii,
  \\
  \Ri\tensor{\Rii}'\tensor{\Rii}''&= \Ri\Pone\tensor\Pii\tensor\Rii
\end{align*}
(with $R=\Ri\tensor\Rii=\Pone\tensor\Pii$, and so on).

\subsubsection{Drinfeld map}\label{sec:M}
With the above $R$-matrix, we introduce the so-called
M[onodromy] ``matrix''
\begin{equation*}
  M=R_{21}R\in\algU\tensor\algU.
\end{equation*}
It is known to give rise to the Drinfeld map from the space $\Ch$ of
$q$-characters to the center $\Zcenter$ of $\algU$:
\begin{equation}\label{dr-map}
  \drmap:\Ch\to\Zcenter:
  \alpha\mapsto \alpha\tensor\id(M)
  =\alpha(M_1)M_2,
\end{equation}
where we standardly write $M=M_1\tensor M_2$.  We recall that $\Ch$ is
the space of elements of $\algU^*$ that are invariant under the left
coadjoint action of $\algU$.

\begin{prop}\label{prop:M-matrix}
  For the $R$-matrix in~\bref{thm:R}, the M-matrix can be written as
  \begin{equation}\label{M-matrix}
    M=\bar{M}\bar{\rho}\bar{R},
  \end{equation}
  where
  \begin{align*}
    \bar{M}&= \sum_{a=0}^{p - 1} \ffrac{\q^{-\frac{3}{2} a^2 - \half
        a} (\q - \q^{-1})^a\!}{\qfac{a}}\, F^{a} K^a\tensor E^{a}
    K^{-a} \, \bigl(1\tensor 1 + (\q - \q^{-1}) B k\tensor C
    k^{-1}\bigr)
    \\
    &\quad\times \bigl(1\tensor 1 - (\q^2 - 1) \Bar{B} K k\tensor
    \Bar{C} K^{-1} k^{-1}\bigr) \bigl(1\tensor 1 + \q^{2} (\q -
    \q^{-1})^3 \Bar{\Bar{B}} K k^{2}\tensor \Bar{\Bar{C}} K^{-1}
    k^{-2}\bigr)
    \\
    \bar{\rho}&= \ffrac{1}{p^2} \sum_{i=0}^{p - 1}\sum_{j=0}^{p -
      1}\sum_{i'=0}^{p - 1}\sum_{j'=0}^{p - 1} \q^{2 i' j + 2 i j' - 4
      i i'} K^{2 i} k^{2 j}\tensor K^{2 i'} k^{2 j'}.
  \end{align*}
\end{prop}
We show this in~\bref{M-calc}.

\subsubsection{}
The quasitriangular structure corresponding to $\tDelta$
in~\bref{2ndHopf} is given by
\begin{equation*}
  \widetilde{R}=\Phi_{21}^{-1} R \Phi,
\end{equation*}
as readily follows from $R\Delta(x)=\Delta^{\text{op}}(x)R$.  Hence,
the corresponding M-matrix is
\begin{equation*}
  \widetilde{M}=\Phi^{-1} M \Phi.
\end{equation*}
We want to compare the two Drinfeld maps $\drmap:\Ch\to\Zcenter$ and
$\widetilde{\drmap}:\widetilde{\Ch}\to\Zcenter$, where
$\widetilde{\Ch}$ is the space of elements that are invariant under
the ``tilded'' coadjoint action.

The two Drinfeld maps turn out to coincide in the sense that the
diagram
\begin{equation}\label{MM-diag}
  \xymatrix@C10pt{
    \Ch\ar_(.5){\drmap}[dr]\ar[rr]&&\widetilde{\Ch}\ar^(.5){\widetilde{\drmap}}[dl]
    \\
    &\Zcenter&
  }
\end{equation}
is commutative, where the horizontal arrow is a linear isomorphism
$\beta\mapsto(\beta\rightact\xi)=\beta(\xi?)$, where $\xi=S(U^{-1}) U$
and $U=S(\Phi_1)\Phi_2$.  We prove this in the general case
in~\bref{thm:MM}.

With the above $\Phi$, it readily follows that
\begin{equation}\label{Uixi}
  U = 1 - (\q - \q^{-1}) B C k^{-1}\quad\text{and}\quad
  \xi=k^{-2}.
\end{equation}

\section{$\algU$ modules and Yetter--Drinfeld $\Nich(X)$
  modules}\label{sec:equivalence}
We show that the category of $\algU$ modules is equivalent to a
category of Yetter--Drinfeld modules of the corresponding Nichols
algebra.  The exact statement is in~\bref{thm:equivalence} below.  We
begin with briefly recalling the relation between $\algU$ and Nichols
algebras.

\subsection{$\balgU$ as a double bosonization}
The Hopf algebra $\algU$ is a double bosonization \cite{[nich-sl2-1]}
(also see~\cite{[HY],[ARS],[HS-double]}) of the Nichols algebra
$\Nich(X)$ of a two-dimensional braided vector space $X$ with a basis
$(B,F)$ such that the corresponding braiding matrix is $Q_{\asymm}$
in~\eqref{qij-first} .\footnote{We select the first braiding matrix
  in~\eqref{qij-first} because we mainly describe $\algU$ with the
  generators chosen as $B$, $F$, $C$, $E$, $k$, $K$, not those that
  are skew-primitive with respect to the second coproduct
  in~\bref{2ndHopf}.} \ This $\Nich(X)$ is the
quotient~\cite{[Ag-0804-standard]}
\begin{equation}\label{eq:asymm-quotient}
  \Nich(X)=T(X)/ \bigl(\Fi\Fii^2 - \qint{2}\Fii\Fi\Fii + \Fii^2 \Fi,
  \ \Fi^2,\ \Fii^{p}\bigr)
\end{equation}
for $p\geq 3$ (and $\Nich(X)=T(X)/ \bigl(\Fi^2,\ \Fi \Fii \Fi \Fii -
\Fii \Fi \Fii \Fi,\ \Fii^{2}\bigr)$ if $p=2$), with $\dim\Nich(X)=4p$.
Constructing its double bosonization requires choosing a cocommutative
ordinary Hopf algebra $H$ such that $\Nich(X)\in\HHyd$.  Specifically,
we take $H=\oZ_{2p}\tensor\oZ_{2p}$, with the generators $ g_1\equiv
k$ and $g_2\equiv K$ acting on the $\Nich(X)$ generators $F_1\equiv
\Fi$ and $F_2\equiv \Fii$ as $g_i\dotact F_j = q_{i,j} F_j$, where
$q_{i,j}$ are the entries of the braiding matrix; $H$ coacts as
$F_i\mapsto g_i\tensor F_i$. \ The double bosonization
\begin{equation}\label{algU}
  \algU = \Nich(X^*)\mathbin{{\Smash}'}\Nich(X)\Smash H
\end{equation}
(where the right-hand side is a tensor product in $\HHyd$) then
contains the Hopf subalgebra $\Nich(X)\Smash H$ given by the standard,
``single,'' bosonization~\cite{[Radford-bos]}, and similarly for
$\Nich(X^*)\mathbin{{\Smash}'} H$, but with the $H$ action and
coaction changed by composing each with the antipode (hence the
prime); the cross relations are $\dF_i F_j - q_{j,i} F_j \dF_i =
\delta_{i,j}(1 - (g_j)^2)$ in terms of the dual basis $\dF_i$ in
$X^*$.  After a suitable change of basis, this yields
our~$\algU$~\cite{[nich-sl2-1]}.

Naturally, $H$ is identified with the Hopf subalgebra in $\algU$
generated by $k$ and $K$. \ Moreover, we refine~\eqref{UU-decomp} by
introducing the linear space isomorphism
\begin{equation}\label{UUU-decomp}
  \Uminus \tensor H\tensor U_{>} \to \algU,
\end{equation}
where $\Uminus$ is the subalgebra in $\algU$ generated by $F$ and $B$.
Each of the algebras $\Uminus$, $H$, and $U_{>}$ is a Hopf algebra.

\begin{dfn}\label{dfn:catC}
  The category $\catCrho$ is the category $\HHyd$ of
  Yetter--Drinfeld $H$-modules with the coaction
  \begin{equation}\label{eq:rho-action}
    \delta:z\mapsto z\mone\tensor z\zero
    =\rhoii\tensor\rhoi\dotact z,
  \end{equation}
  where $\rho$ is given by~\eqref{the-rho}.
\end{dfn}
Evidently, the braiding in $\catCrho$ is given by
\begin{equation}\label{eq:C-braiding}
  u\tensor v\mapsto u\mone\dotact v\tensor u\zero=\rhoii\dotact
  v\tensor\rhoi\dotact u.
\end{equation}

\begin{Thm}\label{thm:equivalence}
  The category $\catCU$ 
  is equivalent as a braided monoidal category to the category
  $\BByd_{\catCrho}$ of Yetter--Drinfeld $\Nich(X)$-modules in
  $\catCrho$.
\end{Thm}
 
We prove this in the rest of this section, in several steps.  We first
construct a transmutation functor
$\functorT:\catCU\to\BByd_{\catCrho}$
(cf.~\cite{[Mj-braided],[Mj-trans-rank]}).  Its inverse, a double
bosonization functor $\functorDB:\BByd_{\catCrho}\to\catCU$, is then
constructed as a composition of the bosonization functor $\functorB$
and the functor $\functorDD$ sending Yetter--Drinfeld modules to
modules of the Drinfeld double.  Both $\functorT$ and $\functorDB$ are
vector-space-preserving functors.  We then verify that their
compositions are indeed equivalent to the corresponding identity
functors.

\begin{Prop}\label{thm:fun-new}
  There is a vector-space-preserving braided monoidal functor
  $\functorT$ sending each $\algU$-module into a Yetter--Drinfeld
  $\Nich(X)$-module in $\catCrho$.
\end{Prop}
We prove this in~\bref{fun-new-algebra}--\bref{fun-new-monoidal}.
Whenever possible (and this is most often possible), we do not assume
the ``ground'' Hopf algebra $H$ to be cocommutative (nor is it assumed
commutative).  This in fact clarifies the calculations, highlighting
the ``true'' reasons why certain identities hold.  We frequently use
the universal $R$-matrix properties in~\bref{R-properties}, for
$\rho\in H\tensor H$ in particular.

\subsubsection{$\myboldsymbol{\Uminus}$ is an algebra in
  $\bcatCrho$}\label{fun-new-algebra}
We make $\Uminus$ into an algebra object in $\HHyd$ by defining the
left (adjoint) action and left coaction as\footnote{We do not indicate
  the embedding $\iota:H\to\algU$ explicitly, and write $h' x S(h'')$
  for what should be $\iota(h)\leftii x$, etc.; accordingly, we do not
  use a special symbol for the antipode in $H$, and so on.}
\begin{align*}
  h\leftii x &= h' x S(h''),
  \\
  \delta:x&\mapsto x\mone\tensor\ x\zero = \rhoii\tensor \rhoi\leftii
  x.
\end{align*}
The Yetter--Drinfeld axiom $(h\leftii x)\mone\tensor(h\leftii x)\zero
=h' x\mone S(h''')\tensor h''\leftii x\zero$ is then immediate to
verify.
That the product $\Uminus\tensor\Uminus\to\Uminus$ is a
$\catCrho$-morphism is also evident.

\subsubsection{$\myboldsymbol{\Uminus}$ is a Hopf algebra in
  $\bcatCrho$}\label{U-in-HH}
We next define the coaction
\begin{equation}\label{bDelta}
  \bDelta:\Uminus\to\Uminus\tensor\Uminus:
  x\mapsto x\bup{1}\tensor x\bup{2}=x'S(\rhoii)\tensor \rhoi\leftii x''.
\end{equation}
A priori, the right-hand side is an element not of
$\Uminus\tensor\Uminus$ but of $U_{\leq}\tensor\Uminus$\,. \ However,
it is immediately verified for the generators $A=F,\,B$ that
\begin{equation*}
  \bDelta(A) = A\tensor 1 + 1\tensor A \in \Uminus\tensor\Uminus,
\end{equation*}
and we then conclude that
$\bDelta(\Uminus)\subset\Uminus\tensor\Uminus$ from the main axiom of
braided Hopf algebras,
\begin{equation}\label{bDelta-mult}
  \bDelta(x y) = x\bup{1}(x\bup{2}{}\mone\leftii y\bup{1})
  \tensor x\bup{2}{}\zero y\bup{2}.
\end{equation}
This last identity is a relatively standard
statement~\cite{[Mj-braided],[Mj-trans-rank]}, but we prove it here
for completeness:
\begin{align*}
  \text{r.h.s.\ of~\eqref{bDelta-mult}}
  &= x' S(\rhoii)(\sigmaii\leftii y' S(\tauii))
  \tensor (\sigmai\rhoi\leftii x'')(\taui\leftii y'')
  \\[-2pt]
  &= x' S(\rhoii')(\rhoii''\leftii y' S(\tauii))
  \tensor (\rhoi\leftii x'')(\taui\leftii y'')
  \\[-2pt]
  &= x' y' S(\rhoii\tauii)
  \tensor (\rhoi\leftii x'')(\taui\leftii y'')
  \\[-2pt]
  &= x' y' S(\rhoii)
  \tensor (\rhoi'\leftii x'')(\rhoi''\leftii y'')
  \\[-2pt]
  &= x' y' S(\rhoii)
  \tensor (\rhoi\leftii x'' y'')
  = \text{l.h.s.\ of~\eqref{bDelta-mult}}.
\end{align*}

Also, $\bDelta$ is an $\HHyd$ morphism.  For the $H$-coaction, we have
to prove that
\begin{equation*}
  x\mone\tensor x\zero{}\bup{1}\tensor x\zero{}\bup{2}
  =x\bup{1}{}\mone\,x\bup{2}{}\mone\tensor x\bup{1}{}\zero\tensor x\bup{2}{}\zero
\end{equation*}
where we calculate
\begin{align*}
  \text{r.h.s.}
  &= \tauii\sigmaii\tensor(\taui\leftii x' S(\rhoii))
  \tensor(\sigmai\rhoi\leftii x'')
  \\[-2pt]
  &= \tauii\tensor(\taui'\leftii x' S(\rhoii))
  \tensor(\taui''\rhoi\leftii x'')
  \\[-2pt]
  &= \tauii\tensor(\taui' x' S(\taui''\rhoii))
  \tensor(\taui'''\rhoi\leftii x'')
  \\[-2pt]
  &= \tauii\tensor(\taui' x' S(\rhoii \taui'''))
  \tensor(\rhoi\taui''\leftii x''),  
\end{align*}
and
\begin{align*}
  \text{l.h.s.}
  &= \tauii\tensor (\taui' x S(\taui''))'S(\rhoii)\tensor
  (\rhoi\leftii(\taui' x S(\taui''))'')
  \\[-2pt]
  &= \tauii\tensor \taui' x' S(\taui'''') S(\rhoii)\tensor
  (\rhoi\leftii\taui'' x'' S(\taui''')),
\end{align*}
which is the same.  For the $H$ action, we readily see that
$h\leftii\bDelta(x)=\bDelta(h\leftii x)$
using the fact that $H$ is cocommutative.
We thus conclude that $\Uminus$ is a Hopf algebra object
in~$\catCrho$.

\begin{rem}
  As in~\cite{[AS-pointed]}, we can characterize $\Uminus$ inside
  $U_{\leq}$ as
  \begin{equation*}
    \Uminus=\{x\in U_{\leq}\mid x'\tensor\pi(x'')=x\tensor 1\}
  \end{equation*}
  with a projection (Hopf-algebra map) $\pi:U_{\leq}\to H$.  Then
  \begin{equation*}
    \pi(x')\tensor x'' = \rhoii\tensor\rhoi\leftii x.
  \end{equation*}
\end{rem}

\subsubsection{$\myboldsymbol{\Uminus={}}\bNich\myboldsymbol{(X)}$}
In $\catCrho$, the braiding matrix for the generators $B$ and $F$ is
exactly $Q_{\asymm}$ in~\eqref{qij-first}.  Hence, $\Uminus$ endowed
with the coproduct $\bDelta$ (and with the antipode $x\mapsto \rhoii
S(\rhoi\leftii x)$, which we do not discuss in any further detail for
brevity) is the Nichols algebra $\Nich(X)$ of the two-dimensional
braided linear space $X$ with the braiding matrix~$Q_{\asymm}$.

We write $\Nich(X)$ instead of $\Uminus$ in what follows.

\subsubsection{$\balgU$ modules are objects of $\catCrho$}
Let $\repZ\in\catCU$, with a $\algU$ action $x\tensor z\mapsto
x\dotact z$.  By restriction, $\repZ$ is an $H$-module.  We make
$\repZ$ into an object in $\catCrho$ by endowing it with the
$H$-coaction~\eqref{eq:rho-action}.

\subsubsection{$\balgU$ modules are objects of
  $\bBByd$}\label{fun-new-modules}
On each $\algU$-module $\repZ$, which is a $\Nich(X)$-module by
restriction, we also define the $\Nich(X)$-coaction
\begin{align}\nonumber
  \bdelta:\repZ&\to\Nich(X)\tensor\repZ 
  \\
  \label{bdelta}
  z&\mapsto z\bmone\tensor z\bzero= \Rii
  S(\rhoii)\tensor\rhoi\Ri\dotact z =\rhoii\leftii\bRii\tensor
  \rhoi\bRi\dotact z,
\end{align}
where $\bar{R}$ is defined in~\eqref{barR-factored}; the second form
shows that $\bdelta$ indeed maps to $\Nich(X)\tensor\repZ$.  The
coaction property $z\bmone\tensor z\bzero{}\bmone\tensor
z\bzero{}\bzero = z\bmone{}\bup{1}\tensor z\bmone{}\bup{2}\tensor
z\bzero$ is established by direct calculation.  We must also verify
that $\bdelta$ is an $\HHyd$ morphism.  We calculate, writing
$\rho=\rhoi\tensor\rhoii=\sigmai\tensor\sigmaii=\taui\tensor\tauii
=\dots$ (and assuming neither cocommutativity nor commutativity of
$H$):
\begin{align*}
  h(\bdelta z)
  &=(h'\leftii\Rii S(\rhoii))\tensor h''\rhoi\Ri\dotact z
  \\[-2pt]
  &=h'\Rii S(h''\rhoii)\tensor h'''\rhoi\Ri\dotact z
  \\[-2pt]
  &=h'\Rii S(\rhoii h''')\tensor \rhoi h''\Ri\dotact z
  \\[-2pt]
  &=\Rii h'' S(\rhoii h''')\tensor \rhoi \Ri h'\dotact z
  \\[-2pt]
  &=\Rii \tensor \rhoi \Ri h\dotact z = \bdelta(h\dotact z),
\end{align*}
and
\begin{align*}
  z\bmone{}\mone z\bzero{}\mone\tensor z\bmone{}\zero\tensor
  z\bzero{}\zero
  &= \rhoii\sigmaii\tensor(\rhoi\leftii \Rii S(\tauii))\tensor
  \sigmai\taui\Ri\dotact z
  \\[-2pt]
  &= \rhoii\tensor(\rhoi'\leftii \Rii S(\tauii))\tensor
  \rhoi''\taui\Ri\dotact z
  \\[-2pt]
  &= \rhoii\tensor \rhoi' \Rii S(\rhoi''\tauii))\tensor
  \rhoi'''\taui\Ri\dotact z
  \\[-2pt]
  &= \rhoii\tensor \rhoi' \Rii S(\tauii\rhoi''')\tensor
  \taui\rhoi''\Ri\dotact z
  \\[-2pt]
  &= \rhoii\tensor  \Rii \rhoi'' S(\tauii\rhoi''')\tensor
  \taui\Ri\rhoi'\dotact z
  \\[-2pt]
  &= \rhoii\tensor  \Rii  S(\tauii)\tensor
  \taui\Ri\rhoi\dotact z
  = z\mone\tensor z\zero{}\bmone\tensor z\zero{}\bzero.
\end{align*}

We next verify the Yetter--Drinfeld axiom relating the $\Nich(X)$
action and coaction, which in the standard graphic notation (see,
e.g., \cite{[Besp-TMF+]}) is expressed as
\begin{equation*}
  \begin{tangles}{l}
    \hstr{75}\step\fobject{\Nich(X)}\step[2]\fobject{\repZ}\\
    \hstr{75}\vstr{90}\cd\step\id\\
    \hstr{75}\vstr{50}\id\step[2]\hx\\
    \hstr{75}\vstr{90}\lu[2]
    \hstr{75}\step\hd\\
    \hstr{75}\vstr{90}\ld[2]\step\ddh\\
    \hstr{75}\vstr{50}\id\step[2]\hx\\
    \hstr{75}\vstr{90}\cu\step\id
  \end{tangles}\ \ = \ \
  \begin{tangles}{l}
    \hstr{75}\step\fobject{\Nich(X)}\step[3]\fobject{\repZ}\\
    \hstr{75}\cd\step\ld\\
    \hstr{75}\id\step[2]\hx\step\id\\
    \hstr{75}\cu\step\lu
  \end{tangles}
\end{equation*}
With the above $\bDelta$, $\bdelta$, and braiding, this is equivalent
to
\begin{multline*}
  \Rii S(\sigmaii)(\piii\taui\rhoi\leftii x'')\tensor
  \pii\sigmai\Ri x' S(\rhoii)\tauii\dotact z
  \\
  = x' S(\rhoii)(\sigmaii\leftii \Rii S(\tauii))\tensor
  (\sigmai\rhoi\leftii x'')\taui\Ri\dotact z
\end{multline*}
for $x\in\Nich(X)$ and $z\in\repZ$. \ We prove the last identity:
\begin{align*}
  \text{l.h.s.}
  &=\Rii S(\sigmaii)(\piii\rhoi\leftii x'')\tensor
  \pii\sigmai\Ri x' S(\rhoii')\rhoii''\dotact z
  \\[-2pt]
  &=\Rii S(\sigmaii)(\piii\leftii x'')\tensor
  \pii\sigmai\Ri x' \dotact z
  \\[-2pt]
  &=\Rii S(\sigmaii')(\sigmaii''\leftii x'')\tensor
  \sigmai\Ri x' \dotact z
  \\[-2pt]
  &=\Rii x'' S(\sigmaii)\tensor \sigmai\Ri x' \dotact z
  \\[-2pt]
  &=x' \Rii S(\sigmaii)\tensor \sigmai x''\Ri \dotact z,
  \\[-6pt]
  \intertext{but}
  \mbox{}\\[-1.3\baselineskip]
  \text{r.h.s.}
  &=x' S(\rhoii')(\rhoii''\leftii \Rii S(\tauii))\tensor
  (\rhoi\leftii x'')\taui\Ri\dotact z
  \\[-2pt]
  &=x' \Rii S(\rhoii\tauii)\tensor
  (\rhoi\leftii x'')\taui\Ri\dotact z
  \\[-2pt]
  &=x' \Rii S(\rhoii)\tensor
  (\rhoi'\leftii x'')\rhoi''\Ri\dotact z,
\end{align*}
which is the same.

\subsubsection{}\label{fun-new-monoidal}
The functor $\functorT$ is monoidal.  Let $\repY$ and $\repZ$ be two
$\algU$ modules.  When viewed as $\Nich(X)$ modules, their tensor
product carries the $\Nich(X)$ action that actually evaluates as the
$\algU$ action on a tensor product of its modules:
\begin{align*}
  x\dotact(y\tensor z)
  &= (x\bup{1}\rhoii\dotact y)\tensor
  (\rhoi\leftii x\bup{2})\dotact z
  \\[-2pt]
  &=(x'S(\tauii)\rhoii\dotact y)\tensor
  (\rhoi\taui\leftii x'')\dotact z
  \\[-2pt]
  &=(x'S(\rhoii')\rhoii''\dotact y)\tensor
  (\rhoi\leftii x'')\dotact z
  \\[-2pt]
  &=x'\dotact y\tensor x''\dotact z.
\end{align*}
Essentially the same calculation shows that the functor is braided:
braiding $c_{\repY,\repZ}$ in the category $\BByd_{\catCrho}$
evaluates as the braiding in $\catCU$:
\begin{equation*}
  (y\bmone\rhoii\dotact z)\tensor \rhoi\dotact y\bzero
  =\Rii\dotact z\tensor\Ri\dotact y.
\end{equation*}

\subsection{The functor $\bfunctorB$}\label{sec:FDB1}
To establish a functor inverse to $\functorT$, we first construct a
functor $\functorB:\BByd_{\catCrho}\to\BBHHYD$, which essentially
amounts to ``a Radford formula for comodules.''

In this subsection, we do not need to assume the existence of a
universal $R$-matrix $\rho\in H\tensor H$; therefore, $H$ is an
arbitrary Hopf algebra (with bijective antipode).

\subsubsection{Reminder: Radford's biproduct}
If $H$ is an ordinary Hopf algebra and $\cR$ is an algebra object in
$\HHyd$ (a braided Hopf algebra), then Radford's
biproduct~\cite{[Radford-bos]} (bosonization~\cite{[Majid-bos]})
defines the structure of an ordinary Hopf algebra on the smash product
$\cR\Smash H$.  Its multiplication, comultiplication, and antipode are
\begin{align}\label{R-smash}
  (r\Smash h)(s\Smash g) &= r(h'\leftii s)\Smash h'' g,
  \\
  \label{R-Delta}
  \pmb{\Delta}:
  r\Smash h
  \mapsto
  (r\Smash h)\zup{1}\tensor(r\Smash h)\zup{2}&=
  (r\bup{1}\Smash r\bup{2}{}\mone h')\tensor
  (r\bup{2}{}\zero\Smash h''),
  \\
  \label{Rad-S}
  \pmb{S}(r\Smash h) &= (1\Smash \hA(r\mone h))
  \bigl(S(r\zero)\Smash 1\bigr)
\end{align}
(where $\hA$ is the antipode of $H$ and $S$ is the antipode of $\cR$).

\subsubsection{}Our aim is to extend bosonization to Yetter--Drinfeld
modules.  In the same setting as above, let $\repY$ be a
Yetter--Drinfeld $\cR$-module in $\HHyd$.  We write $r\leftact y$ and
$h.y$ for the $\cR$ and $H$ actions, and $y\mapsto y\bmone\tensor
y\bzero\in\cR\tensor \cY$ and $y\mapsto y\mone\tensor y\zero\in
H\tensor\cY$ for the coactions.

\begin{prop}\label{prop:bos-YD}
  Let $\cR$ be a Hopf algebra in $\catC=\HHyd$ and
  $\cY\in\RRYD_{\catC}$.  Then $\repY$ is a Yetter--Drinfeld
  $\cR\Smash H$-module with the $\cR\Smash H$ action
  \begin{align}\label{smash-action}
    (r\Smash h)\leftact y &= r\leftact(h.y),
    \\
    \intertext{and coaction}
    \label{smash-coaction}
    \pmb{\delta}:y\mapsto y\zmone{}\tensor y\zzero
    &=(y\bmone\Smash y\bzero{}\mone)\tensor y\bzero{}\zero
    \in\cR\Smash H\tensor\cY.
  \end{align}
  This defines a vector-space-preserving monoidal braided functor
  \begin{equation*}
    \functorB:\RRYD_{\catC}\to\RHRHYD
  \end{equation*}
  to the category of left--left Yetter--Drinfeld $\cR\Smash
  H$-modules.
\end{prop}

\subsubsection{Proof of~\bref{prop:bos-YD}}
First, the coaction property is readily verified
for~\eqref{smash-coaction} using the fact that the $\cR$-coaction is
an $H$-comodule morphism.
Next, we must show the Yetter--Drinfeld condition in $(\cR\Smash
H)\tensor \cY$,
\begin{align}\label{2prove}
 \bigl((r\Smash h)\zup{1}\leftact y\bigr)\zmone (r\Smash h)\zup{2}
 \tensor
 \bigl((r\Smash h)\zup{1}\leftact y\bigr)\zzero\kern-120pt
 \\
 \notag
 &=(r\Smash h)\zup{1} y\zmone\tensor (r\Smash h)\zup{2}\leftact y\zzero,
 \\[-6pt]
 \intertext{where the right-hand side is simple:}
 &=(r\bup{1}\Smash r\bup{2}{}\mone h')(y\bmone\Smash y\bzero{}\mone)
 \tensor (r\bup{2}{}\zero\Smash h'')\leftact y\bzero{}\zero.
 \notag
\end{align}
As regards the left-hand side, we first note that for $z\in \cY$,
$s\in \cR$, and $g\in H$,
\begin{multline*}
  z\zmone(s\Smash g)\tensor z\zzero=
  (z\bmone\Smash z\bzero{}\mone)(s\Smash g)
  \tensor z\bzero{}\zero
  \\
  =
  (z\bmone (z\bzero{}\mone{}'\leftii s)\Smash z\bzero{}\mone{}'' g)
  \tensor z\bzero{}\zero
  =(z\bmone (z\bzero{}\mone{}\leftii s)\Smash z\bzero{}\zero{}\mone{} g)
  \tensor z{}\bzero{}\zero{}\zero.
\end{multline*}
Using this, we readily see directly from the definitions that
\begin{multline*}
  \text{l.h.s.\ of~\eqref{2prove}} ={}\\
  (\bigl(r\bup{1}\leftact(r\bup{2}{}\mone h' . y)\bigr){}\bmone (\bigl(r\bup{1}\leftact(r\bup{2}{}\mone h' . y)\bigr){}\bzero{}\mone{}\leftii r\bup{2}{}\zero)\Smash \bigl(r\bup{1}\leftact(r\bup{2}{}\mone h' . y)\bigr){}\bzero{}\zero{}\mone{} h'')
  \\
  \tensor \bigl(r\bup{1}\leftact(r\bup{2}{}\mone h' . y)\bigr){}\bzero{}\zero{}\zero,
\end{multline*}
where we next use the condition that $\cY$ is a Yetter--Drinfeld
$\cR$-module, which is
\begin{multline*}
  \bigl(r\bup{1}\leftact(r\bup{2}{}\mone.y)\bigr){}\bmone
  \bigl(r\bup{1}\leftact(r\bup{2}{}\mone.y)\bigr){}\bzero{}\mone\leftii
  r\bup{2}{}\zero
  \tensor
  \bigl(r\bup{1}\leftact(r\bup{2}{}\mone.y)\bigr){}\bzero{}\zero
  \\
  =r\bup{1}(r\bup{2}{}\mone\leftii y\bmone)
  \tensor
  (r\bup{2}{}\zero\leftact y\bzero),
\end{multline*}
whence
\begin{align*}
  \text{l.h.s.\ of~\eqref{2prove}}
  &=
  r\bup{1}(r\bup{2}{}\mone \leftii (h'.y)\bmone)
  \Smash
  (r\bup{2}{}\zero\leftact(h' . y)\bzero)\mone h''
  \tensor (r\bup{2}{}\zero\leftact(h' . y)\bzero)\zero.
  \\
  \intertext{We continue by using the fact that the $\cR$-action
    $\leftact$ is a morphism of $H$-comodules,}
  &=
  r\bup{1}(r\bup{2}{}\mone \leftii (h'.y)\bmone)
  \Smash
  r\bup{2}{}\zero{}\mone(h' . y)\bzero{}\mone h''
  \tensor
  r\bup{2}{}\zero{}\zero\leftact(h' . y)\bzero{}\zero
  \\
  \intertext{and that the $\cR$ coaction is a morphism of $H$-modules,}
  &=
  r\bup{1}(r\bup{2}{}\mone \leftii (h'\leftii y\bmone))
  \Smash
  r\bup{2}{}\zero{}\mone(h'' . y\bzero)\mone h'''
  \tensor
  r\bup{2}{}\zero{}\zero\leftact(h'' . y\bzero)\zero
  \\
  &=
  r\bup{1}(r\bup{2}{}\mone{}' h'\leftii y\bmone)
  \Smash
  r\bup{2}{}\mone{}''(h'' . y\bzero)\mone h'''
  \tensor
  r\bup{2}{}\zero\leftact(h'' . y\bzero)\zero,
  \\
  \intertext{after which the $H$-Yetter--Drinfeld condition for $\cY$
    yields}
  &=
  r\bup{1}(r\bup{2}{}\mone{}' h'\leftii y\bmone)
  \Smash
  r\bup{2}{}\mone{}'' h'' y\bzero{}\mone
  \tensor
  r\bup{2}{}\zero\leftact (h'' . y\bzero{}\zero)
  \\
  &=\text{r.h.s.\ of~\eqref{2prove}}.
\end{align*}
This shows that we have a functor.

The functor is monoidal: for $\cY,\cZ\in\RRYD$, the action of $\cR$ on
a tensor product $\cY\tensor\cZ$ is given by the throughout map in
\begin{equation*}
  r\tensor(y\tensor z) \mapsto r\bup{1}\tensor r\bup{2}\tensor
  y\tensor z
  \mapsto
  r\bup{1}\tensor r\bup{2}{}\mone.y\tensor r\bup{2}{}\zero\tensor z
  \mapsto
  \bigl(r\bup{1}\leftact(r\bup{2}{}\mone.y)\bigr)
  \tensor r\bup{2}{}\zero\leftact z.
\end{equation*}
But with $\cY,\cZ\in\RHRHYD$, the $\cR\Smash H$ action on the tensor
product is
\begin{align*}
  (r\Smash h)\tensor(y\tensor z)
  \mapsto
  (&(r\Smash h)\zup{1}\leftact  y)\tensor((r\Smash h)\zup{2}\leftact z)
  \\
  &=(r\bup{1}\Smash r\bup{2}{}\mone h')
  \leftact  y)\tensor((r\bup{2}{}\zero\Smash h'')\leftact z)
  \\
  &=\bigl(r\bup{1}\leftact(r\bup{2}{}\mone h' . y)\bigr)
  \tensor\bigl(r\bup{2}{}\zero\leftact( h''. z)\bigr),
\end{align*}
and hence $r\Smash 1$ acts the same as $r$ in the preceding formula.

The functor is braided: for $y\in\cY$ and $z\in\cZ$, with
$\cY,\cZ\in\RRYD$, the braiding is
\begin{equation*}
  c^{\RRyd}_{\cY,\cZ} : y\tensor z\mapsto
  \bigl(y\bmone\leftact(y\bzero{}\mone.z)\bigr)\tensor y\bzero{}\zero.
\end{equation*}
On the other hand, when $\cY$ and $\cZ$ are viewed as objects in
$\RHRHYD$, the braiding is
\begin{align*}
  c^{\RHRHyd}_{\cY,\cZ} : y\tensor z\mapsto
  (y\zmone\leftact z)\tensor y\zzero
  &= (y\bmone\Smash y\bzero{}\mone)\leftact z\tensor y\bzero{}\zero\\
  &= \bigl(y\bmone\leftact(y\bzero{}\mone.z)\bigr)\tensor y\bzero{}\zero,
\end{align*}
which is the same.

For $\cY\tensor\cZ\in\RRYD$, the coaction is the throughout map in
\begin{equation*}
  y\tensor z
  \mapsto
  y\bmone\tensor y\bzero\tensor z\bmone\tensor z\bzero
  \mapsto
  y\bmone\tensor (y\bzero{}\mone\leftii z\bmone)\tensor
  y\bzero{}\zero\tensor z\bzero
  \mapsto
  y\bmone (y\bzero{}\mone\leftii z\bmone)\tensor
  y\bzero{}\zero\tensor z\bzero,
\end{equation*}
which lies in $\cR\tensor\cY\tensor\cZ$.  Now for
$\cY\tensor\cZ\in\RHRHYD$, the coaction is simply
\begin{align*}
  y\zmone z\zmone\tensor y\zzero\tensor z\zzero
  &=
  (y\bmone\Smash y\bzero{}\mone)(z\bmone\Smash z\bzero{}\mone)
  \tensor y\bzero{}\zero\tensor z\bzero{}\zero\\
  &=
  \bigl(y\bmone(y\bzero{}\mone'\leftii z\bmone)
  \Smash
  (y\bzero{}\mone'' z\bzero{}\mone)\bigr)
  \tensor y\bzero{}\zero\tensor z\bzero{}\zero,
\end{align*}
which is in $\cR\Smash H\tensor\cY\tensor\cZ$.  Applying
$\id\tensor\varepsilon\tensor\id\tensor\id$ restores the previous
formula.

\begin{rem}
  A particular case of~\bref{prop:bos-YD} is well known: for
  $\cY=(\cR,\leftii)$ with the left adjoint action, we
  write~\eqref{smash-action} as (with $r,p\in\cR$)
  \begin{align*}
    (r\Smash 1)\leftii p &=r\bup{1}(r\bup{2}{}\mone
    . p)S(r\bup{2}{}\zero)
    \\
    &=(r\bup{1}\Smash r\bup{2}{}\mone{}')(p\Smash 1) (1\Smash
    \hA(r\bup{2}{}\mone{}''))(S(r\bup{2}{}\zero)\Smash 1)
    \\
    &=r\zup{1}(p\Smash 1)\pmb{S}(r\zup{2})
  \end{align*}
  (see~\eqref{Rad-S} for the antipode), which is a restriction of the
  (``nonbraided'') left adjoint action of $\cR\Smash H$ on itself.
\end{rem}

\subsection{The functor $\bfunctorDB$}\label{sec:FDB}
To construct a ``double bosonization'' functor $\functorDB$ inverse to
$\functorT$, we compose $\functorB$ with the standard functor
establishing the equivalence~\cite{[Mj]} of Yetter--Drinfeld modules
with modules of the Drinfeld double $\DD(\Nich(X)\Smash
H)=\DD(U_{\leq})$, $\functorDD:\BBHHYD\to\catCD$.

Applied to any $\repY\in\BByd_{\catCrho}$, \
$\functorDB=\functorDD\functorB$ then produces a $\algU$-module
because the action of the generators of $H^*\in\DD(U_{\leq})$ is
completely determined by the action of the generators of
$H\in\DD(U_{\leq})$, which means that the module is actually a module
of $\algU\simeq\DD(U_{\leq})/(H^*\sim H)$ (the quotient by relations
expressing $H^*$ through $H$).  Indeed, for $\mu\in H^*$, its action
on $\repY$ viewed as a Yetter--Drinfeld $U_{\leq}$-module is
standardly given by
\begin{align*}
  \mu\drsh y = \eval{\mu}{\pmb{S}^{-1}(y\zmone)}y\zzero &=
  \eval{\varepsilon\tensor\mu}{ \pmb{S}^{-1}(y\bmone\Smash
    y\bzero{}\mone)} y\bzero{}\zero
  \\
  &= \eval{\mu}{S^{-1}(y\mone)}y\zero
  \\
  &= \eval{\mu}{S^{-1}(\rhoii)}\rhoi\dotact y.
\end{align*}
It follows that the generators of $H^*$ (see~\bref{app:DD-rel}) act as
\begin{alignat}{2}\label{Lell}
    L\drsh z &= K\dotact z,\qquad &\ell\drsh z &= k\dotact z,
\end{alignat}
and hence we have a functor
\begin{equation*}
  \functorDB:\BByd_{\catCrho}\to\catCU.
\end{equation*}

\subsection{Composing the functors}\label{sec:TDB}
To verify that $\functorDB \functorT\sim\one_{\catCU}$, we first
calculate $\functorB\functorT$.  Applied to $\Uminus\subset\algU$,
$\functorT$ gives a braided Hopf algebra with
coproduct~\eqref{bDelta}, $x\bup{1}\tensor x\bup{2}=x'S(\rhoii)\tensor
\rhoi\leftii x''$ ($x\in\Uminus$); applying $\functorB$---using
Radford's formula~\eqref{R-Delta}---we obtain the coproduct
$\pmb{\Delta}:U_{\leq}\to U_{\leq}\tensor U_{\leq}$ that evaluates as
\begin{equation*}
  \pmb{\Delta}(x\Smash h)=x' h'\tensor x'' h'',
\end{equation*}
which is the original coproduct on $U_{\leq}$\,.  For
$\repY\in\catCU$, similarly, $\functorT$ produces
coproduct~\eqref{bdelta}, $y\bmone\tensor y\bzero=\Rii
S(\rhoii)\tensor \rhoi\Ri\dotact y$; to further apply $\functorB$, we
substitute the last formula in the ``Radford formula for comodules,''
Eq.~\eqref{smash-coaction}, to obtain the $\Uminus$ coaction
\begin{equation*}
  \pmb{\delta}(y) = \Rii\tensor \Ri\dotact y\in U_{\leq}\tensor\repY.
\end{equation*}
The resulting $U_{\leq}$ action, obviously, is the restriction of the
$\algU$ action.  To the Yetter--Drinfeld $U_{\leq}$-module $\repY$
thus obtained, we now apply $\functorDD$, making it into a module of
the Drinfeld double $\DD(U_{\leq})$.  Then $U_{\leq}^*$\, (a Hopf
subalgebra in $\DD(U_{\leq})$) acts as
\begin{equation*}
  \mu\drsh y = \eval{\mu}{\pmb{S}^{-1}(y\zmone)}y\zzero
  =\eval{\mu}{S^{-1}(\Rii)}\Ri\dotact y
  =\eval{\mu}{\Rii}S(\Ri)\dotact y.
\end{equation*}
With the $R$-matrix in~\bref{thm:R} and with the duality worked out
in~\bref{sec:duality}, we find, along with~\eqref{Lell}, that
\begin{alignat*}{2}
  E\drsh y &= E\dotact y,\qquad &C\drsh y &= C\dotact y.
\end{alignat*}
Comparing with the  formulas in \bref{app:DD-rel} shows that
the resultant $\DD(U_{\leq})$-module is in fact a $\algU$-module,
naturally isomorphic to the original $\repY$.

We show similarly that
$\functorT\functorDB\sim\one_{^{\Nich(X)}_{\Nich(X)}\YDname_{\catCrho}}$.
Starting from $Z\in\BByd_{\catCrho}$ (where we write the action of
both $H$ and $\Uminus$ as $x\tensor z\mapsto x\dotact z$ for
simplicity, and the coaction as $z\mapsto z\bmone\tensor z\bzero$) and
applying $\functorT\functorDB$, we arrive at the Yetter--Drinfeld
$\Uminus$-module with the coaction
\begin{align*}
  z\mapsto{}&\rhoii\leftii\Rii\tensor\eval{\bRi}{
    \pmb{S}^{-1}(z\bmone\Smash\tauii)}\rhoi\taui\dotact z\bzero
  \\
  &=\rhoii''\leftii\Rii\tensor\eval{\bRi}{
    \pmb{S}^{-1}(z\bmone\Smash\rhoii')}\rhoi\dotact z\bzero
  \\
  &=\rhoii''\leftii\Rii\tensor\eval{\bRi}{
    S^{-1}(\rhoii')\leftii S^{-1}(z\bmone)}
  \rhoi\dotact z\bzero
  \\
  &=\rhoii''\leftii\Rii\tensor\eval{\rhoii'\leftii\bRi}{
     S^{-1}(z\bmone)}
  \rhoi\dotact z\bzero
  \\
  &=\Rii\tensor\eval{\bRi}{S^{-1}(z\bmone)}
  z\bzero
  =z\bmone\tensor 
  z\bzero.
\end{align*}
Hence, we have constructed an inverse functor $\functorDB$ to the
functor $\functorT$ in~\bref{thm:fun-new}.  The categories $\catCU$
and $\BByd_{\catCrho}$ are equivalent.  In the rest of the paper, we
study the modules of $\algU$.

\section{\textbf{Simple modules of $\balgU$}}\label{sec:simples}
The following theorem is of course contained in
\cite{[ARS],[HS-double]}; we spell out the details here in order to
fix our notation and conventions.
\begin{Thm}
  The algebra $\algU$ has $4p^2$ simple modules, labeled as
  \begin{gather*}
    \repZ^{\alpha,\beta}_{s,r},\qquad \alpha,\beta=\pm,\quad
    s=1,\dots,p,\quad r=0,\dots,p-1,
  \end{gather*}
  with the dimensions
  \begin{equation}\label{dim-simp-mod}
    \dim\repZ^{\alpha,\beta}_{s,r}=
    \begin{cases}
      2s-1,&1\leq s\leq p,\quad r=0,\\
      2s+1,& 1\leq s\leq p-1,\quad r=s,\\
      4s,&1\leq s\leq p-1,\quad r\neq0,s,\\
      4p,& s=p,\quad 1\leq r\leq p-1.
    \end{cases}
  \end{equation}
  The modules in the first two lines are atypical, and the others
  typical.
\end{Thm}
The cases occurring in the theorem can be illustrated by the diagram
in Fig.~\ref{module-pattern}.
\begin{figure}[tbh]
  \centering
  \begin{gather*}
    \xymatrix@C=0pt@R=0pt{ &&
      \\
      &&
      \\
      &*{}\ar@{{}{-}{}}[];[rrrrrrrrrr]&&&&&&*{}&&*{}&&*{}\\
      {\scriptstyle p-1}&&&&&&&&*{\Zpp}&&*{\Zpr}&\\
      &&&&&*{}\ar@{{}{-}{}}[];[rrrr]&&&&*{}&&\\
      &&&&&&*{\Zss}&&&&*{\Zpr}&\\
      &&&*{}\ar@{{}{-}{}}[];[rrrr]&&&&*{}\ar@{{}{-}{}}[];[uuuu]&&&&\\
      &&&&*{\Zss}&&&&&&*{\Zpr}&\\
      &*{}\ar@{{}{-}{}}[];[rrrr]&&&&*{}\ar@{{}{-}{}}[];[uuuu]&&&&&&\\
      &&*{\Zss}&&&&&&&&*{\Zpr}
      \\
      &*{}\ar@{{}{-}{}}[];[rrrrrrrrrr]&&&&&&&&*{}&&*{}
      \\
      {\scriptstyle 0}&&*{\Zoz}&&*{\Zso}&&*{\Zso}&&*{\Zso}&&*{\Zpo}&
      \\
      &*{}\ar[uuuuuuuuuuuu]^(.95)r
      \ar@{{}{-}{}}[];[rrrrrrrrrr]*{}\ar@{{}{-}{}}[];[uuuuuuuuuu]
      \ar[rrrrrrrrrrrr]^(.95)s&&*{}\ar@{{}{-}{}}[];[uuuuuu]&&&&&&*{}\ar@{{}{-}{}}[];[uuuuuuuuuu]&&*{}\ar@{{}{-}{}}[];[uuuuuuuuuu]&&
      \\
      &&{\scriptstyle 1}&&&&&&&&{\scriptstyle p}& }
  \end{gather*}
  \caption{\small The $p^2$ simple $\algU$ modules, labeled by $s$ and
    $r$.}
  \label{module-pattern}
\end{figure}
Typical modules are the ``bulk'' of the diagram and the greatest part
of the right column.  In what follows, we refer to these last modules
(the fourth case in~\eqref{dim-simp-mod}) as Steinberg modules.  The
notation in the figure distinguishes more cases than just those
in~\eqref{dim-simp-mod} for our later purposes.

\subsection{Constructing simple
  $\balgU$-modules}\label{simple-U-modules}
We construct the simple $\algU$ modules explicitly, by ``gluing
together'' some modules of the algebra $\Uresk$ defined in
\eqref{Ures} and \eqref{Ustar}.  We let simple $\Uresk$ modules be
denoted by $\repX^{\alpha,\beta}_{s,r}$, where $\alpha,\beta=\pm$,
$s=1,\dots,p$, and $r=0,\dots,p-1$.  As a $\Ures$-module, an
$\repX^{\alpha,\beta}_{s,r}$ with any $\beta$ and $r$ is isomorphic to
$\repX^{\alpha}_{s}$ (see~\bref{app:SL2-irrep}). In particular, there
is a highest-weight vector $\ket{\alpha,s,\beta,r}_0=0$ such that
\begin{gather*}
  E \ket{\alpha,s,\beta,r}_0=0 \quad\text{and}\quad
  K\ket{\alpha,s,\beta,r}_0=\alpha\q^{s-1}\ket{\alpha,s,\beta,r}_0,
\end{gather*}
and the second Cartan generator of $\algU$ acts on this vector as
\begin{gather*}
  k\ket{\alpha,s,\beta,r}_0=\beta\q^{-r}\ket{\alpha,s,\beta,r}_0.
\end{gather*}
It follows that $\dim\repX^{\alpha,\beta}_{s,r}=s$.

\subsubsection{$\Uresk$ decompositions and bases of simple
  modules}\label{simple-list}
For each of the four cases in~\eqref{dim-simp-mod}, we now list the
$\Uresk$ decompositions of simple $\algU$ modules, specify the
corresponding choice of basis, and describe how the $\Uresk$
constituents are glued together by the fermionic $\algU$ generators
$B$ and $C$.  \ There are two sorts of basis vectors $\ZC{~}_m$ and
$\ZB{~}_m$ for atypical modules and two more, $\ZU{~}_m$ and
$\ZD{~}_m$, for typical modules; in all cases,
\begin{equation*}
  C\ZC{~}_m=0\quad\text{and}\quad B\ZB{~}_m=0
\end{equation*}
in any simple $\algU$ module.

\begin{description}\addtolength{\itemsep}{6pt}
\item[Atypical modules, $\myboldsymbol{1\leq s\leq p}$,
  $\myboldsymbol{r=0}$] as $\Uresk$ modules, these simple $\algU$
  modules decompose as
  \begin{equation}\label{Zs0-decomp}
    \repZ^{\alpha,\beta}_{s,0}=
    \repX^{\alpha,\beta}_{s,0}\oplus\repX^{\alpha,\beta}_{s-1,p-1},
    \quad 2\leq s\leq p,
  \end{equation}
  which is illustrated in Fig.~\ref{fig:ZZZ},\;left, and as
  $\repZ^{\alpha,\beta}_{1,0}=\repX^{\alpha,\beta}_{1,0}$
  (one-dimensional modules).
  \begin{figure}[tb]
    \centering
    \begin{equation*}
      \xymatrix@C=10pt@R=5pt{
        \blacktriangleleft&&&&&&&\blacktriangleright&
        \\
        &\blacktriangleright&&&&&\blacktriangleleft\ar[ur]^{B}\ar[dd]_{F}&&
        \\
        \blacktriangleleft&&&&&&&\blacktriangleright&
        \\
        &\blacktriangleright&&&&&\blacktriangleleft&&
        \\
        \blacktriangleleft&&&&&&&\blacktriangleright&
        \\
        &\blacktriangleright&&&&&\blacktriangleleft&&
        \\
        \blacktriangleleft&&&&&&&\blacktriangleright&
        \\
        &\blacktriangleright&&&&&\blacktriangleleft&&
        \\
        \blacktriangleleft&&&&&&&\blacktriangleright&
      }\qquad\qquad\qquad
      \xymatrix@C=10pt@R=5pt{
        &\blacktriangleup\ \ &\\
        \blacktriangleleft&&\blacktriangleright\\
        &\blacktriangleup\
        \blacktriangledown&\\
        \blacktriangleleft&&\blacktriangleright\\
        &\blacktriangleup\
        \blacktriangledown&\\
        \blacktriangleleft&&\blacktriangleright\\
        &\blacktriangleup\
        \blacktriangledown&\\
        \blacktriangleleft&&\blacktriangleright\\
        &\blacktriangleup\ \ &
      }
    \end{equation*}
    \caption[Simple $\algU$ modules]{\small\textsc{Left:} An atypical
      module $\repZ^{\alpha,\beta}_{s,0}$ (with $s=5$).  Each vertical
      column is a $\Uresk$ module in~\eqref{Zs0-decomp}.  The top
      state is $\ZC{\alpha,s,\beta,0}_0$ and the bottom,
      $\ZC{\alpha,s,\beta,0}_{s-1}$. \ \textsc{Middle:} An atypical
      module $\repZ^{\alpha,\beta}_{s,s}$ (with $s=4$).  Each vertical
      column is a $\Uresk$ module in~\eqref{Zss-decomp}.  The top
      state is $\ZB{\alpha,s,\beta,s}_0$ and the bottom,
      $\ZB{\alpha,s,\beta,s}_s$. \ \textsc{Right:} The typical module
      $\repZ^{\alpha,\beta}_{s,r}$ ($s=4$).  Each column is a
      $\Uresk$-module in~\eqref{four-sub}. \ The directions in which
      the generators map are common for all modules.}
    \label{fig:ZZZ}
  \end{figure}
  We choose a basis in $\repZ^{\alpha,\beta}_{s,0}$ in accordance with
  this decomposition, as
  \begin{equation*}
    \bigl(
    \ZC{\alpha,s,\beta,0}_n\in\repX^{\alpha,\beta}_{s,0}\bigr)_{
      0\leq n\leq s-1},\qquad
    \bigl(
    \ZB{\alpha,s,\beta,0}_m\in\repX^{\alpha,\beta}_{s-1,p-1}\bigr)_{
      0\leq m\leq s-2}.
  \end{equation*}
  The arrows in the notation for basis vectors refer to the
  visualization of $\repZ^{\alpha,\beta}_{s,0}$ as in
  Fig.~\ref{fig:ZZZ},\;left.  The fermionic generators relate the two
  types of vectors as
  \begin{alignat*}{2}
    B \ZC{\alpha,s,\beta,0}_n&=-\qint{n}\ZB{\alpha,s,\beta,0}_{n-1},
    \qquad
    &C\ZB{\alpha,s,\beta,0}_m&=\beta\ZC{\alpha,s,\beta,0}_{m+1}.
  \end{alignat*}

  Here and hereafter, we set $\ket{\alpha,s,\beta,r}^{\bullet}_m=0$
  whenever $m$ is outside the range indicated for vectors of a given
  type.

\item[Atypical modules, $\myboldsymbol{1\leq s\leq p-1}$,
  $\myboldsymbol{r=s}$] the module decomposes as
  \begin{equation}\label{Zss-decomp}
    \repZ^{\alpha,\beta}_{s,s}=\repX^{\alpha,\beta}_{s,s}\oplus
    \repX^{\alpha,-\beta}_{s+1,s},
  \end{equation}
  as is illustrated in Fig.~\ref{fig:ZZZ},\;middle.  We choose a basis
  in $\repZ^{\alpha,\beta}_{s,s}$ accordingly:
  \begin{equation*}
    \bigl(
    \ZC{\alpha,s,\beta,s}_n\in\repX^{\alpha,\beta}_{s,s}\bigr)_{
      0\leq n\leq s-1},\qquad
    \bigl(
    \ZB{\alpha,s,\beta,s}_m\in\repX^{\alpha,-\beta}_{s+1,s}\bigr)_{
      0\leq m\leq s}.
  \end{equation*}
  The fermions act as
  \begin{alignat*}{2}
    B \ZC{\alpha,s,\beta,s}_n&=\qint{s-n}\ZB{\alpha,s,\beta,s}_{n},
    \qquad
    &C \ZB{\alpha,s,\beta,s}_m&=\beta\ZC{\alpha,s,\beta,s}_{m}.
  \end{alignat*}

\item[Typical modules, $1\myboldsymbol{\leq s\leq p-1}$,
  $\myboldsymbol{r\neq0,s}$] the modules decompose as
  \begin{equation}\label{bulk-decomp}
    \repZ^{\alpha,\beta}_{s,r}=\repX^{\alpha,\beta}_{s,r}
    \oplus\repX^{\alpha,-\beta}_{s+1,r}
    \oplus\repX^{\alpha,-\beta}_{s-1,r-1}
    \oplus\repX^{\alpha,\beta}_{s,r-1},
    \quad 2\leq s\leq p-1,
  \end{equation}
  which is illustrated in Fig.~\ref{fig:ZZZ},\;right, and as
  $\repZ^{\alpha,\beta}_{1,r}=\repX^{\alpha,\beta}_{1,r}
  \oplus\repX^{\alpha,-\beta}_{2,r}\oplus\repX^{\alpha,\beta}_{1,r-1}$.
  We choose a basis in $\repZ^{\alpha,\beta}_{s,r}$ accordingly, such
  that the bases in respective $\Uresk$ submodules are
  \begin{multline}\label{four-sub}
    \bigl(\ZC{\alpha,s,\beta,r}_j\bigr)_{0\leq j\leq s-1}, \
    \bigl(\ZU{\alpha,s,\beta,r}_m\bigr)_{0\leq m\leq s},
    \\
    \bigl(\ZD{\alpha,s,\beta,r}_n \bigr)_{0\leq n\leq s-2}, \
    \bigl(\ZB{\alpha,s,\beta,r}_j\bigr)_{0\leq j\leq s-1}.
  \end{multline}
  The fermions glue the $\Uresk$ modules together as
  \begin{alignat*}{2}
    B \ZC{\alpha, s, \beta, r}_{j} &= \fffrac{\qint{j}}{\qint{s}}
    \ZD{\alpha, s, \beta, r}_{j - 1} + \beta \fffrac{\qint{r}\qint{s -
        j}}{\qint{s}} \ZU{\alpha, s, \beta, r}_{j},\kern-200pt
    \\
    B \ZU{\alpha, s, \beta, r}_{m} &= \qint{m} \ZB{\alpha, s, \beta,
      r}_{m - 1}, & C \ZU{\alpha, s, \beta, r}_{m} &= \ZC{\alpha, s,
      \beta, r}_{m},
    \\
    B \ZD{\alpha, s, \beta, r}_{n} &= \beta \qint{r}
    \qint{n\!+\!1\!-\!s} \ZB{\alpha, s, \beta, r}_{n}, \quad& C
    \ZD{\alpha, s, \beta, r}_{n} &= \beta \qint{r\!-\!s} \ZC{\alpha,
      s, \beta, r}_{n + 1},
    \\
    C \ZB{\alpha, s, \beta, r}_{j} &= \fffrac{1}{\qint{s}} \ZD{\alpha,
      s, \beta, r}_{j} + \beta \fffrac{\qint{s - r}}{\qint{s}}
    \ZU{\alpha, s, \beta, r}_{j + 1}.\kern-100pt
  \end{alignat*}

\item[Steinberg modules, $\myboldsymbol{s=p}$, $\myboldsymbol{1\leq
    r\leq p-1}$] the decomposition is
  \begin{equation}\label{stein-decomp}
    \repZ^{\alpha,\beta}_{p,r} =\repX^{\alpha,\beta}_{p,r}\oplus
    \repP^{\alpha,-\beta}_{p-1,r-1}\oplus \repX^{\alpha,\beta}_{p,r-1},
  \end{equation}
  where $\repP^{\alpha,\beta}_{p-1,r}$ is a projective
  $\Uresk$-module.  This can be illustrated with much the same diagram
  as in Fig.~\ref{fig:ZZZ},\;right, but with two differences in the
  middle columns: first, there are exactly $p$ $\blacktriangleup$ and
  $p$ $\blacktriangledown$ vectors and, second, they are not a direct
  sum of two simple $\Uresk$ modules; instead, the action of $F$
  and~$E$ generators glues them together~as
  
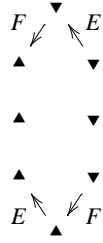
\begin{equation*}
    \xymatrix@C=0pt@R=8pt{
      &\blacktriangledown\ar_F[ld]&\\
      \blacktriangleup
      &&
      \blacktriangledown\ar_E[ul]\\
      \blacktriangleup
      &&
      \blacktriangledown\\
      \blacktriangleup&&
      \blacktriangledown\ar^{F}[dl]\\
      &\blacktriangleup\ar^{E}[ul]&
    }
  \end{equation*}
  (the vectors are conventionally separated horizontally; we remind
  the reader that this is only the middle part of a diagram similar to
  the Fig.~\ref{fig:ZZZ},\;right).  We choose the basis in
  $\repZ^{\alpha,\beta}_{p,r}$ in accordance
  with~\eqref{stein-decomp}, as
  \begin{equation*}
    \ZC{\alpha,p,\beta,r}_n,\quad
    \ZU{\alpha,p,\beta,r}_n,\quad
    \ZD{\alpha,p,\beta,r}_n,\quad
    \ZB{\alpha,p,\beta,r}_n
  \end{equation*}
  where $n=0,\dots,p-1$.  The fermions act as 
  \begin{alignat*}{2}
    B \ZC{\alpha, p, \beta, r}_{n} &= \qint{n - 1} \ZU{\alpha, p,
      \beta, r}_{n} + \qint{n} \ZD{\alpha, p, \beta, r}_{n - 1},
    \kern-120pt &
    \\
    B \ZU{\alpha, p, \beta, r}_{n} &= -\qint{n} \ZB{\alpha, p, \beta,
      r}_{n - 1}, \quad& C \ZU{\alpha, p, \beta, r}_{n} &= -\beta
    \qint{r} \ZC{\alpha, p, \beta, r}_{n},
    \\
    B \ZD{\alpha, p, \beta, r}_{n} &= \qint{n} \ZB{\alpha, p, \beta,
      r}_{n}, & C \ZD{\alpha, p, \beta, r}_{n} &= \beta \qint{r - 1}
    \ZC{\alpha, p, \beta, r}_{n + 1},
    \\
    C \ZB{\alpha, p, \beta, r}_{n} &= \beta \qint{r - 1} \ZU{\alpha,
      p, \beta, r}_{n + 1} + \beta \qint{r} \ZD{\alpha, p, \beta,
      r}_{n}.\kern-120pt
  \end{alignat*}
\end{description}

The remaining formulas for the $\algU$ action are collected
in~\bref{action-on-simple}.

\subsection{Casimirs from simple modules}
We recall the Drinfeld map~\eqref{dr-map}.  If $\repZ$ is a
$\algU$-module, then the (``quantum'') trace operation
\begin{equation}\label{q-trace}
  \mathsf{Tr}_{\repZ}:A\mapsto \mathrm{tr}_{\repZ}(A K^p k^2):
  \algU\to\oC,
\end{equation}
where $\mathrm{tr}_{\repZ}$ evaluates the ordinary trace in any chosen
basis in $\repZ$, defines an element of $\Ch$.  We note that if
$\gamma=\mathsf{Tr}_{\repZ}$, then it follows from~\eqref{Uixi} that
$\widetilde{\Ch}\ni\widetilde{\gamma}:A\mapsto \mathrm{tr}_{\repZ}(A
K^p)$.  Calculating with~\eqref{M-matrix}, we see that traces over
three-dimensional representations are mapped by the Drinfeld map into
the Casimir elements:
\begin{align*}
  \drmap:\mathsf{Tr}_{\repZ^{\alpha,\beta}_{1,1}} &\mapsto -\alpha^p
  \q^{-1} (\q - \q^{-1})^2 C_1,
  \\
  \drmap:\mathsf{Tr}_{\repZ^{\alpha,\beta}_{2,0}} &\mapsto -\alpha^p
  \q (\q - \q^{-1})^2 C_2.
\end{align*}

\section{\textbf{$\bExti$ spaces for simple
    $\balgU$-modules}}\label{sec:ext}
We now describe the $\Exti$ groups for simple $\algU$-modules.  We
recall that for two modules $\repZ_2$ and $\repZ_1$,
$\Exti(\repZ_1,\repZ_2)$ is a linear space with the basis identified
with nontrivial short exact sequences
\begin{equation*}
  0\to\repZ_2\to\repZ_1\oright\repZ_2\to\repZ_1\to0
\end{equation*}
modulo a certain equivalence relation~\cite{[McL]}.

The action of any $\algU$ generator $A$ on $\repZ_1\oleft\repZ_2$ is
given by
\begin{equation*}
  \rho_{A}=\rho^{(0)}_{A}+\xi_A
\end{equation*}
where $\rho^{(0)}_{A}$ is the direct sum of the actions of $\algU$
generators on simple modules and
$\xi_A=\xi^{\repZ_1,\repZ_2}_A:\repZ_1\to\repZ_2$ are linear maps
(also linear in $A$).  We list the $\xi_A$ for each extension in what
follows.

\subsection{$\bExti$ spaces for typical simple
  modules}\label{sec:Ext1-typ} The nontrivial
$\Exti(\repZ_1,\repZ_2)$ spaces for the typical simple $\algU$-modules
are one-dimensional.  These are $\Exti(\repZ^{\alpha, \beta}_{s,
  r},\repZ^{-\alpha, -\beta}_{p - s, p + r - s})$ and
$\Exti(\repZ^{\alpha, \beta}_{s, r},\repZ^{-\alpha, \beta}_{p - s, p +
  r - s})$ for each pair $(s,r)$ such that
\begin{gather*}
  1\leq s\leq p-1,\quad 1\leq r\leq p-1,\quad s\neq r.
\end{gather*}
To avoid notational complications, we here adopt the convention that
$\repZ^{\alpha, \beta}_{s, p + r}=\repZ^{\alpha, -\beta}_{s, r}$ for
$1\leq r\leq p-1$.

We fix a basis element in each of the $\Exti$ spaces, writing
\begin{align}\label{basis-f}
  \Exti(\repZ^{\alpha, \beta}_{s, r}, \repZ^{-\alpha, -\beta}_{p - s,
    p + r - s}) &= \{f_{p - s, p + r - s}\},
  \\
  \intertext{and}
  \label{basis-e}
  \Exti(\repZ^{\alpha, \beta}_{s, r}, \repZ^{-\alpha, \beta}_{p - s, p
    + r - s}) &= \{e_{p - s, p + r - s}\}.
\end{align}
Here and hereafter, $\alpha$ and $\beta$ are omitted in the
right-hands sides in order not to overburden the notation; they are in
all cases easily reconstructed from the context.  The maps
$\repZ^{\alpha, \beta}_{s, r}\to\repZ^{-\alpha, -\beta}_{p - s, p + r
  - s}$ by the $\algU$ generators associated with~\eqref{basis-f} are
\begin{align}
  \xi_{F} : \ZU{\alpha, s, \beta, r}_{m} &\mapsto \delta_{m, s}
  \ZD{-\alpha, p - s, -\beta, p + r - s}_{0},
  \nonumber\\
  \xi_{C} : \ZU{\alpha, s, \beta, r}_{m} &\mapsto -\beta \qint{r}
  \delta_{m, s} \ZC{-\alpha, p - s, -\beta, p + r - s}_{0},
  \nonumber\\
  \xi_{F} : \ZD{\alpha, s, \beta, r}_{m} &\mapsto \qint{s - r}
  \qint{r} \delta_{m, s - 2} \ZU{-\alpha, p - s, -\beta, p + r -
    s}_{0},
  \nonumber\\[-.5\baselineskip]
  \mbox{}\label{xi-maps-f}\\[-.5\baselineskip]
  \xi_{F} : \ZC{\alpha, s, \beta, r}_{m} &\mapsto -\beta \delta_{m, s
    - 1} \qint{r} \ZC{-\alpha, p - s, -\beta, p + r - s}_{0},
  \nonumber\\
  \xi_{F} : \ZB{\alpha, s, \beta, r}_{m} &\mapsto \beta \qint{s - r}
  \delta_{m, s - 1} \ZB{-\alpha, p - s, -\beta, p + r - s}_{0},
  \nonumber\\
  \xi_{C} : \ZB{\alpha, s, \beta, r}_{m} &\mapsto \delta_{m, s -
    1}\fffrac{\qint{r}\qint{s - r}}{\qint{s}} \ZU{-\alpha, p - s,
    -\beta, p + r - s}_{0} \nonumber
\end{align}
(and zero otherwise).  The maps $\repZ^{\alpha, \beta}_{s,
  r}\to\repZ^{-\alpha, \beta}_{p - s, p + r - s}$ by the $\algU$
generators associated with~\eqref{basis-e} are
\begin{align}
  \xi_{E} : \ZU{\alpha, s, \beta, r}_{m} &\mapsto \fffrac{\qint{s +
      1}}{\qint{r - s} \qint{r}} \delta_{m, 0} \ZD{-\alpha, p - s,
    \beta, p + r - s}_{p - s - 2},
  \nonumber\\
  \xi_{E} : \ZC{\alpha, s, \beta, r}_{m} &\mapsto \beta
  \fffrac{\qint{s}}{\qint{r - s}} \delta_{m, 0} \ZC{-\alpha, p - s,
    \beta, p + r - s}_{p - s - 1},
  \nonumber\\
  \xi_{E} : \ZD{\alpha, s, \beta, r}_{m} &\mapsto \qint{s - 1}
  \delta_{m, 0} \ZU{-\alpha, p - s, \beta, p + r - s}_{p - s},
  \nonumber\\[-.5\baselineskip]
  \mbox{}\label{xi-maps-e}\\[-.5\baselineskip]
  \xi_{B} : \ZC{\alpha, s, \beta, r}_{m} &\mapsto
  \alpha\fffrac{1}{\qint{s}} \delta_{m, 0} \ZU{-\alpha, p - s, \beta,
    p + r - s}_{p - s},
  \nonumber\\
  \xi_{B} : \ZU{\alpha, s, \beta, r}_{m} &\mapsto -\alpha
  \beta\fffrac{1}{\qint{r}} \delta_{m, 0} \ZB{-\alpha, p - s, \beta, p
    + r - s}_{p - s - 1},
  \nonumber\\
  \xi_{E} : \ZB{\alpha, s, \beta, r}_{m} &\mapsto -\beta
  \fffrac{\qint{s}}{\qint{r}} \delta_{m, 0} \ZB{-\alpha, p - s, \beta,
    p + r - s}_{p - s - 1}.  \nonumber
\end{align}
That these formulas make $\repZ^{\alpha, \beta}_{s, r}\oleft
\repZ^{-\alpha, -\beta}_{p - s, p + r - s}$ and $\repZ^{\alpha,
  \beta}_{s, r}\oleft \repZ^{-\alpha, \beta}_{p - s, p + r - s}$ into
$\algU$ modules is verified directly.

\subsection{$\bExti$ spaces for atypical simple
  $\balgU$-modules}\label{sec:Ext1-atyp}
The $\Exti(\repZ_1,\repZ_2)$ groups for atypical simple
$\algU$-modules are either $1$- or $0$-dimen\-sional.  The nonzero
$\Exti$ spaces can be arranged into several series, in addition to
which there are a few ``exceptional'' cases, where one of the modules
involved is $\repZ^{\alpha,\beta}_{1,0}$.  These cases can be
conveniently absorbed into the series by adopting the convention that
$\repZ^{\alpha,\beta}_{0,0}=\repZ^{\alpha,-\beta}_{1,0}$ and setting
$\ZB{\alpha,0,\beta,0}_0=\ZC{\alpha, 1,-\beta, 0}_{0}$ and
$\ZC{\alpha,0,\beta,0}_m=0$, $m\neq0$ for the basis vectors.  Then the
extensions are defined by the following formulas, where we choose a
basis vector in each space, placing it in curly brackets, and specify
how the $\algU$ generators map from $\repZ_1$ to $\repZ_2$ (again
omitting the $\alpha$ and $\beta$ indices from the notation)
\begin{align*}
  \myatop{
    \Exti(\repZ^{\alpha,\beta}_{s,0},\repZ^{\alpha,-\beta}_{s+1,0})=\{b_{s+1}\},}{
    1\leq s\leq p-1
  }\quad&
  \begin{aligned}
    \xi_{B}&:\ZC{\alpha,s,\beta,0}_m\mapsto-\qint{s - m}\ZC{\alpha, s+1,-\beta, 0}_m,\\
    \xi_{B}&:\ZB{\alpha,s,\beta,0}_m\mapsto \qint{s - m-1}\ZB{\alpha,
      s+1, -\beta, 0}_m,
  \end{aligned}
  \\[4pt]
  \myatop{
    \Exti(\repZ^{\alpha,\beta}_{s,0},\repZ^{\alpha,-\beta}_{s-1,0})=\{c_{s-1}\},}{
    2\leq s\leq p
  }\quad&
  \begin{aligned}
    \xi_{C}&:\ZC{\alpha,s,\beta,0}_m\mapsto \ZC{\alpha, s-1, -\beta, 0}_m,\\
    \xi_{C}&:\ZB{\alpha,s,\beta,0}_m\mapsto \ZB{\alpha, s-1, -\beta,
      0}_m,
  \end{aligned}
  \\[4pt]
  \myatop{
    \Exti(\repZ^{\alpha,\beta}_{s,s},\repZ^{-\alpha,\beta}_{p-s,0})=\{f_{p-s}\},}{
    0\leq s\leq p-1
  }\quad&
  \begin{aligned}
    \xi_{F}&:\ZC{\alpha,s,\beta,s}_{s-1}\mapsto \ZC{-\alpha, p-s, \beta, 0}_0,\\
    \xi_{F}&:\ZB{\alpha,s,\beta,s}_{s}\mapsto \ZB{-\alpha, p-s, \beta, 0}_0,\\
    \xi_{C}&:\ZB{\alpha,s,\beta,0}_{s}\mapsto\beta \ZC{-\alpha, p-s,
      \beta, 0}_0,
  \end{aligned}
  \\[4pt]
  \myatop{
    \Exti(\repZ^{\alpha,\beta}_{s,s},\repZ^{-\alpha,-\beta}_{p-s,0})=\{e_{p-s}\},}{
    0\leq s\leq p-1
  }\quad&
  \begin{aligned}
    \xi_{E}&:\ZC{\alpha,s,\beta,s}_0\mapsto \qint{s} 
    \ZC{-\alpha, p-s, -\beta, 0}_{p-s-1},\\
    \xi_{E}&:\ZB{\alpha,s,\beta,s}_0\mapsto -\qint{s +
      1} 
    \ZB{-\alpha, p-s, -\beta, 0}_{p-s-2},
  \end{aligned}
  \\[4pt]
  \myatop{
    \Exti(\repZ^{\alpha,\beta}_{s,s},\repZ^{\alpha,-\beta}_{s-1,s-1})=\{\bar{b}_{s-1}\},}{
    1\leq s\leq p-1
  }\quad&
  \begin{aligned}
    \xi_{B}&:\ZC{\alpha,s,\beta,s}_m\mapsto-\qint{m}\ZC{\alpha, s-1,-\beta, s-1}_{m-1},\\
    \xi_{B}&:\ZB{\alpha,s,\beta,s}_m\mapsto\qint{m}\ZB{\alpha, s-1,
      -\beta, s-1}_{m-1},
  \end{aligned}
  \\[4pt]
  \myatop{
    \Exti(\repZ^{\alpha,\beta}_{s,s},\repZ^{\alpha,-\beta}_{s+1,s+1})=\{\bar{c}_{s+1}\},}{
    0\leq s\leq p-2
  }\quad&
  \begin{aligned}
    \xi_{C}&:\ZC{\alpha,s,\beta,s}_m\mapsto \ZC{\alpha, s+1, -\beta, s+1}_{m+1},\\
    \xi_{C}&:\ZB{\alpha,s,\beta,s}_m\mapsto \ZB{\alpha, s+1, -\beta,
      s+1}_{m+1},
  \end{aligned}
  \\[4pt]
  \myatop{
    \Exti(\repZ^{\alpha,\beta}_{s,0},\repZ^{-\alpha,-\beta}_{p-s,p-s})=\{\bar{f}_{p-s}\},}{
    1\leq s\leq p
  }\quad&
  \begin{aligned}
    \xi_{F}&:\ZC{\alpha,s,\beta,0}_{s-1}\mapsto \ZC{-\alpha, p-s, -\beta, p-s}_0,\\
    \xi_{F}&:\ZB{\alpha,s,\beta,0}_{s-2}\mapsto -\ZB{-\alpha, p-s,
      -\beta, 0}_0,
  \end{aligned}
  \\[4pt]
  \myatop{
    \Exti(\repZ^{\alpha,\beta}_{s,0},\repZ^{-\alpha,\beta}_{p-s,p-s})=\{\bar{e}_{p-s}\},}{
    1\leq s\leq p
  }\quad&
  \begin{aligned}
    \xi_{E}&:\ZC{\alpha,s,\beta,0}_0\mapsto\qint{s} 
    \ZC{-\alpha, p-s, \beta, p-s}_{p-s-1},\\
    \xi_{E}&:\ZB{\alpha,s,\beta,0}_0\mapsto \qint{s - 1}\ZB{-\alpha, p-s, \beta, p-s}_{p-s},\\
    \xi_{B}&:\ZC{\alpha,s,\beta,0}_{0}\mapsto-\alpha\ZB{-\alpha, p-s,
      \beta, p-s}_{p-s},
  \end{aligned}
\end{align*}
with all other generators mapping by zero in each case.

\subsection{Linkage classes}\label{sec:link-cl}
It follows that the simple $\algU$ modules are divided into linkage
classes as follows.
\begin{enumerate}
\item There are $4(p-1)$ Steinberg linkage classes, labeled by
  $\alpha=\pm$, $\beta=\pm$, and $1\leq r\leq p-1$; each class
  contains a single module $\repZ^{\alpha,\beta}_{p,r}$.

\item\label{lc-typ} There are $(p-1)(p-2)$ typical linkage classes,
  labeled by $(\alpha, s, r)$ with $\alpha=\pm$ and $1\leq r < s\leq
  p-1$; each such class contains four simple modules
  $\repZ^{\alpha,\pm}_{s,r}$, $\repZ^{-\alpha,\pm}_{p-s,p+r-s}$.

\item There is one atypical linkage class containing $4(2p-1)$ simple
  modules: $\repZ^{\pm,\pm}_{s,s}$ with $1\leq s\leq p-1$ and
  $\repZ^{\pm,\pm}_{s,0}$ with $1\leq s\leq p$ (with uncorrelated
  signs in either case).
  
\end{enumerate}

\section{\textbf{Projective $\balgU$ modules}}\label{sec:proj}
\begin{Thm}
  Projective modules of $\algU$ are exhausted by the following list
  \textup{(}with $\alpha,\beta=\pm$ in all cases\textup{)}.
  \begin{itemize}
  \item Steinberg modules $\repQ^{\alpha,\beta}_{p,r}$ with $1\leq
    r\leq p-1$, which are the simple $4p$-dimensional modules
    described~\bref{simple-list}.
  \item Typical modules $\repQ^{\alpha,\beta}_{s,r}$ with $1\leq s\leq
    p-1$, $1\leq r\leq p-1$ and $r\neq s$, described
    in~\bref{proj:typical}; each has four simple subquotients and
    dimension~$8p$.
  \item Atypical modules:
    \begin{itemize}
    \item The $12p$-dimensional $\repQ^{\alpha,\beta}_{p,0}$ modules
      described in~\bref{proj:po}, with $12$ simple subquotients each.

    \item The $12p$-dimensional modules
      $\repQ^{\alpha,\beta}_{p-1,p-1}$ described in~\bref{proj:pp},
      with $12$ simple subquotients.

    \item The $16p$-dimensional $\repQ^{\alpha,\beta}_{s,0}$ modules
      with $1\leq s\leq p-1$ described in~\bref{proj:so}, with $16$
      simple subquotients.

    \item The $16p$-dimensional modules $\repQ^{\alpha,\beta}_{s,s}$
      with $1\leq s\leq p-2$ described in~\bref{proj:ss}, with $16$
      simple subquotients.

    \item The $24p$-dimensional modules $\repQ^{\alpha,\beta}_{1,0}$
      described in~\bref{proj:10}, with $24$ simple subquotients.
    \end{itemize}
  \end{itemize}
\end{Thm}

We construct each projective module $\repQ$ explicitly by choosing a
basis and defining the $\algU$ action in that basis.  The set of basis
vectors is the union of the bases of all simple subquotients.  This
means that for each projective $\repQ$, we choose and fix a linear
space isomorphism
\begin{equation}\label{direct-sum}
  \mu:\repQ\to\brepQ=\bigoplus_{\iota}\repZ_{\iota},
\end{equation}
with the direct sum of all simple subquotients of~$\repQ$.  In
accordance with the direct sum decomposition, we also write
$\mu=\sum_{\iota}\mu_{\iota}$.

\subsection{The $\algU$ action and graphs}\label{action+graphs}
We explain how we construct the action of $\algU$ generators on
projective modules.  For a generator $A$, its action on $v\in\repQ$ is
given by
\begin{equation}\label{rho-all}
  \rho_{A}(v) = \rho^{(0)}_{A}(v) + \rho^{(1)}_{A}(v) + \rho^{(2)}_{A}(v),
\end{equation}
with the following ingredients.
\begin{enumerate}\addtocounter{enumi}{-1}
\item $\rho^{(0)}_{A}=\mu^{-1}\ccirc\bar{\rho}_{A}\ccirc\mu$, where
  $\bar{\rho}$ is the direct sum of $\algU$ actions on simple modules.

\item To define $\rho^{(1)}_A$, we recall that in~\bref{sec:Ext1-typ}
  and~\bref{sec:Ext1-atyp}, we chose a collection of maps $\xi_{A}$
  between simple $\algU$ modules, which we now write using a more
  detailed notation, as $\xi^{i,j}_A:\repZ_i\to\repZ_j$ (for some
  simple modules $\repZ_i$ and $\repZ_j$).  Then
  \begin{equation}\label{rho1}
    \rho^{(1)}_{A}
    =\mu^{-1}
    \ccirc
    \sum_{\iota,\kappa}c_{\kappa,\iota}\xi^{\kappa,\iota}_{A}
    \ccirc
    \mu_{\kappa},
  \end{equation}
  where $c_{\kappa,\iota}\in\oC$ are some coefficients, depending on a
  pair of simple subquotients in the projective module in
  question.\footnote{\textit{Not} on pairs of isomorphism classes of
    simple $\algU$ modules.}

  For a projective module $\repQ$, the $\rho^{(1)}_A$ maps can be
  represented as a directed graph with the set of vertices given by
  the $\repZ_{\iota}$ in~\eqref{direct-sum} and the set of edges
  corresponding to the nonzero products
  $c_{\kappa,\iota}\xi^{\kappa,\iota}_{A}$ (nonzero for at least one
  $A$).  The edge is directed from $\repZ_{\kappa}$ to
  $\repZ_{\iota}$.  We construct such graphs in what follows,
  decorating the edges with $c_{\kappa,\iota}$.

\item Finally,
  \begin{equation}\label{rho2}
    \rho^{(2)}_A =\mu^{-1} \ccirc
    \sum_{\iota,\kappa}\eta^{\kappa,\iota}_{A} \ccirc \mu_{\kappa},
  \end{equation}
  where $\eta^{\kappa,\iota}_{A}:\repZ_{\kappa}\to\repZ_{\iota}$ are
  linear maps (also linear in $A$) and $\repZ_{\iota}$ is a descendant
  but not a child of $\repZ_{\kappa}$ in the graph defined by the
  $c_{\mu,\nu}\xi^{\mu,\nu}_{A}$.  These maps are needed
  for~\eqref{rho-all} to be a $\algU$ action.
\end{enumerate}

Proving the existence of the projective cover $\repQ$ of a simple
module $\repZ_{*}$ amounts to finding the coefficients $c_{\mu,\nu}$
and the maps $\eta^{\kappa,\iota}_{A}$ such that the graph has a root
vertex given by $\repZ_{*}$ and Eq.~\eqref{rho-all} is a $\algU$
action (and the resultant module is maximal indecomposable).  We solve
for such $c_{\mu,\nu}$ and $\eta^{\kappa,\iota}_{A}$ in what follows.
The solution is not unique due to the freedom of taking linear
combinations of (the respective basis vectors in) isomorphic
subquotients, but the existence of a solution proves the existence of
the corresponding module; that it is maximal indecomposable is then
verified by inspection in each particular case.

\subsection{Constructing the projective modules}
We proceed with projective modules starting with the simple ones
in~\bref{proj:steinberg} and ending with those having $24$ simple
subquotients in~\bref{proj:10}.

\addtocounter{subsubsection}{-1}%
\subsubsection{} We need the obvious notion of \textit{level} in a
directed graph with a root vertex ${*}$.  That vertex is assigned
level~$1$, every child of ${*}$ has level~$2$, and so on (in the
graphs we are dealing with, this defines the level of each vertex
uniquely).

\subsubsection{Simple projective modules: $\myboldsymbol{s=p}$, \
  $\myboldsymbol{1\leq r\leq p-1}$}\label{proj:steinberg} These are
the simple (``Steinberg'') modules discussed in and after
Eq.~\eqref{stein-decomp}.  Each of them is also projective.

\subsubsection{Projective covers of typical simple modules with
  $\myboldsymbol{1\leq s\leq p-1}$,
  $\myboldsymbol{r\neq0,s}$}\label{proj:typical}
The ``bulk'' of the diagram in Fig.~\ref{module-pattern} yields
$4(p-1)(p-2)$ modules $\repQ^{\alpha,\beta}_{s,r}$, each of dimension
$8p$, labeled by
\begin{equation*}
  1\leq s\leq p-1,\quad 1\leq r\leq p-1,\quad r\neq s,\quad
  \alpha,\beta=\pm.
\end{equation*}
Their graphs (see~\bref{action+graphs}) are very simple:
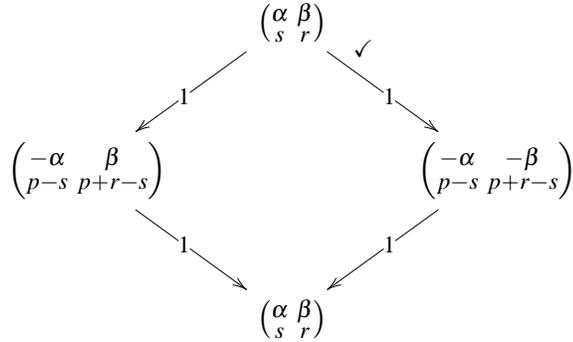
\begin{equation}\label{schem-proj}
  \xymatrix@=30pt{
    &&{\left(\!\begin{smallmatrix}\alpha &\beta\\s&r\end{smallmatrix}\!\right)\kern-4pt}
    \ar[dl]|1
    \ar[dr]^(.3){\checkmark}|1
    &\\
    &{\left(\!\begin{smallmatrix}-\alpha &\beta\\p-s&p+r-s\end{smallmatrix}\!\right)\kern-4pt}
    \ar[dr]|1
    &
    &
    {\left(\!\begin{smallmatrix}-\alpha &-\beta\\p-s&p+r-s\end{smallmatrix}\!\right)\kern-4pt}
    \ar[dl]|1
    \\
    &&
    {\left(\!\begin{smallmatrix}\alpha &\beta\\s&r\end{smallmatrix}\!\right)\kern-4pt}
    &
  }
\end{equation}
Here, $\left(\!\begin{smallmatrix}\alpha
    &\beta\\s&r\end{smallmatrix}\!\right)\kern-4pt=\repZ^{\alpha,\beta}_{s,r}$
are the simple subquotients,\footnote{We omit the uninformative
  $\repZ$ for conformity with similar, but much more complicated
  graphs in what follows, where extra symbols would complicate the
  picture even more.} and the arrows are the (basis elements in) the
corresponding $\Exti$ and the units on the arrows are the coefficients
$c^{\kappa,\iota}$ standing in front of these elements
(see~\eqref{rho1}).

We consider the arrow with a checkmark as an example.  From the
subquotients that it connects, we see that the relevant extension is
$f_{p - s, p + r - s}$ in~\eqref{basis-f}; the corresponding $\xi_A$
piece of the action of $\algU$ generators is therefore given
by~\eqref{xi-maps-f} times the coefficient~$1$.  For the opposite
arrow in the diagram, the same maps \eqref{xi-maps-f} should be used
with $\alpha\to-\alpha$, $s\to p-s$, and $r\to p+r-s$, and similarly
for the other arrows.

It remains to specify the $\eta_A$ piece of the action such that
Eqs.~\eqref{rho-all}--\eqref{rho2} define a $\algU$ action.  The
choice of the $\eta_A$ maps in a given basis is not unique (the
different solutions being mapped into one another by basis changes),
and we choose $\eta_A$ to be nonzero only for $A= B, E $, and only
acting on the top subquotient.  Its basis vectors~\eqref{four-sub},
now denoted as
\begin{equation*}
  \ZC{\alpha,s,\beta,r}_j\TOP,
  \quad
  \ZU{\alpha,s,\beta,r}_m\TOP,
  \quad
  \ZD{\alpha,s,\beta,r}_n\TOP,
  \quad
  \ZB{\alpha,s,\beta,r}_j\TOP,
\end{equation*}
are than mapped into the respective vectors in the bottom subquotient,
denoted as
\begin{equation*}
  \ZC{\alpha,s,\beta,r}_j\!\BOT,
  \quad
  \ZU{\alpha,s,\beta,r}_m\BOT,
  \quad
  \ZD{\alpha,s,\beta,r}_n\BOT,
  \quad
  \ZB{\alpha,s,\beta,r}_j\!\BOT.
\end{equation*}
The nonzero maps $\eta_{B}$ and $\eta_{E}$ are given by
\begin{align*}
  \eta_{B} \ZC{\alpha, s, \beta, r}_{j}\TOP &=
  \alpha \Gamma_{s, r, j} \ZD{\alpha, s, \beta, r}_{j - 1}\BOT -
  \alpha \beta \qint{r - s}\Gamma_{s, r, j} \ZU{\alpha, s, \beta,
    r}_{j}\BOT,
  \\
  \eta_{B} \ZU{\alpha, s, \beta, r}_{m}\TOP &=
  \alpha \qint{s} \Gamma_{s, r, m} \ZB{\alpha, s, \beta, r}_{m -
    1}\BOT,
  \\
  \eta_{B} \ZD{\alpha, s, \beta, r}_{n}\TOP &= \alpha \beta \qint{s}
  \qint{r - s}\Gamma_{s, r, n + 1} \ZB{\alpha, s, \beta, r}_{n}\BOT,
  \\
  \eta_{E} \ZU{\alpha, s, \beta, r}_{m}\TOP &=
  \qint{s + 1} \ZU{\alpha, s, \beta, r}_{m - 1}\BOT,
  \\
  \eta_{E} \ZC{\alpha, s, \beta, r}_{j}\TOP &=
  \qint{s} \ZC{\alpha, s, \beta, r}_{j - 1}\BOT,
  \\
  \eta_{E} \ZD{\alpha, s, \beta, r}_{n}\TOP &=
  \qint{s - 1} \ZD{\alpha, s, \beta, r}_{n - 1}\BOT,
  \\
  \eta_{E} \ZB{\alpha, s, \beta, r}_{j}\TOP &=
  \qint{s} \ZB{\alpha, s, \beta, r}_{j - 1}\BOT,
\end{align*}
where
\begin{equation*}
  \Gamma_{s, r, m} =
  -\ffrac{\qint{m - 1}}{\qint{s}}
  + \half\Bigl(-(\q^r + \q^{-r}) \ffrac{\qint{s - 1}}{\qint{s}}
  + \qint{r} \ffrac{\qint{2 s - 1} - 3}{\qint{s}^2}\Bigr)
  \ffrac{\qint{m}}{\qint{r - s}}.
\end{equation*}
Direct calculation shows that with the coefficients $1$
in~\eqref{schem-proj}, all relations for the $\algU$ generators are
satisfied.

\subsubsection{$\Zpo$: Projective cover of
  $\boldsymbol{\repZ}^{\myboldsymbol{\alpha,\beta}}_{\myboldsymbol{p,0}}$}
\label{proj:po}
The projective cover $\repQ^{\alpha,\beta}_{p,0}$ of
$\repZ^{\alpha,\beta}_{p,0}$ (denoted by $\zpo$ in
Fig.~\ref{module-pattern}) has $12$ simple subquotients and is
$12p$-dimensional.  The corresponding graph is shown in
Fig.~\ref{fig:proj-po}.
\begin{figure}[tbhp]
  \centering\footnotesize
  \begin{gather*}
    \xymatrix@C=25pt@R=40pt{\poiApBZ="poiApBZ"
      \\
      \poiieemAIBZ="poiieemAIBZ"&& \poiiefIApmimBZ="poiiefIApmimBZ"&&
      \poiiffmAImBZ="poiiffmAImBZ"
      \\
      \poiiiCcApBZ="poiiiCcApBZ"&& \poiiiEImAIBI="poiiiEImAIBI"&&
      \poiiiEfcApmiBpmi="poiiiEfcApmiBpmi"&&
      \poiiiFImAImBI="poiiiFImAImBI"
      \\
      \poivEEmAIBZ="poivEEmAIBZ"&& \poivBBApmimBZ="poivBBApmimBZ"&&
      \poivFFmAImBZ="poivFFmAImBZ"
      \\
      \povbotApBZ="povbotApBZ"
      \ar@{->}|(.45){\arlabel{\kern5pt-\beta}} "poiApBZ";
      "poiieemAIBZ" \ar@{->}|(.45){\arlabel{\kern5pt1}} "poiApBZ";
      "poiiefIApmimBZ" \ar@{->}|(.45){\arlabel{\kern5pt\alpha \beta}}
      "poiApBZ"; "poiiffmAImBZ" \ar@{->}|(.45){\arlabel{\kern5pt\alpha
          \beta}} "poiieemAIBZ"; "poiiiCcApBZ"
      \ar@{->}|(.70){\arlabel{\kern5pt-\beta}} "poiiefIApmimBZ";
      "poiiiCcApBZ" \ar@{->}|(.38){\arlabel{\kern5pt-\beta}}
      "poiiffmAImBZ"; "poiiiCcApBZ"
      \ar@{->}|(.40){\arlabel{\kern5pt-\alpha \beta}} "poiieemAIBZ";
      "poiiiEImAIBI" \ar@{->}|(.35){\arlabel{\kern5pt-\alpha}}
      "poiiefIApmimBZ"; "poiiiEImAIBI"
      \ar@{->}|(.70){\arlabel{\kern5pt1}} "poiieemAIBZ";
      "poiiiEfcApmiBpmi" \ar@{->}|(.35){\arlabel{\kern5pt-\alpha}}
      "poiiffmAImBZ"; "poiiiEfcApmiBpmi"
      \ar@{->}|(.38){\arlabel{\kern5pt-\alpha}} "poiiefIApmimBZ";
      "poiiiFImAImBI" \ar@{->}|(.45){\arlabel{\kern5pt-\beta}}
      "poiiffmAImBZ"; "poiiiFImAImBI"
      \ar@{->}|(.35){\arlabel{\kern5pt1}} "poiiiCcApBZ";
      "poivBBApmimBZ" \ar@{->}|(.35){\arlabel{\kern5pt1}}
      "poiiiEImAIBI"; "poivBBApmimBZ"
      \ar@{->}|(.35){\arlabel{\kern5pt1}} "poiiiFImAImBI";
      "poivBBApmimBZ" \ar@{->}|(.50){\arlabel{\kern5pt-\alpha \beta}}
      "poiiiCcApBZ"; "poivEEmAIBZ"
      \ar@{->}|(.65){\arlabel{\kern5pt-1}} "poiiiEImAIBI";
      "poivEEmAIBZ" \ar@{->}|(.35){\arlabel{\kern5pt1}}
      "poiiiEfcApmiBpmi"; "poivEEmAIBZ"
      \ar@{->}|(.63){\arlabel{\kern5pt-\alpha \beta}} "poiiiCcApBZ";
      "poivFFmAImBZ" \ar@{->}|(.35){\arlabel{\kern5pt1}}
      "poiiiEfcApmiBpmi"; "poivFFmAImBZ"
      \ar@{->}|(.45){\arlabel{\kern5pt-1}} "poiiiFImAImBI";
      "poivFFmAImBZ" \ar@{->}|(.45){\arlabel{\kern5pt-1}}
      "poivBBApmimBZ"; "povbotApBZ"
      \ar@{->}|(.45){\arlabel{\kern5pt1}} "poivEEmAIBZ"; "povbotApBZ"
      \ar@{->}|(.45){\arlabel{\kern5pt1}} "poivFFmAImBZ";
      "povbotApBZ"}
  \end{gather*}
  \caption[Graph of the projective module
  $\repQ^{\alpha,\beta}_{p,0}$.]{\small\ \ Graph of the projective
    module $\repQ^{\alpha,\beta}_{p,0}$}
  \label{fig:proj-po}
\end{figure}
Each simple subquotient $\repZ^{\alpha,\beta}_{s,r}$ is identified by
its parameters $\left(\!\begin{smallmatrix}\alpha &\beta
    \\s&r\end{smallmatrix}\!\right)$ and is in addition labeled by
$\mylabel{}{\ell}{n}$, where $\ell$ is the level and $n$ consecutively
labels subquotients within each level.

As previously, each link $\repZ_1\to\repZ_2$ corresponds to the basis
element in $\Exti(\repZ_1,\repZ_2)$ (which, we recall, is
one-dimensional) times the coefficient placed on the link.  For
example, consider the level-2-to-level-3 link
$\myxy{\poiiefIApmimBZ}\xrightarrow{-\alpha}\myxy{\poiiiFImAImBI}$,
which of course stands for
$\repZ^{\alpha,-\beta}_{p-1,0}\xrightarrow{-\alpha}\repZ^{-\alpha,-\beta}_{1,1}$.
It corresponds to the element $\bar{e}_{1}$ in~\bref{sec:Ext1-atyp}
(the last in the list) in accordance with the somewhat truncated
notation used there, which ignores the upper indices.  The coefficient
$-\alpha$ on the link means that $ E $ and $ B $ map from
$\repZ^{\alpha,-\beta}_{p-1,0}$ to $\repZ^{-\alpha,-\beta}_{1,1}$~as
\begin{align*}
  \xi_{E}&:\ZC{\alpha,p-1,-\beta,0}_0\mapsto
  -\alpha\cdot\ZC{-\alpha, 1, -\beta, 1}_{0},\\
  \xi_{E}&:\ZB{\alpha,p-1,-\beta,0}_0\mapsto
  -\alpha\cdot\qint{-2}\ZB{-\alpha, 1, -\beta, 1}_{1},\\
  \xi_{B}&:\ZC{\alpha,p-1,-\beta,0}_{0}\mapsto
  -\alpha\cdot(-\alpha)\ZB{-\alpha, 1, -\beta, 1}_{1},
\end{align*}
and these are the only maps by the $\algU$ generators between the
simple subquotients $\repZ^{\alpha,-\beta}_{p-1,0}$ and
$\repZ^{-\alpha,-\beta}_{1,1}$ of $\repQ^{\alpha,\beta}_{p,0}$.

The same $\repZ^{\alpha,-\beta}_{p-1,0}$ is also linked in the graph
to $\repZ^{-\alpha,\beta}_{1,1}$, with a link decorated by $-\alpha$;
this is $\bar{f}_{1}$ from the list in~\bref{sec:Ext1-atyp}, and the
nonzero maps by $\algU$ generators are therefore given by
\begin{align*}
  \xi_{F}&:\ZC{\alpha,p-1,-\beta,0}_{p-2}\mapsto -\alpha\cdot\ZC{-\alpha, 1, \beta, 1}_0,\\
  \xi_{F}&:\ZB{\alpha,p-1,-\beta,0}_{p-3}\mapsto
  -\alpha\cdot(-)\ZB{-\alpha, 1, \beta, 0}_0.
\end{align*}

In addition to the graph, we have to specify the $\eta_A$ maps such
that Eqs.~\eqref{rho-all}--\eqref{rho2} define a $\algU$ action.  This
part of the action of $\algU$ generators can be chosen as follows:
\begin{align*}
  \eta_E\ZB{\alpha,p,\beta,0}_{n,\mylabel{}{1}{0}}&= -\alpha
  \ZB{\alpha,p,\beta,0}_{n-1,\mylabel{Cc}{3}{0}} -\beta
  \ZC{\alpha,p\!-\!1,\beta,p\!-\!1}_{n-1,\mylabel{Efc}{3}{2}},
  \\
  \eta_B\ZC{\alpha,p,\beta,0}_{n,\mylabel{}{1}{0}}&= -\beta \qint{n}
  \ZB{\alpha,p,\beta,0}_{n-1,\mylabel{bot}{5}{0}}
  +\qint{n\!+\!1}\ZB{\alpha,p,\beta,0}_{n-1,\mylabel{Cc}{3}{0}}
  \\
  &\quad{}+\alpha \beta \qint{n\!+\!1} \ZC{\alpha,p\!-\!1,\beta,p\!-\!1}_{n-1,\mylabel{Efc}{3}{2}},
  \\
  \eta_B\ZB{\alpha,p,\beta,0}_{n,\mylabel{}{1}{0}}&= \alpha \beta
  \qint{n+2} \ZB{\alpha,p\!-\!1,\beta,p\!-\!1}_{n,\mylabel{Efc}{3}{2}},
  \\
  \eta_C\ZC{\alpha,p,\beta,0}_{n,\mylabel{}{1}{0}}&= -\alpha \delta_{n,p-1} \ZC{-\alpha,1,\beta,1}_{0,\mylabel{E1}{3}{1}},
  \\
  \eta_C\ZB{\alpha,p,\beta,0}_{n,\mylabel{}{1}{0}}&= \alpha \delta_{n,p-2} \ZB{-\alpha,1,\beta,1}_{0,\mylabel{E1}{3}{1}},
  \\
  \eta_E\ZC{\alpha,p\!-\!1,-\beta,0}_{n,\mylabel{ef1}{2}{1}}&=
  -\alpha \ZC{\alpha,p\!-\!1,-\beta,0}_{n-1,\mylabel{BB}{4}{1}},
  \\
  \eta_E\ZB{\alpha,p\!-\!1,-\beta,0}_{n,\mylabel{ef1}{2}{1}}&=
  -\alpha \qint{2}\ZB{\alpha,p\!-\!1,-\beta,0}_{n-1,\mylabel{BB}{4}{1}}
  \\
  \eta_B\ZC{\alpha,p\!-\!1,-\beta,0}_{n,\mylabel{ef1}{2}{1}}&=
  \qint{n\!+\!1}\ZB{\alpha,p\!-\!1,-\beta,0}_{n-1,\mylabel{BB}{4}{1}}
  -\qint{n}\ZC{\alpha,p,\beta,0}_{n,\mylabel{bot}{5}{0}},
  \\
  \eta_B\ZB{\alpha,p\!-\!1,-\beta,0}_{n,\mylabel{ef1}{2}{1}}&=
  \qint{n\!+\!1}\ZB{\alpha,p,\beta,0}_{n,\mylabel{bot}{5}{0}},
  \\
  \eta_E\ZB{\alpha,p,\beta,0}_{n,\mylabel{Cc}{3}{0}}&=-\alpha \beta
  \ZB{\alpha,p,\beta,0}_{n-1,\mylabel{bot}{5}{0}},
  \\
  \eta_B\ZC{\alpha,p,\beta,0}_{n,\mylabel{Cc}{3}{0}}&=\beta
  \qint{n\!+\!1} \ZB{\alpha,p,\beta,0}_{n-1,\mylabel{bot}{5}{0}},
  \\
  \eta_B\ZC{-\alpha,1,\beta,1}_{n,\mylabel{E1}{3}{1}}&= \alpha \delta_{0,n}\ZC{\alpha,p,\beta,0}_{p-1,\mylabel{bot}{5}{0}},
  \\
  \eta_B\ZB{-\alpha,1,\beta,1}_{n,\mylabel{E1}{3}{1}}&= \alpha \delta_{0,n}\ZB{\alpha,p,\beta,0}_{p-2,\mylabel{bot}{5}{0}},
  \\
  \eta_E\ZC{\alpha,p\!-\!1,\beta,p\!-\!1}_{n,\mylabel{Efc}{3}{2}}&=
  \ZB{\alpha,p,\beta,0}_{n-1,\mylabel{bot}{5}{0}}.
\end{align*}
With these $\eta_A$ and with the $c_{\iota,\kappa}$ read off from the
graph, Eqs.~\eqref{rho-all}--\eqref{rho2} define an $\algU$ action, as
can be verified directly.  Inspection shows that the resulting module
is indecomposable and maximal.

\subsubsection{$\Zpp$: Projective cover of
  $\boldsymbol{\repZ}^{\myboldsymbol{\alpha,\beta}}_{\myboldsymbol{p-1,p-1}}$}
\label{proj:pp}
The projective cover $\repQ^{\alpha,\beta}_{p-1,p-1}$ of
$\repZ^{\alpha,\beta}_{p-1,p-1}$ (denoted by $\zpp$
in Fig.~\ref{module-pattern}) also has $12$ subquotients and is
$12p$-dimensional.  Its graph is shown in Fig.~\ref{fig:projpipi},
with the same notation and conventions as for the preceding projective
module.
\begin{figure}[tbhp]
  \centering\footnotesize
  \begin{gather*}
    \xymatrix@C=25pt@R=39pt{\ppiApmiBpmi="ppiApmiBpmi"
      \\
      \ppiiccApmiimBpmii="ppiiccApmiimBpmii"&&
      \ppiieemAIBZ="ppiieemAIBZ"&& \ppiiffmAImBZ="ppiiffmAImBZ"
      \\
      \ppiiiEemAIImBZ="ppiiiEemAIImBZ"&&
      \ppiiiEfbApmiBpmi="ppiiiEfbApmiBpmi"&&
      \ppiiiFfmAIIBZ="ppiiiFfmAIIBZ"&& \ppiiiBbApBZ="ppiiiBbApBZ"
      \\
      \ppivCCApmiimBpmii="ppivCCApmiimBpmii"&&
      \ppivEEmAIBZ="ppivEEmAIBZ"&& \ppivFFmAImBZ="ppivFFmAImBZ"
      \\
      \ppvbotApmiBpmi="ppvbotApmiBpmi"
      \ar@{->}|(.45){\arlabel{\kern5pt\beta}}
      "ppiApmiBpmi";"ppiiccApmiimBpmii"
      \ar@{->}|(.45){\arlabel{\kern5pt1}} "ppiApmiBpmi"; "ppiieemAIBZ"
      \ar@{->}|(.45){\arlabel{\kern5pt-\alpha}} "ppiApmiBpmi";
      "ppiiffmAImBZ" \ar@{->}|(.65){\arlabel{\kern5pt-\alpha}}
      "ppiieemAIBZ"; "ppiiiBbApBZ" \ar@{->}|(.45){\arlabel{\kern5pt1}}
      "ppiiffmAImBZ"; "ppiiiBbApBZ"
      \ar@{->}|(.45){\arlabel{\kern5pt1}}
      "ppiiccApmiimBpmii";"ppiiiEemAIImBZ"
      \ar@{->}|(.65){\arlabel{\kern5pt\beta}}
      "ppiieemAIBZ";"ppiiiEemAIImBZ"
      \ar@{->}|(.65){\arlabel{\kern5pt\alpha}} "ppiiccApmiimBpmii";
      "ppiiiEfbApmiBpmi" \ar@{->}|(.67){\arlabel{\kern5pt1}}
      "ppiieemAIBZ"; "ppiiiEfbApmiBpmi"
      \ar@{->}|(.50){\arlabel{\kern5pt-\alpha}}
      "ppiiffmAImBZ";"ppiiiEfbApmiBpmi"
      \ar@{->}|(.45){\arlabel{\kern5pt1}} "ppiiccApmiimBpmii";
      "ppiiiFfmAIIBZ" \ar@{->}|(.38){\arlabel{\kern5pt\alpha \beta}}
      "ppiiffmAImBZ";"ppiiiFfmAIIBZ"
      \ar@{->}|(.45){\arlabel{\kern5pt1}} "ppiiiEemAIImBZ";
      "ppivCCApmiimBpmii" \ar@{->}|(.35){\arlabel{\kern5pt\beta}}
      "ppiiiEfbApmiBpmi"; "ppivCCApmiimBpmii"
      \ar@{->}|(.42){\arlabel{\kern5pt1}} "ppiiiFfmAIIBZ";
      "ppivCCApmiimBpmii" \ar@{->}|(.35){\arlabel{\kern5pt1}}
      "ppiiiBbApBZ"; "ppivEEmAIBZ" \ar@{->}|(.35){\arlabel{\kern5pt1}}
      "ppiiiEemAIImBZ"; "ppivEEmAIBZ"
      \ar@{->}|(.68){\arlabel{\kern5pt\alpha}} "ppiiiEfbApmiBpmi";
      "ppivEEmAIBZ" \ar@{->}|(.45){\arlabel{\kern5pt1}} "ppiiiBbApBZ";
      "ppivFFmAImBZ" \ar@{->}|(.50){\arlabel{\kern5pt\alpha}}
      "ppiiiEfbApmiBpmi"; "ppivFFmAImBZ"
      \ar@{->}|(.67){\arlabel{\kern5pt1}} "ppiiiFfmAIIBZ";
      "ppivFFmAImBZ" \ar@{->}|(.45){\arlabel{\kern5pt1}}
      "ppivCCApmiimBpmii"; "ppvbotApmiBpmi"
      \ar@{->}|(.45){\arlabel{\kern5pt1}} "ppivEEmAIBZ";
      "ppvbotApmiBpmi" \ar@{->}|(.45){\arlabel{\kern5pt1}}
      "ppivFFmAImBZ"; "ppvbotApmiBpmi"}
  \end{gather*}
  \caption[Graph of the projective module
  $\repQ^{\alpha,\beta}_{p-1,p-1}$.]{\small\ \ Graph of the projective
    module $\repQ^{\alpha,\beta}_{p-1,p-1}$.}
  \label{fig:projpipi}
\end{figure}

The $\eta$ piece of the action by the $\algU$ generators needed
in~\eqref{rho-all} is as follows.  On the basis vectors of the top
subquotient, we have
\begin{align*}
  \eta_F\ZC{\alpha,p\!-\!1,\beta,p\!-\!1}_{n,\mylabel{}{1}{0}}&=-\qint{2} \delta _{n,p-2}
  \ZC{-\alpha,1,\beta,0}_{0,\mylabel{EE}{4}{1}} 
 ,\\
  \eta_E\ZC{\alpha,p\!-\!1,\beta,p\!-\!1}_{n,\mylabel{}{1}{0}}&=-\alpha \ZB{\alpha,p,\beta,0}_{n-1,\mylabel{Bb}{3}{3}} -\qint{2}\ZC{\alpha,p\!-\!1,\beta,p\!-\!1}_{n-1,\mylabel{bot}{5}{0}}
  \\
  &\quad{}+\ZC{\alpha,p\!-\!1,\beta,p\!-\!1}_{n-1,\mylabel{Efb}{3}{1}}-\qint{2}
  \delta_{0,n} \ZC{-\alpha,1,-\beta,0}_{0,\mylabel{FF}{4}{2}} 
 ,\\
  \eta_B\ZC{\alpha,p\!-\!1,\beta,p\!-\!1}_{n,\mylabel{}{1}{0}}&=-\alpha \beta \qint{n\!+\!1} \ZC{\alpha,p\!-\!2,-\beta,p\!-\!2}_{n-1,\mylabel{CC}{4}{0}}\\
  &\quad{}+\alpha \qint{n} \ZB{\alpha,p\!-\!1,\beta,p\!-\!1}_{n,\mylabel{Efb}{3}{1}}-\alpha \beta
  \delta_{0,n} \ZC{-\alpha,2,\beta,0}_{1,\mylabel{Ff}{3}{2}},\\
  \eta_B\ZB{\alpha,p\!-\!1,\beta,p\!-\!1}_{n,\mylabel{}{1}{0}}&=\alpha \beta \qint{n\!+\!1}
  \ZB{\alpha,p\!-\!2,-\beta,p\!-\!2}_{n-1,\mylabel{CC}{4}{0}} -\alpha
  \beta \delta _{0,n} \ZB{-\alpha,2,\beta,0}_{0,\mylabel{Ff}{3}{2}},\\
  \eta_C\ZB{\alpha,p\!-\!1,\beta,p\!-\!1}_{n,\mylabel{}{1}{0}}&=\alpha \beta \ZC{\alpha,p\!-\!1,\beta,p\!-\!1}_{n,\mylabel{bot}{5}{0}} -\beta
  \qint{2}\delta _{n,p-1} \ZC{-\alpha,1,\beta,0}_{0,\mylabel{EE}{4}{1}}.  
\end{align*}
The $\eta_{A}$ maps from level-$2$ vectors are
\begin{align*}
  \eta_E\ZC{\alpha,p\!-\!2,-\beta,p\!-\!2}_{n,\mylabel{cc}{2}{0}}&=
  \qint{2}\ZC{\alpha,p\!-\!2,-\beta,p\!-\!2}_{n-1,\mylabel{CC}{4}{0}}
 ,\\
  \eta_E\ZB{\alpha,p\!-\!2,-\beta,p\!-\!2}_{n,\mylabel{cc}{2}{0}}&=
  \ZB{\alpha,p\!-\!2,-\beta,p\!-\!2}_{n-1,\mylabel{CC}{4}{0}}
 ,\\
  \eta_B\ZC{\alpha,p\!-\!2,-\beta,p\!-\!2}_{n,\mylabel{cc}{2}{0}}&=
  \alpha \qint{n\!+\!1} \ZB{\alpha,p\!-\!2,-\beta,p\!-\!2}_{n,\mylabel{CC}{4}{0}}
 ,\\
  \eta_B\ZC{-\alpha,1,\beta,0}_{n,\mylabel{ee}{2}{1}}&= \alpha
  \qint{2}\delta _{0,n} \ZB{\alpha,p\!-\!1,\beta,p\!-\!1}_{p-1,\mylabel{bot}{5}{0}},
\end{align*}
and those from level 3,
\begin{align*}
  \eta_E\ZB{\alpha,p,\beta,0}_{n,\mylabel{Bb}{3}{3}}&=\ZC{\alpha,p\!-\!1,\beta,p\!-\!1}_{n-1,\mylabel{bot}{5}{0}}
 ,\\
  \eta_B\ZC{\alpha,p,\beta,0}_{n,\mylabel{Bb}{3}{3}}&=-\alpha
  \qint{n\!+\!1} \ZC{\alpha,p\!-\!1,\beta,p\!-\!1}_{n-1,\mylabel{bot}{5}{0}}
 ,\\
  \eta_B\ZB{\alpha,p,\beta,0}_{n,\mylabel{Bb}{3}{3}}&=-\alpha
  \qint{n+2} \ZB{\alpha,p\!-\!1,\beta,p\!-\!1}_{n,\mylabel{bot}{5}{0}}
 ,\\
  \eta_E\ZC{\alpha,p\!-\!1,\beta,p\!-\!1}_{n,\mylabel{Efb}{3}{1}}&=
  \alpha \ZC{\alpha,p\!-\!1,\beta,p\!-\!1}_{n-1,\mylabel{bot}{5}{0}}
 ,\\
  \eta_B\ZC{\alpha,p\!-\!1,\beta,p\!-\!1}_{n,\mylabel{Efb}{3}{1}}&=
  \qint{n}\ZB{\alpha,p\!-\!1,\beta,p\!-\!1}_{n,\mylabel{bot}{5}{0}}
 ,\\
  \eta_C\ZC{-\alpha,2,\beta,0}_{n,\mylabel{Ff}{3}{2}}&= \delta_{1,n}\ZC{\alpha,p\!-\!1,\beta,p\!-\!1}_{0,\mylabel{bot}{5}{0}}
 ,\\
  \eta_C\ZB{-\alpha,2,\beta,0}_{n,\mylabel{Ff}{3}{2}}&= -\delta_{0,n}\ZB{\alpha,p\!-\!1,\beta,p\!-\!1}_{0,\mylabel{bot}{5}{0}}.
\end{align*}
Direct calculation shows that the above $\eta_A$ and the
$c_{\iota,\kappa}$ read off from the graph ensure that
Eqs.~\eqref{rho-all}--\eqref{rho2} define a projective $\algU$
module.\enlargethispage{\baselineskip}

\subsubsection{$\Zso$: Projective cover of
  $\boldsymbol{\repZ}^{\myboldsymbol{\alpha,\beta}}_{\myboldsymbol{s,0}}$
  for $\myboldsymbol{2\leq s\leq p-1}$.}\label{proj:so}
The projective cover $\repQ^{\alpha,\beta}_{s,0}$ of
$\repZ^{\alpha,\beta}_{s,0}$ (any of the $\zso$
in Fig.~\ref{module-pattern}) with $2\leq s\leq p-1$ has $16$ simple
subquotients and is $16p$-dimensional.  Its graph is shown in
Fig.~\ref{fig:proj-so}.
\afterpage{%
  \begin{landscape}
    \begin{figure}[tbh]
      \centering\footnotesize
      \begin{gather*}
        \xymatrix@R=25pt@C=20pt{
          &&&&&\soiAsBZ="soiAsBZ"&&&&&\\
          && \soiietmApmsmBpms="soiietmApmsmBpms"&&
          \soiiftmApmsBpms="soiiftmApmsBpms"&&
          \soiitlAsmimBZ="soiitlAsmimBZ"&&
          \soiitrAsimBZ="soiitrAsimBZ"&&\\
          &&&&&&&&&&\\
          \\
          \soiiiBtAsBZ="soiiiBtAsBZ"&&
          \soiiiElmApimsBpims="soiiiElmApimsBpims"&&
          \soiiiErmApmsmiBpmsmi="soiiiErmApmsmiBpmsmi"&&
          \soiiiFlmApimsmBpims="soiiiFlmApimsmBpims"&&
          \soiiiFrmApmsmimBpmsmi="soiiiFrmApmsmimBpmsmi"&&
          \soiiiTbAsBZ="soiiiTbAsBZ"\\
          \\
          &&&&&&&&&&\\
          && \soivBLAsmimBZ="soivBLAsmimBZ"&&
          \soivBRAsimBZ="soivBRAsimBZ"&&
          \soivEBmApmsmBpms="soivEBmApmsmBpms"&&
          \soivFBmApmsBpms="soivFBmApmsBpms"&&\\
          &&&&& \sovbotAsBZ="sovbotAsBZ"&&&&&
          \ar@{->}|(.45){\arlabel{\kern5pt\qint{s-1} \qint{s}^3}}
          "soiAsBZ"; "soiietmApmsmBpms"
          \ar@{->}|(.45){\arlabel{\kern5pt-\qint{s}^3}} "soiAsBZ";
          "soiiftmApmsBpms"
          \ar@{->}|(.45){\arlabel{\kern5pt\qint{s}^4}} "soiAsBZ";
          "soiitlAsmimBZ" \ar@{->}|(.45){\arlabel{\kern5pt\qint{s-1}
              \qint{s}^2}} "soiAsBZ"; "soiitrAsimBZ"
          \ar@{->}|(.45){\arlabel{\kern5pt-1}} "soiietmApmsmBpms";
          "soiiiBtAsBZ" \ar@{->}|(.70){\arlabel{\kern5pt\qint{s-1}}}
          "soiiftmApmsBpms"; "soiiiBtAsBZ"
          \ar@{->}|(.71){\arlabel{\kern5pt\alpha \beta}}
          "soiitrAsimBZ"; "soiiiBtAsBZ"
          \ar@{->}|(.35){\arlabel{\kern5pt-\qint{s}}}
          "soiietmApmsmBpms"; "soiiiElmApimsBpims"
          \ar@{->}|(.20){\arlabel{\kern5pt-\qint{s-1}}}
          "soiitlAsmimBZ"; "soiiiElmApimsBpims"
          \ar@{->}|(.47){\arlabel{\kern5pt-\qint{s-1}}}
          "soiietmApmsmBpms"; "soiiiErmApmsmiBpmsmi"
          \ar@{->}|(.80){\arlabel{\kern5pt\qint{s-1} \qint{s}}}
          "soiitrAsimBZ"; "soiiiErmApmsmiBpmsmi"
          \ar@{->}|(.53){\arlabel{\kern5pt-\qint{s}}}
          "soiiftmApmsBpms"; "soiiiFlmApimsmBpims"
          \ar@{->}|(.75){\arlabel{\kern5pt1}} "soiitlAsmimBZ";
          "soiiiFlmApimsmBpims"
          \ar@{->}|(.75){\arlabel{\kern5pt-\qint{s-1}^2}}
          "soiiftmApmsBpms"; "soiiiFrmApmsmimBpmsmi"
          \ar@{->}|(.30){\arlabel{\kern5pt\qint{s-1} \qint{s}}}
          "soiitrAsimBZ"; "soiiiFrmApmsmimBpmsmi"
          \ar@{->}|(.70){\arlabel{\kern5pt1}} "soiitlAsmimBZ";
          "soiiiTbAsBZ" \ar@{->}|(.45){\arlabel{\kern5pt\qint{s}}}
          "soiitrAsimBZ"; "soiiiTbAsBZ"
          \ar@{->}|(.45){\arlabel{\kern5pt-\qint{s}}} "soiiiBtAsBZ";
          "soivBLAsmimBZ" \ar@{->}|(.72){\arlabel{\kern5pt-1}}
          "soiiiElmApimsBpims"; "soivBLAsmimBZ"
          \ar@{->}|(.25){\arlabel{\kern5pt\qint{s-1}}}
          "soiiiFlmApimsmBpims"; "soivBLAsmimBZ"
          \ar@{->}|(.71){\arlabel{\kern5pt\alpha \beta}}
          "soiiiTbAsBZ"; "soivBLAsmimBZ"
          \ar@{->}|(.62){\arlabel{\kern5pt-\qint{s-1}}} "soiiiBtAsBZ";
          "soivBRAsimBZ" \ar@{->}|(.25){\arlabel{\kern5pt1}}
          "soiiiErmApmsmiBpmsmi"; "soivBRAsimBZ"
          \ar@{->}|(.83){\arlabel{\kern5pt1}} "soiiiFrmApmsmimBpmsmi";
          "soivBRAsimBZ" \ar@{->}|(.25){\arlabel{\kern5pt1}}
          "soiiiElmApimsBpims"; "soivEBmApmsmBpms"
          \ar@{->}|(.46){\arlabel{\kern5pt1}} "soiiiErmApmsmiBpmsmi";
          "soivEBmApmsmBpms"
          \ar@{->}|(.80){\arlabel{\kern5pt\qint{s-1}}} "soiiiTbAsBZ";
          "soivEBmApmsmBpms"
          \ar@{->}|(.48){\arlabel{\kern5pt\qint{s-1}}}
          "soiiiFlmApimsmBpims"; "soivFBmApmsBpms"
          \ar@{->}|(.72){\arlabel{\kern5pt1}} "soiiiFrmApmsmimBpmsmi";
          "soivFBmApmsBpms"
          \ar@{->}|(.45){\arlabel{\kern5pt\qint{s-1}}} "soiiiTbAsBZ";
          "soivFBmApmsBpms" \ar@{->}|(.45){\arlabel{\kern5pt1}}
          "soivBLAsmimBZ"; "sovbotAsBZ"
          \ar@{->}|(.45){\arlabel{\kern5pt1}} "soivBRAsimBZ";
          "sovbotAsBZ" \ar@{->}|(.45){\arlabel{\kern5pt1}}
          "soivEBmApmsmBpms"; "sovbotAsBZ"
          \ar@{->}|(.45){\arlabel{\kern5pt1}} "soivFBmApmsBpms";
          "sovbotAsBZ" }
      \end{gather*}
      \caption[Subquotients of the projective module
      $\repQ^{\alpha,\beta}_{s,0}$.]{\small\ \ Subquotients of the
        projective module $\repQ^{\alpha,\beta}_{s,0}$ for $2\leq
        s\leq p-1$.}
      \label{fig:proj-so}
    \end{figure}
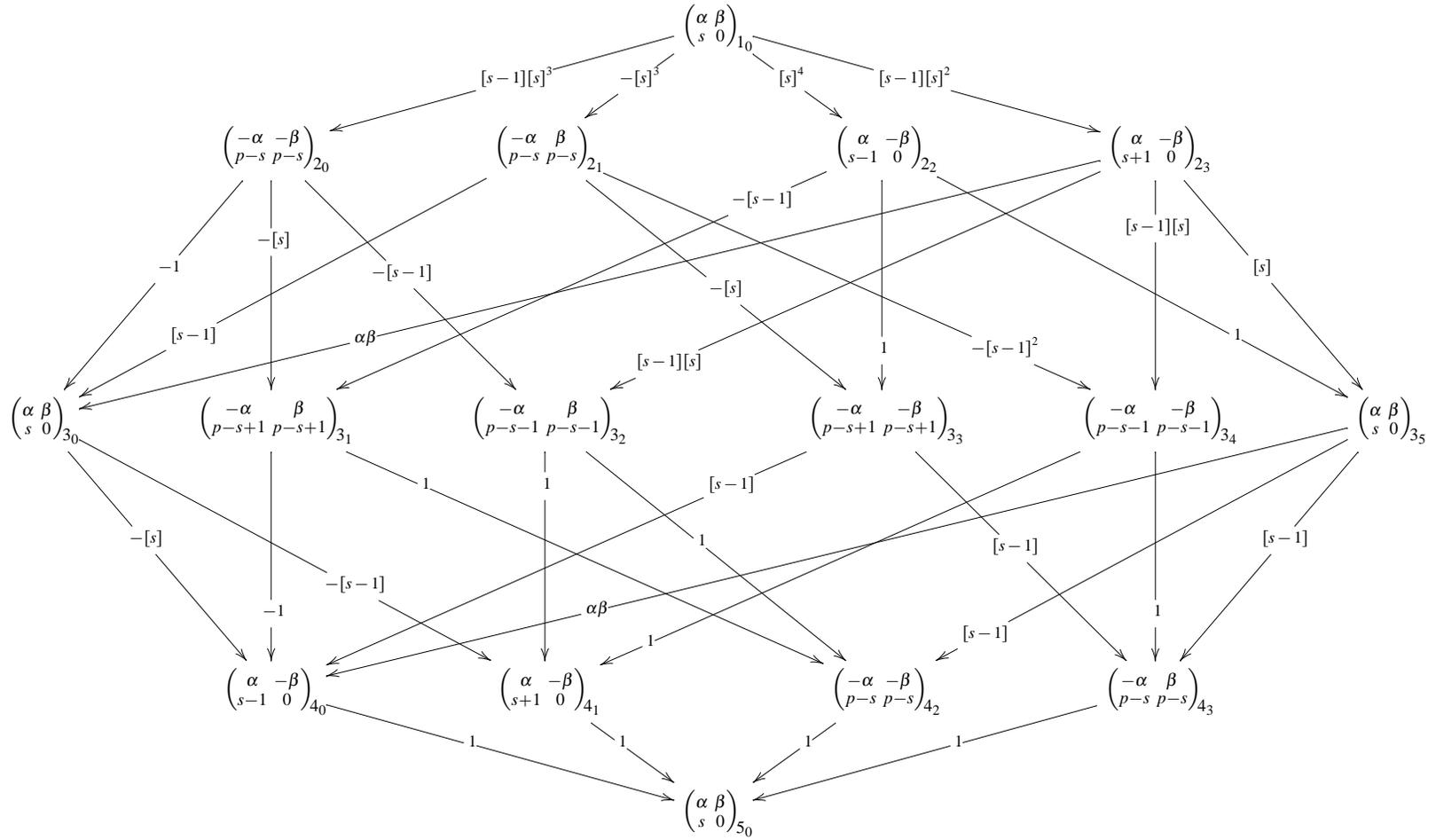
  \end{landscape}%
}%
A feature not encountered in the previous graphs is the occurrence of
links from a given subquotient leading to isomorphic subquotients on
the next level.  Two (or more) such links
\begin{equation*}
  \xymatrix@R=25pt@C=0pt{
    &{\left(\!\begin{smallmatrix}\alpha &\beta
          \\s&r\end{smallmatrix}\!\right)_{\mylabel{}{\ell}{n}}\kern-10pt}
    \ar[dl]|{c_1}\ar[dr]|{c_2}&
    \\
    {\kern10pt\left(\!\begin{smallmatrix}\alpha' &\beta'
          \\s'&r'\end{smallmatrix}\!\right)_{\mylabel{}{(\ell+1)}{n_1}}\kern-20pt}&&{\left(\!\begin{smallmatrix}\alpha' &\beta'
          \\s'&r'\end{smallmatrix}\!\right)_{\mylabel{}{(\ell+1)}{n_2}}\kern-10pt}
  }
\end{equation*}
mean that the $\algU$ generators map via $c_1\xi_A^{1}+c_2\xi_A^{2}$,
where $\xi_A^{i}:\left(\!\begin{smallmatrix}\alpha &\beta
    \\s&r\end{smallmatrix}\!\right)_{\mylabel{}{\ell}{n}}\to
\left(\!\begin{smallmatrix}\alpha' &\beta'
    \\s'&r'\end{smallmatrix}\!\right)_{\mylabel{}{\ell+1}{n_i}}\!$ are
the maps in~\bref{sec:Ext1-atyp} associated with each link.  In other
words, the corresponding basis vectors in the two isomorphic
subquotients occur in linear combinations with $c_1$ and $c_2$ as the
coefficients.\footnote{The isomorphic targets could of course be
  identified differently, so as to ``split'' any given triple, but
  this cannot be done for all such triples simultaneously.  We did not
  attempt to ``optimize'' the possible linear combinations, in
  particular because there seems to be no well-defined optimum.}

The $\eta_{A}$ maps from the top-level subquotient are
\begin{align*}
  \eta_E\ZC{\alpha,s,\beta,0}_{n,\mylabel{}{1}{0}}&=-\qint{s\!-\!1}
  \qint{s}^4 \ZC{\alpha,s,\beta,0}_{n-1,\mylabel{Bt}{3}{0}}
 ,\\
  \eta_E\ZB{\alpha,s,\beta,0}_{n,\mylabel{}{1}{0}}&=-\qint{s\!-\!1}^2
  \qint{s}^3 \ZB{\alpha,s,\beta,0}_{n-1,\mylabel{Bt}{3}{0}}
 ,\\
  \eta_B\ZC{\alpha,s,\beta,0}_{n,\mylabel{}{1}{0}}&= \beta
  \qint{n}\qint{s}^3 \ZB{\alpha,s,\beta,0}_{n-1,\mylabel{Tb}{3}{5}}
  +\alpha \qint{s\!-\!1} \qint{s}^2 \qint{s\!-\!n}\ZB{\alpha,s,\beta,0}_{n-1,\mylabel{Bt}{3}{0}}
  \\
  &\quad{}+\alpha \beta \qint{s\!-\!1} \qint{s}^3 \delta_{0,n}\ZB{-\alpha,p\!-\!s,\beta,p\!-\!s}_{p-s,\mylabel{FB}{4}{3}}
  \\
  &\quad{}-\alpha\qint{n}\qint{s\!-\!1}^2 \qint{s}^2
  \ZC{\alpha,s\!+\!1,-\beta,0}_{n,\mylabel{BR}{4}{1}}
  \\
  &\quad{}+\alpha\beta\qint{s}^2(\qint{s\!-\!1} \qint{s\!-\!n}
  -\qint{n}\qint{2 s\!-\!1}) \ZB{\alpha,s,\beta,0}_{n-1,\mylabel{bot}{5}{0}}
 ,\\
  \eta_B\ZB{\alpha,s,\beta,0}_{n,\mylabel{}{1}{0}}&= \alpha
  \qint{n\!+\!1} \qint{s\!-\!1}^2 \qint{s}^2\ZB{\alpha,s\!+\!1,-\beta,0}_{n,\mylabel{BR}{4}{1}}
 ,\\
  \eta_C\ZC{\alpha,s,\beta,0}_{n,\mylabel{}{1}{0}}&=\qint{1\!-\!s}
  \qint{s}^4 \delta _{n,s-1}\ZC{-\alpha,p\!-\!s\!+\!1,\beta,p\!-\!s\!+\!1}_{0,\mylabel{El}{3}{1}}
  \\
  &\quad{}-\alpha \qint{s}^3 \ZC{\alpha,s\!-\!1,-\beta,0}_{n,\mylabel{BL}{4}{0}}
 ,\\
  \eta_C\ZB{\alpha,s,\beta,0}_{n,\mylabel{}{1}{0}}&=
  \qint{s\!-\!1}\qint{s}^4 \delta _{n,s-2}
  \ZB{-\alpha,p\!-\!s\!+\!1,\beta,p\!-\!s\!+\!1}_{0,\mylabel{El}{3}{1}}
  \\
  &\quad{}-\alpha \qint{s}^3 \ZB{\alpha,s\!-\!1,-\beta,0}_{n,\mylabel{BL}{4}{0}}+\alpha
  \qint{s}^2\ZC{\alpha,s,\beta,0}_{n+1,\mylabel{bot}{5}{0}}
\end{align*}
and those from level 2,
\begin{align*}
  \eta_B\ZC{-\alpha,p\!-\!s,-\beta,p\!-\!s}_{n,\mylabel{et}{2}{0}}&=
  -\alpha \qint{s\!-\!1} \delta_{0,n} \ZC{\alpha,s\!+\!1,-\beta,0}_{s,\mylabel{BR}{4}{1}}
  \\
  &\quad{}-\beta \qint{n+s} \ZB{-\alpha,p\!-\!s,-\beta,p\!-\!s}_{n,\mylabel{EB}{4}{2}}
 ,\\
  \eta_B\ZB{-\alpha,p\!-\!s,-\beta,p\!-\!s}_{n,\mylabel{et}{2}{0}}&= -\alpha \qint{s\!-\!1} \delta _{0,n} \ZB{\alpha,s\!+\!1,-\beta,0}_{s-1,\mylabel{BR}{4}{1}}
 ,\\
  \eta_B\ZC{-\alpha,p\!-\!s,\beta,p\!-\!s}_{n,\mylabel{ft}{2}{1}}&= \beta \qint{s\!-\!1} \qint{n+s} \ZB{-\alpha,p\!-\!s,\beta,p\!-\!s}_{n,\mylabel{FB}{4}{3}}
 ,\\
  \eta_E\ZC{\alpha,s\!-\!1,-\beta,0}_{n,\mylabel{tl}{2}{2}}&=
  \qint{s\!-\!1}^2\ZC{\alpha,s\!-\!1,-\beta,0}_{n-1,\mylabel{BL}{4}{0}}
 ,\\
  \eta_E\ZB{\alpha,s\!-\!1,-\beta,0}_{n,\mylabel{tl}{2}{2}}&=
  \qint{s-2}\qint{s\!-\!1} \ZB{\alpha,s\!-\!1,-\beta,0}_{n-1,\mylabel{BL}{4}{0}}
 ,\\
  \eta_B\ZC{\alpha,s\!-\!1,-\beta,0}_{n,\mylabel{tl}{2}{2}}&= \alpha \qint{\!n\!+\!1\!-\!s} \ZB{\alpha,s\!-\!1,-\beta,0}_{n-1,\mylabel{BL}{4}{0}}-\alpha \qint{n}\ZC{\alpha,s,\beta,0}_{n,\mylabel{bot}{5}{0}},\\
  \eta_B\ZB{\alpha,s\!-\!1,-\beta,0}_{n,\mylabel{tl}{2}{2}}&= \alpha \qint{n\!+\!1}\ZB{\alpha,s,\beta,0}_{n,\mylabel{bot}{5}{0}}
 ,\\
  \eta_E\ZC{\alpha,s\!+\!1,-\beta,0}_{n,\mylabel{tr}{2}{3}}&=
  \qint{s\!-\!1}\qint{s} \qint{s\!+\!1} \ZC{\alpha,s\!+\!1,-\beta,0}_{n-1,\mylabel{BR}{4}{1}}
 ,\\
  \eta_E\ZB{\alpha,s\!+\!1,-\beta,0}_{n,\mylabel{tr}{2}{3}}&=
  \qint{s\!-\!1}\qint{s}^2 \ZB{\alpha,s\!+\!1,-\beta,0}_{n-1,\mylabel{BR}{4}{1}}
 ,\\
  \eta_B\ZC{\alpha,s\!+\!1,-\beta,0}_{n,\mylabel{tr}{2}{3}}&= -\alpha \qint{s\!-\!1} \qint{s\!-\!n} \ZB{\alpha,s\!+\!1,-\beta,0}_{n-1,\mylabel{BR}{4}{1}}
 ,\\
  \eta_C\ZC{\alpha,s\!+\!1,-\beta,0}_{n,\mylabel{tr}{2}{3}}&=
  \qint{s\!-\!1}\qint{s} \delta _{n,s} \ZC{-\alpha,p\!-\!s,-\beta,p\!-\!s}_{0,\mylabel{EB}{4}{2}}
 ,\\
  \eta_C\ZB{\alpha,s\!+\!1,-\beta,0}_{n,\mylabel{tr}{2}{3}}&=
  -\qint{s\!-\!1}\qint{s} \delta _{n,s-1} \ZB{-\alpha,p\!-\!s,-\beta,p\!-\!s}_{0,\mylabel{EB}{4}{2}}
\end{align*}
and from level 3,
\begin{align*}
  \eta_B\ZC{\alpha,s,\beta,0}_{n,\mylabel{Bt}{3}{0}}&=-\beta \qint{n}
  \ZB{\alpha,s,\beta,0}_{n-1,\mylabel{bot}{5}{0}}
 ,\\
  \eta_B\ZC{-\alpha,p\!-\!s\!+\!1,\beta,p\!-\!s\!+\!1}_{n,\mylabel{El}{3}{1}}&=
  \alpha \delta_{0,n} \ZC{\alpha,s,\beta,0}_{s-1,\mylabel{bot}{5}{0}}
 ,\\
  \eta_B\ZB{-\alpha,p\!-\!s\!+\!1,\beta,p\!-\!s\!+\!1}_{n,\mylabel{El}{3}{1}}&=
  \alpha \delta_{0,n} \ZB{\alpha,s,\beta,0}_{s-2,\mylabel{bot}{5}{0}}
 ,\\
  \eta_E\ZC{\alpha,s,\beta,0}_{n,\mylabel{Tb}{3}{5}}&=\qint{s\!-\!1}
  \qint{s} \ZC{\alpha,s,\beta,0}_{n-1,\mylabel{bot}{5}{0}}
 ,\\
  \eta_E\ZB{\alpha,s,\beta,0}_{n,\mylabel{Tb}{3}{5}}&=\qint{s\!-\!1}^2
  \ZB{\alpha,s,\beta,0}_{n-1,\mylabel{bot}{5}{0}}
 ,\\
  \eta_B\ZC{\alpha,s,\beta,0}_{n,\mylabel{Tb}{3}{5}}&=-\alpha
  \qint{s\!-\!1\!-\!n} \ZB{\alpha,s,\beta,0}_{n-1,\mylabel{bot}{5}{0}}.
\end{align*}

\subsubsection{$\Zss$: Projective cover of
  $\boldsymbol{\repZ}^{\myboldsymbol{\alpha,\beta}}_{\myboldsymbol{s,s}}$
  for $\myboldsymbol{1\leq s\leq p-2}$}\label{proj:ss}
For each $1\leq s\leq p-2$, the projective module
$\repQ^{\alpha,\beta}_{s,s}$ (any of the $\zss$
in Fig.~\ref{module-pattern}) is $16p$-dimensional and has $16$
subquotients.  Its graph is shown in Fig.~\ref{fig:proj-ss}.
\afterpage{%
  \begin{landscape}
    \begin{figure}[htbp]
      \centering\footnotesize
      \begin{gather*}
        \xymatrix@R=25pt@C=20pt{ &&&&&\ssiAsBs="ssiAsBs"
          \\
          &&\ssiibcemApmsBZ="ssiibcemApmsBZ"&&
          \ssiibcfmApmsmBZ="ssiibcfmApmsmBZ"&&
          \ssiibcefAsmimBsmi="ssiibcefAsmimBsmi"&&
          \ssiibefAsimBsi="ssiibefAsimBsi"
          \\
          \\
          \\
          \ssiiiBcAsBs="ssiiiBcAsBs"&&
          \ssiiiCemApimsmBZ="ssiiiCemApimsmBZ"&&
          \ssiiiBemApmsmimBZ="ssiiiBemApmsmimBZ"&&
          \ssiiiCfmApimsBZ="ssiiiCfmApimsBZ"&&
          \ssiiiBfmApmsmiBZ="ssiiiBfmApmsmiBZ"&&
          \ssiiiCefAsBs="ssiiiCefAsBs"
          \\
          \\
          \\
          && \ssivCCAsmimBsmi="ssivCCAsmimBsmi"&&
          \ssivBBAsimBsi="ssivBBAsimBsi"&&
          \ssivEEmApmsBZ="ssivEEmApmsBZ"&&
          \ssivFFmApmsmBZ="ssivFFmApmsmBZ"
          \\
          &&&&&\ssvbotAsBs="ssvbotAsBs"
          \ar@{->}|(.45){\arlabel{\kern5pt\qint{s} \qint{s+1}}}
          "ssiAsBs"; "ssiibcemApmsBZ"
          \ar@{->}|(.45){\arlabel{\kern5pt\qint{s+1}}} "ssiAsBs";
          "ssiibcefAsmimBsmi"
          \ar@{->}|(.45){\arlabel{\kern5pt-\qint{s} \qint{s+1}}}
          "ssiAsBs"; "ssiibcfmApmsmBZ"
          \ar@{->}|(.45){\arlabel{\kern5pt\qint{s}^2}} "ssiAsBs";
          "ssiibefAsimBsi" \ar@{->}|(.45){\arlabel{\kern5pt-1}}
          "ssiibcemApmsBZ"; "ssiiiBcAsBs"
          \ar@{->}|(.45){\arlabel{\kern5pt-\alpha \beta}}
          "ssiibcefAsmimBsmi"; "ssiiiBcAsBs"
          \ar@{->}|(.45){\arlabel{\kern5pt1}} "ssiibcfmApmsmBZ";
          "ssiiiBcAsBs" \ar@{->}|(.45){\arlabel{\kern5pt\qint{s}}}
          "ssiibcemApmsBZ"; "ssiiiBemApmsmimBZ"
          \ar@{->}|(.35){\arlabel{\kern5pt\qint{s+1}}}
          "ssiibefAsimBsi"; "ssiiiBemApmsmimBZ"
          \ar@{->}|(.45){\arlabel{\kern5pt-\qint{s}}}
          "ssiibcfmApmsmBZ"; "ssiiiBfmApmsmiBZ"
          \ar@{->}|(.70){\arlabel{\kern5pt\qint{s+1}}}
          "ssiibefAsimBsi"; "ssiiiBfmApmsmiBZ"
          \ar@{->}|(.35){\arlabel{\kern5pt\qint{s+1}}}
          "ssiibcemApmsBZ"; "ssiiiCemApimsmBZ"
          \ar@{->}|(.55){\arlabel{\kern5pt\qint{s} \qint{s+1}}}
          "ssiibcefAsmimBsmi"; "ssiiiCemApimsmBZ"
          \ar@{->}|(.45){\arlabel{\kern5pt\qint{s}}}
          "ssiibcefAsmimBsmi"; "ssiiiCefAsBs"
          \ar@{->}|(.45){\arlabel{\kern5pt1}} "ssiibefAsimBsi";
          "ssiiiCefAsBs" \ar@{->}|(.70){\arlabel{\kern5pt\qint{s}
              \qint{s+1}}} "ssiibcefAsmimBsmi"; "ssiiiCfmApimsBZ"
          \ar@{->}|(.45){\arlabel{\kern5pt\qint{s+1}}}
          "ssiibcfmApmsmBZ"; "ssiiiCfmApimsBZ"
          \ar@{->}|(.62){\arlabel{\kern5pt-\qint{s}}} "ssiiiBcAsBs";
          "ssivBBAsimBsi" \ar@{->}|(.73){\arlabel{\kern5pt1}}
          "ssiiiBemApmsmimBZ"; "ssivBBAsimBsi"
          \ar@{->}|(.45){\arlabel{\kern5pt1}} "ssiiiBfmApmsmiBZ";
          "ssivBBAsimBsi" \ar@{->}|(.55){\arlabel{\kern5pt-\alpha
              \beta}} "ssiiiCefAsBs"; "ssivBBAsimBsi"
          \ar@{->}|(.45){\arlabel{\kern5pt-\qint{s+1}}} "ssiiiBcAsBs";
          "ssivCCAsmimBsmi" \ar@{->}|(.65){\arlabel{\kern5pt1}}
          "ssiiiCemApimsmBZ"; "ssivCCAsmimBsmi"
          \ar@{->}|(.45){\arlabel{\kern5pt1}} "ssiiiCfmApimsBZ";
          "ssivCCAsmimBsmi" \ar@{->}|(.50){\arlabel{\kern5pt1}}
          "ssiiiBemApmsmimBZ"; "ssivEEmApmsBZ"
          \ar@{->}|(.60){\arlabel{\kern5pt1}} "ssiiiCemApimsmBZ";
          "ssivEEmApmsBZ" \ar@{->}|(.77){\arlabel{\kern5pt\qint{s+1}}}
          "ssiiiCefAsBs"; "ssivEEmApmsBZ"
          \ar@{->}|(.65){\arlabel{\kern5pt-1}} "ssiiiBfmApmsmiBZ";
          "ssivFFmApmsmBZ"
          \ar@{->}|(.45){\arlabel{\kern5pt\qint{s+1}}} "ssiiiCefAsBs";
          "ssivFFmApmsmBZ" \ar@{->}|(.55){\arlabel{\kern5pt1}}
          "ssiiiCfmApimsBZ"; "ssivFFmApmsmBZ"
          \ar@{->}|(.45){\arlabel{\kern5pt1}} "ssivBBAsimBsi";
          "ssvbotAsBs" \ar@{->}|(.45){\arlabel{\kern5pt1}}
          "ssivCCAsmimBsmi"; "ssvbotAsBs"
          \ar@{->}|(.45){\arlabel{\kern5pt1}} "ssivEEmApmsBZ";
          "ssvbotAsBs" \ar@{->}|(.45){\arlabel{\kern5pt1}}
          "ssivFFmApmsmBZ"; "ssvbotAsBs"}
      \end{gather*}
      \caption[Graph of the projective module
      $\repQ^{\alpha,\beta}_{s,s}$]{\small\ \ Graph of the projective
        module $\repQ^{\alpha,\beta}_{s,s}$ for $1\leq s\leq p-2$.}
      \label{fig:proj-ss}
    \end{figure}
  \end{landscape}
}%
with all the previous conventions in force.

The $\eta_{A}$ piece of the action of $\algU$ generators on
$\repQ^{\alpha,\beta}_{s,s}$ is as follows.  On the basis vectors of
the top-level subquotient, we have
\begin{align*}
  \eta_E\ZC{\alpha,s,\beta,s}_{n,\mylabel{}{1}{0}}&=
  -\qint{s}^2\qint{s\!+\!1} \ZC{\alpha,s,\beta,s}_{n-1,\mylabel{Bc}{3}{0}} + 2 \beta \qint{s}^2
  \qint{s\!+\!1} \ZC{\alpha,s,\beta,s}_{n-1,\mylabel{bot}{5}{0}}
 ,\\
  \eta_E\ZB{\alpha,s,\beta,s}_{n,\mylabel{}{1}{0}}&=
  -\qint{s}\qint{s\!+\!1}^2 \ZB{\alpha,s,\beta,s}_{n-1,\mylabel{Bc}{3}{0}} + 3 \beta \qint{s}
  \qint{s\!+\!1}^2 \ZB{\alpha,s,\beta,s}_{n-1,\mylabel{bot}{5}{0}}
 ,\\
  \eta_B\ZC{\alpha,s,\beta,s}_{n,\mylabel{}{1}{0}}&= -\alpha
  \qint{n} \qint{s\!+\!1} \ZB{\alpha,s,\beta,s}_{n,\mylabel{Bc}{3}{0}} -\alpha
  \qint{s}\qint{s\!+\!1}^2 \delta _{0,n}
  \ZC{-\alpha,p\!-\!s\!+\!1,\beta,0}_{p-s,\mylabel{Cf}{3}{3}}
  \\
  &\quad{}-\beta \qint{s} \qint{s\!-\!n} \ZB{\alpha,s,\beta,s}_{n,\mylabel{Cef}{3}{5}}-\alpha \qint{n\!+\!1}\qint{s}
  \qint{s\!+\!1}^2 \ZC{\alpha,s\!-\!1,-\beta,s\!-\!1}_{n-1,\mylabel{CC}{4}{0}}
  \\
  &\quad{}+\alpha \beta \qint{s\!+\!1}
  (\qint{n}-\qint{s}\qint{s\!-\!1\!-\!n}) \ZB{\alpha,s,\beta,s}_{n,\mylabel{bot}{5}{0}}
 ,\\
  \eta_B\ZB{\alpha,s,\beta,s}_{n,\mylabel{}{1}{0}}&= \alpha
  \qint{n\!+\!1} \qint{s}\qint{s\!+\!1}^2\ZB{\alpha,s\!-\!1,-\beta,s\!-\!1}_{n-1,\mylabel{CC}{4}{0}}
  \\
  &\quad{}-\alpha \qint{s}\qint{s\!+\!1}^2 \delta _{0,n}
  \ZB{-\alpha,p\!-\!s\!+\!1,\beta,0}_{p-s-1,\mylabel{Cf}{3}{3}}.
\end{align*}
On basis vectors of the level-$2$ modules, the $\eta$ maps are
\begin{align*}
  \eta_E\ZC{-\alpha,p\!-\!s,\beta,0}_{n,\mylabel{bce}{2}{0}}&= 2
  \beta \qint{s}\delta _{0,n} \ZC{\alpha,s,\beta,s}_{s-1,\mylabel{bot}{5}{0}}
 ,\\
  \eta_E\ZB{-\alpha,p\!-\!s,\beta,0}_{n,\mylabel{bce}{2}{0}}&= 3
  \beta \qint{s\!+\!1} \delta _{0,n} \ZB{\alpha,s,\beta,s}_{s,\mylabel{bot}{5}{0}}
 ,\\
  \eta_B\ZC{-\alpha,p\!-\!s,\beta,0}_{n,\mylabel{bce}{2}{0}}&= \beta
  \qint{n}\ZB{-\alpha,p\!-\!s,\beta,0}_{n-1,\mylabel{EE}{4}{2}} 
  + 2 \alpha \beta \delta_{0,n} \ZB{\alpha,s,\beta,s}_{s,\mylabel{bot}{5}{0}}
 ,\\
  \eta_E\ZC{\alpha,s\!-\!1,-\beta,s\!-\!1}_{n,\mylabel{bcef}{2}{2}}&=
  \qint{s\!-\!1}\qint{s} \qint{s\!+\!1} \ZC{\alpha,s\!-\!1,-\beta,s\!-\!1}_{n-1,\mylabel{CC}{4}{0}}
 ,\\
  \eta_E\ZB{\alpha,s\!-\!1,-\beta,s\!-\!1}_{n,\mylabel{bcef}{2}{2}}&=
  \qint{s}^2\qint{s\!+\!1} \ZB{\alpha,s\!-\!1,-\beta,s\!-\!1}_{n-1,\mylabel{CC}{4}{0}}
 ,\\
  \eta_B\ZC{\alpha,s\!-\!1,-\beta,s\!-\!1}_{n,\mylabel{bcef}{2}{2}}&=
  \alpha \qint{n\!+\!1} \qint{s\!+\!1} \ZB{\alpha,s\!-\!1,-\beta,s\!-\!1}_{n,\mylabel{CC}{4}{0}}
 ,\\
  \eta_C\ZB{\alpha,s\!-\!1,-\beta,s\!-\!1}_{n,\mylabel{bcef}{2}{2}}&=
  \alpha \ZB{\alpha,s,\beta,s}_{n+1,\mylabel{bot}{5}{0}}
 ,\\
  \eta_F\ZC{-\alpha,p\!-\!s,-\beta,0}_{n,\mylabel{bcf}{2}{1}}&= -2
  \beta \delta_{n,p-s-1} \ZC{\alpha,s,\beta,s}_{0,\mylabel{bot}{5}{0}}
 ,\\
  \eta_F\ZB{-\alpha,p\!-\!s,-\beta,0}_{n,\mylabel{bcf}{2}{1}}&= 3
  \beta \delta_{n,p-s-2} \ZB{\alpha,s,\beta,s}_{0,\mylabel{bot}{5}{0}}
 ,\\
  \eta_B\ZC{-\alpha,p\!-\!s,-\beta,0}_{n,\mylabel{bcf}{2}{1}}&=
  -\beta \qint{n}\ZB{-\alpha,p\!-\!s,-\beta,0}_{n-1,\mylabel{FF}{4}{3}}
 ,\\
  \eta_C\ZC{-\alpha,p\!-\!s,-\beta,0}_{n,\mylabel{bcf}{2}{1}}&=
  -\qint{s} \delta_{n,p-s-1} \ZC{\alpha,s\!+\!1,-\beta,s\!+\!1}_{0,\mylabel{BB}{4}{1}}
 ,\\
  \eta_C\ZB{-\alpha,p\!-\!s,-\beta,0}_{n,\mylabel{bcf}{2}{1}}&=
  \qint{s} \delta_{n,p-s-2} \ZB{\alpha,s\!+\!1,-\beta,s\!+\!1}_{0,\mylabel{BB}{4}{1}}
 ,\\
  \eta_E\ZC{\alpha,s\!+\!1,-\beta,s\!+\!1}_{n,\mylabel{bef}{2}{3}}&=
  \qint{s\!+\!1}^2\ZC{\alpha,s\!+\!1,-\beta,s\!+\!1}_{n-1,\mylabel{BB}{4}{1}}
 ,\\
  \eta_E\ZB{\alpha,s\!+\!1,-\beta,s\!+\!1}_{n,\mylabel{bef}{2}{3}}&=
  \qint{s\!+\!1}\qint{s+2} \ZB{\alpha,s\!+\!1,-\beta,s\!+\!1}_{n-1,\mylabel{BB}{4}{1}}
 ,\\
  \eta_B\ZC{\alpha,s\!+\!1,-\beta,s\!+\!1}_{n,\mylabel{bef}{2}{3}}&=
  \alpha \qint{n} \ZB{\alpha,s\!+\!1,-\beta,s\!+\!1}_{n,\mylabel{BB}{4}{1}}-\alpha
  \qint{n\!+\!1}\qint{s\!+\!1} \ZC{\alpha,s,\beta,s}_{n-1,\mylabel{bot}{5}{0}}
 ,\\
  &\quad{}-\alpha \qint{s\!+\!1}\delta _{0,n} \ZC{-\alpha,p\!-\!s,-\beta,0}_{p-s-1,\mylabel{FF}{4}{3}}
 ,\\
  \eta_B\ZB{\alpha,s\!+\!1,-\beta,s\!+\!1}_{n,\mylabel{bef}{2}{3}}&=
  \alpha \qint{n\!+\!1} \qint{s\!+\!1} \ZB{\alpha,s,\beta,s}_{n-1,\mylabel{bot}{5}{0}}
  \\
  &\quad{}-\alpha \qint{s\!+\!1} \delta _{0,n}
  \ZB{-\alpha,p\!-\!s,-\beta,0}_{p-s-2,\mylabel{FF}{4}{3}}
\end{align*}
and on level-$3$ vectors,
\begin{align*}
  \eta_C\ZB{\alpha,s,\beta,s}_{n,\mylabel{Bc}{3}{0}}&=\ZC{\alpha,s,\beta,s}_{n,\mylabel{bot}{5}{0}}
 ,\\
  \eta_E\ZC{\alpha,s,\beta,s}_{n,\mylabel{Cef}{3}{5}}&=
  \qint{s}\qint{s\!+\!1} \ZC{\alpha,s,\beta,s}_{n-1,\mylabel{bot}{5}{0}}
 ,\\
  \eta_E\ZB{\alpha,s,\beta,s}_{n,\mylabel{Cef}{3}{5}}&=
  \qint{s\!+\!1}^2\ZB{\alpha,s,\beta,s}_{n-1,\mylabel{bot}{5}{0}}
 ,\\
  \eta_B\ZC{\alpha,s,\beta,s}_{n,\mylabel{Cef}{3}{5}}&= \alpha
  \qint{n\!+\!1} \ZB{\alpha,s,\beta,s}_{n,\mylabel{bot}{5}{0}}
 ,\\
  \eta_C\ZC{-\alpha,p\!-\!s\!+\!1,\beta,0}_{n,\mylabel{Cf}{3}{3}}&=
  \delta _{n,p-s}\ZC{\alpha,s,\beta,s}_{0,\mylabel{bot}{5}{0}}
 ,\\
  \eta_C\ZB{-\alpha,p\!-\!s\!+\!1,\beta,0}_{n,\mylabel{Cf}{3}{3}}&=
  -\delta _{n,p-s-1}\ZB{\alpha,s,\beta,s}_{0,\mylabel{bot}{5}{0}}.
\end{align*}

\subsubsection{$\Zoz$: Projective cover of
  $\boldsymbol{\repZ}^{\myboldsymbol{\alpha,\beta}}_{\myboldsymbol{1,0}}$}
\label{proj:10}
We finally describe projective covers of the one-dimensional
representations (\,$\zoz$ in Fig.~\ref{module-pattern}).  The
projective module $\repQ^{\alpha,\beta}_{1,0}$ has dimension $24p$ and
is built from $24$ simple subquotients.  Its rather involved graph is
shown in Fig.~\ref{fig:proj-10}.
\afterpage{%
  \begin{landscape}
    \begin{figure}[htbp]
      \centering\footnotesize
      \begin{gather*}
        \mbox{}\xymatrix@R=30pt@C=5pt{ &&&&&&&&&\oziAIBZ="oziAIBZ"
          \\
          &&&& \oziibIIIAIBI="oziibIIIAIBI"&&
          \oziieIImApmBZ="oziieIImApmBZ"&&
          \oziifImApBZ="oziifImApBZ"&&
          \oziieVImApmimBpmi="oziieVImApmimBpmi"&&
          \oziifVmApmiBpmi="oziifVmApmiBpmi"&&
          \oziicIVAIImBZ="oziicIVAIImBZ"
          \\
          \\
          \\
          \oziiibIImApmiBZ="oziiibIImApmiBZ"&&
          \oziiibImApmimBZ="oziiibImApmimBZ"&&
          \oziiisVAIBZ="oziiisVAIBZ"&&
          \oziiiltVIIIAImBZ="oziiiltVIIIAImBZ"&&
          \oziiimtIVAIBZ="oziiimtIVAIBZ"&&
          \oziiisVIAIBZ="oziiisVIAIBZ"&&
          \oziiihtIIIAImBZ="oziiihtIIIAImBZ"&&
          \oziiisVIIAIBZ="oziiisVIIAIBZ"&&
          \oziiicIZmApmiiBpmii="oziiicIZmApmiiBpmii"&&
          \oziiicIXmApmiimBpmii="oziiicIXmApmiimBpmii"
          \\
          \\
          \\
          &&&& \ozivEIImApmBZ="ozivEIImApmBZ"&&
          \ozivBIIIAIBI="ozivBIIIAIBI"&& \ozivFImApBZ="ozivFImApBZ"&&
          \ozivEVImApmimBpmi="ozivEVImApmimBpmi"&&
          \ozivFVmApmiBpmi="ozivFVmApmiBpmi"&&
          \ozivCIVAIImBZ="ozivCIVAIImBZ"
          \\
          &&&&&&&&&\ozvbotAIBZ="ozvbotAIBZ"
          \ar@{->}|(.45){\arlabel{\kern3pt2 \alpha \beta\kern2pt}}
          "oziAIBZ"; "oziibIIIAIBI" \ar@{->}|(.45){\arlabel{\kern3pt-2
              \alpha \beta\kern2pt}} "oziAIBZ"+<22pt,-6pt>;
          "oziicIVAIImBZ" \ar@{->}|(.45){\arlabel{\kern3pt1\kern2pt}}
          "oziAIBZ"; "oziieIImApmBZ"
          \ar@{->}|(.45){\arlabel{\kern5pt1\kern5pt}} "oziAIBZ";
          "oziieVImApmimBpmi"
          \ar@{->}|(.45){\arlabel{\kern3pt1\kern2pt}} "oziAIBZ";
          "oziifImApBZ" \ar@{->}|(.45){\arlabel{\kern3pt-1\kern2pt}}
          "oziAIBZ"+<22pt,-8pt>; "oziifVmApmiBpmi"
          \ar@{->}|(.75){\arlabel{\kern3pt1\kern2pt}} "oziibIIIAIBI";
          "oziiibImApmimBZ" \ar@{->}|(.77){\arlabel{\kern3pt-2 \alpha
              \beta}} "oziifImApBZ"; "oziiibImApmimBZ"
          \ar@{->}|(.45){\arlabel{\kern1pt1}} "oziibIIIAIBI";
          "oziiibIImApmiBZ" \ar@{->}|(.40){\arlabel{\kern3pt2 \alpha
              \beta}} "oziieIImApmBZ"; "oziiibIImApmiBZ"
          \ar@{->}|(.75){\arlabel{\kern2pt1\kern2pt}} "oziicIVAIImBZ";
          "oziiicIZmApmiiBpmii" \ar@{->}|(.75){\arlabel{\kern3pt2
              \alpha \beta}} "oziieVImApmimBpmi";
          "oziiicIZmApmiiBpmii" \ar@{->}|(.45){\arlabel{\kern3pt1}}
          "oziicIVAIImBZ"; "oziiicIXmApmiimBpmii"
          \ar@{->}|(.65){\arlabel{\kern3pt2 \alpha \beta}}
          "oziifVmApmiBpmi"; "oziiicIXmApmiimBpmii"
          \ar@{->}|(.41){\arlabel{\kern3pt2}} "oziifImApBZ";
          "oziiihtIIIAImBZ" \ar@{->}|(.75){\arlabel{\kern3pt2}}
          "oziifVmApmiBpmi"; "oziiihtIIIAImBZ"
          \ar@{->}|(.58){\arlabel{\kern3pt-1}} "oziieIImApmBZ";
          "oziiiltVIIIAImBZ" \ar@{->}|(.80){\arlabel{\kern3pt1}}
          "oziieVImApmimBpmi"; "oziiiltVIIIAImBZ"
          \ar@{->}|(.45){\arlabel{\kern3pt1}} "oziieIImApmBZ";
          "oziiimtIVAIBZ" \ar@{->}|(.80){\arlabel{\kern3pt-1}}
          "oziieVImApmimBpmi"; "oziiimtIVAIBZ"
          \ar@{->}|(.25){\arlabel{\kern3pt-1}} "oziifImApBZ";
          "oziiimtIVAIBZ" \ar@{->}|(.35){\arlabel{\kern3pt-1}}
          "oziifVmApmiBpmi"; "oziiimtIVAIBZ"
          \ar@{->}|(.48){\arlabel{\kern3pt1}} "oziibIIIAIBI";
          "oziiisVAIBZ" \ar@{->}|(.30){\arlabel{\kern3pt1}}
          "oziieIImApmBZ"; "oziiisVAIBZ"
          \ar@{->}|(.48){\arlabel{\kern3pt-1}} "oziieVImApmimBpmi";
          "oziiisVAIBZ" \ar@{->}|(.71){\arlabel{\kern3pt1}}
          "oziifImApBZ"; "oziiisVAIBZ"
          \ar@{->}|(.70){\arlabel{\kern3pt1}} "oziifVmApmiBpmi";
          "oziiisVAIBZ" \ar@{->}|(.55){\arlabel{\kern3pt-1}}
          "oziicIVAIImBZ"; "oziiisVIAIBZ"
          \ar@{->}|(.81){\arlabel{\kern3pt1}} "oziieIImApmBZ";
          "oziiisVIAIBZ" \ar@{->}|(.80){\arlabel{\kern3pt-1}}
          "oziieVImApmimBpmi"; "oziiisVIAIBZ"
          \ar@{->}|(.57){\arlabel{\kern3pt1}} "oziifImApBZ";
          "oziiisVIAIBZ" \ar@{->}|(.47){\arlabel{\kern3pt1}}
          "oziifVmApmiBpmi"; "oziiisVIAIBZ"
          \ar@{->}|(.47){\arlabel{\kern3pt1}} "oziicIVAIImBZ";
          "oziiisVIIAIBZ" \ar@{->}|(.70){\arlabel{\kern3pt2}}
          "oziieVImApmimBpmi"; "oziiisVIIAIBZ"
          \ar@{->}|(.70){\arlabel{\kern3pt-2}} "oziifVmApmiBpmi";
          "oziiisVIIAIBZ" \ar@{->}|(.35){\arlabel{\kern3pt1}}
          "oziiibImApmimBZ"; "ozivBIIIAIBI"
          \ar@{->}|(.45){\arlabel{\kern3pt1}} "oziiibIImApmiBZ";
          "ozivBIIIAIBI" \ar@{->}|(.25){\arlabel{\kern3pt\alpha
              \beta}} "oziiisVAIBZ"; "ozivBIIIAIBI"
          \ar@{->}|(.20){\arlabel{\kern3pt\alpha \beta}}
          "oziiisVIAIBZ"; "ozivBIIIAIBI"
          \ar@{->}|(.30){\arlabel{\kern3pt\alpha \beta}}
          "oziiisVIIAIBZ"; "ozivBIIIAIBI"
          \ar@{->}|(.30){\arlabel{\kern3pt1}} "oziiicIZmApmiiBpmii";
          "ozivCIVAIImBZ" \ar@{->}|(.45){\arlabel{\kern3pt1}}
          "oziiicIXmApmiimBpmii"; "ozivCIVAIImBZ"
          \ar@{->}|(.47){\arlabel{\kern3pt-\alpha \beta}}
          "oziiisVIIAIBZ"; "ozivCIVAIImBZ"
          \ar@{->}|(.45){\arlabel{\kern3pt1}} "oziiibIImApmiBZ";
          "ozivEIImApmBZ" \ar@{->}|(.21){\arlabel{\kern3pt-2}}
          "oziiiltVIIIAImBZ"; "ozivEIImApmBZ"
          \ar@{->}|(.35){\arlabel{\kern3pt1}} "oziiimtIVAIBZ";
          "ozivEIImApmBZ" \ar@{->}|(.80){\arlabel{\kern-5pt1}}
          "oziiisVAIBZ"; "ozivEIImApmBZ"
          \ar@{->}|(.45){\arlabel{\kern3pt1}} "oziiicIZmApmiiBpmii";
          "ozivEVImApmimBpmi" \ar@{->}|(.55){\arlabel{\kern3pt2}}
          "oziiiltVIIIAImBZ"; "ozivEVImApmimBpmi"
          \ar@{->}|(.80){\arlabel{\kern3pt-1}} "oziiimtIVAIBZ";
          "ozivEVImApmimBpmi" \ar@{->}|(.25){\arlabel{\kern3pt-1}}
          "oziiisVIAIBZ"; "ozivEVImApmimBpmi"
          \ar@{->}|(.70){\arlabel{\kern3pt1}} "oziiibImApmimBZ";
          "ozivFImApBZ" \ar@{->}|(.80){\arlabel{\kern3pt1}}
          "oziiihtIIIAImBZ"; "ozivFImApBZ"
          \ar@{->}|(.20){\arlabel{\kern3pt-1}} "oziiimtIVAIBZ";
          "ozivFImApBZ" \ar@{->}|(.40){\arlabel{\kern3pt1}}
          "oziiisVAIBZ"; "ozivFImApBZ"
          \ar@{->}|(.82){\arlabel{\kern3pt1}} "oziiicIXmApmiimBpmii";
          "ozivFVmApmiBpmi" \ar@{->}|(.45){\arlabel{\kern3pt-1}}
          "oziiihtIIIAImBZ"; "ozivFVmApmiBpmi"
          \ar@{->}|(.68){\arlabel{\kern3pt1}} "oziiimtIVAIBZ";
          "ozivFVmApmiBpmi" \ar@{->}|(.50){\arlabel{\kern3pt-1}}
          "oziiisVIAIBZ"; "ozivFVmApmiBpmi"
          \ar@{->}|(.45){\arlabel{\kern3pt1\kern2pt}} "ozivBIIIAIBI";
          "ozvbotAIBZ" \ar@{->}|(.45){\arlabel{\kern3pt1\kern2pt}}
          "ozivCIVAIImBZ"; "ozvbotAIBZ"
          \ar@{->}|(.45){\arlabel{\kern3pt1\kern2pt}} "ozivEIImApmBZ";
          "ozvbotAIBZ" \ar@{->}|(.45){\arlabel{\kern3pt1\kern2pt}}
          "ozivEVImApmimBpmi"; "ozvbotAIBZ"
          \ar@{->}|(.45){\arlabel{\kern3pt1\kern2pt}} "ozivFImApBZ";
          "ozvbotAIBZ" \ar@{->}|(.45){\arlabel{\kern3pt1\kern2pt}}
          "ozivFVmApmiBpmi"; "ozvbotAIBZ"}
      \end{gather*}%
      \caption[Graph of the projective module
      $\repQ^{\alpha,\beta}_{1,0}$]{\small Graph of the projective
        module $\repQ^{\alpha,\beta}_{1,0}$}
      \label{fig:proj-10}
    \end{figure}
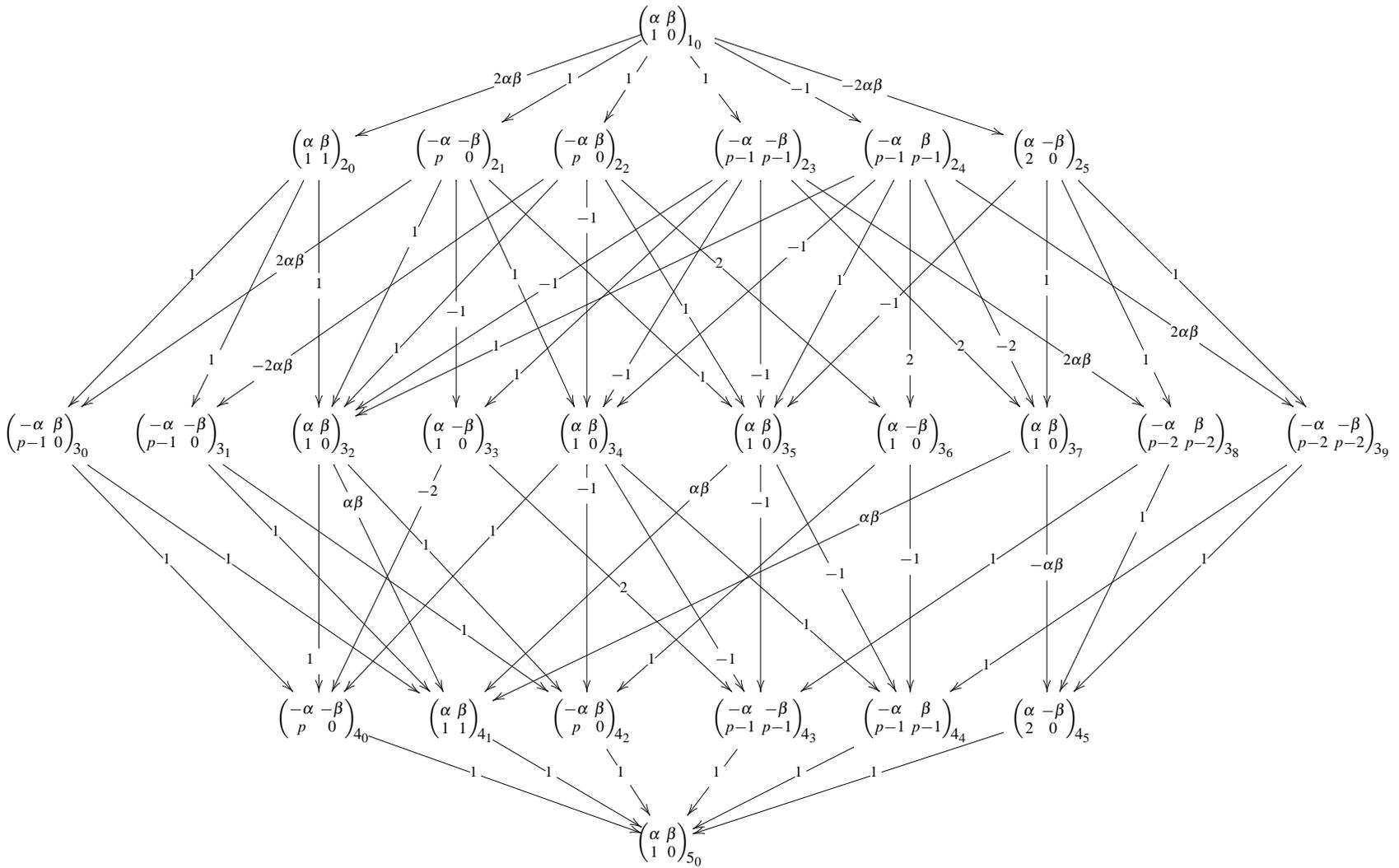
  \end{landscape}%
}%
Maps into linear combinations of isomorphic subquotients, which
already occurred in Figs.~\ref{fig:proj-so} and~\ref{fig:proj-ss},
here involve up to four modules.

The $\eta$ piece of the action of $\algU$ generators is as follows.
The maps from the top are
\begin{align*}
  \eta_F\ZC{\alpha,1,\beta,0}_{0,\mylabel{}{1}{0}}&= 2 \alpha \qint{2}
  \ZC{-\alpha,p\!-\!1,-\beta,p\!-\!1}_{0,\mylabel{E6}{4}{3}}
  \\
  \eta_E\ZC{\alpha,1,\beta,0}_{0,\mylabel{}{1}{0}}&= 2 \alpha \qint{2}
  \ZC{-\alpha,p\!-\!1,\beta,p\!-\!1}_{p-2,\mylabel{F5}{4}{4}},
\end{align*}
and those from level $2$ are
\begin{align*}
  \eta_E\ZB{\alpha,1,\beta,1}_{n,\mylabel{b3}{2}{0}}&=\qint{2} \ZB{\alpha,1,\beta,1}_{n-1,\mylabel{B3}{4}{1}}
 ,\\
  \eta_B\ZC{\alpha,1,\beta,1}_{n,\mylabel{b3}{2}{0}}&=\alpha \delta_{0,n} \ZC{-\alpha,p,\beta,0}_{p-1,\mylabel{F1}{4}{2}}
 ,\\
  \eta_B\ZB{\alpha,1,\beta,1}_{n,\mylabel{b3}{2}{0}}&=\alpha \delta_{0,n} \ZB{-\alpha,p,\beta,0}_{p-2,\mylabel{F1}{4}{2}}
 ,\\
  \eta_E\ZC{\alpha,2,-\beta,0}_{n,\mylabel{c4}{2}{5}}&= \qint{2}\ZC{\alpha,2,-\beta,0}_{n-1,\mylabel{C4}{4}{5}}
 ,\\
  \eta_C\ZC{\alpha,2,-\beta,0}_{n,\mylabel{c4}{2}{5}}&= \delta_{1,n}\ZC{-\alpha,p\!-\!1,-\beta,p\!-\!1}_{0,\mylabel{E6}{4}{3}}
 ,\\
  \eta_C\ZB{\alpha,2,-\beta,0}_{n,\mylabel{c4}{2}{5}}&= -\delta_{0,n}\ZB{-\alpha,p\!-\!1,-\beta,p\!-\!1}_{0,\mylabel{E6}{4}{3}}
 ,\\
  \eta_E\ZB{-\alpha,p,-\beta,0}_{n,\mylabel{e2}{2}{1}}&= 2
  \ZB{-\alpha,p,-\beta,0}_{n-1,\mylabel{E2}{4}{0}}-2\ZC{-\alpha,p\!-\!1,-\beta,p\!-\!1}_{n-1,\mylabel{E6}{4}{3}}
 ,\\
  \eta_B\ZC{-\alpha,p,-\beta,0}_{n,\mylabel{e2}{2}{1}}&= 2 \alpha \qint{n\!+\!1} \ZB{-\alpha,p,-\beta,0}_{n-1,\mylabel{E2}{4}{0}}\\
  &\quad{}-2 \alpha \qint{n\!+\!1} \ZC{-\alpha,p\!-\!1,-\beta,p\!-\!1}_{n-1,\mylabel{E6}{4}{3}}
 ,\\
  \eta_B\ZB{-\alpha,p,-\beta,0}_{n,\mylabel{e2}{2}{1}}&= -2 \alpha
  \qint{n+2} \ZB{-\alpha,p\!-\!1,-\beta,p\!-\!1}_{n,\mylabel{E6}{4}{3}}
 ,\\
  \eta_E\ZC{-\alpha,p\!-\!1,-\beta,p\!-\!1}_{n,\mylabel{e6}{2}{3}}&= 2
  \ZC{-\alpha,p\!-\!1,-\beta,p\!-\!1}_{n-1,\mylabel{E6}{4}{3}}-2\ZB{-\alpha,p,-\beta,0}_{n-1,\mylabel{E2}{4}{0}}
 ,\\
  \eta_B\ZC{-\alpha,p\!-\!1,-\beta,p\!-\!1}_{n,\mylabel{e6}{2}{3}}&= 2
  \beta \delta _{0,n}\ZC{\alpha,2,-\beta,0}_{1,\mylabel{C4}{4}{5}} -2
  \alpha \qint{n}\ZB{-\alpha,p\!-\!1,-\beta,p\!-\!1}_{n,\mylabel{E6}{4}{3}}
 ,\\
  \eta_B\ZB{-\alpha,p\!-\!1,-\beta,p\!-\!1}_{n,\mylabel{e6}{2}{3}}&= 2
  \beta \delta _{0,n}\ZB{\alpha,2,-\beta,0}_{0,\mylabel{C4}{4}{5}}
 ,\\
  \eta_E\ZB{-\alpha,p,\beta,0}_{n,\mylabel{f1}{2}{2}}&= 2
  \ZB{-\alpha,p,\beta,0}_{n-1,\mylabel{F1}{4}{2}}-2\ZC{-\alpha,p\!-\!1,\beta,p\!-\!1}_{n-1,\mylabel{F5}{4}{4}}
 ,\\
  \eta_B\ZC{-\alpha,p,\beta,0}_{n,\mylabel{f1}{2}{2}}&= 2 \alpha \qint{n\!+\!1} \ZB{-\alpha,p,\beta,0}_{n-1,\mylabel{F1}{4}{2}}\\
  &\quad{} -2 \alpha \qint{n\!+\!1} \ZC{-\alpha,p\!-\!1,\beta,p\!-\!1}_{n-1,\mylabel{F5}{4}{4}}
 ,\\
  \eta_B\ZB{-\alpha,p,\beta,0}_{n,\mylabel{f1}{2}{2}}&= -2 \alpha
  \qint{n+2} \ZB{-\alpha,p\!-\!1,\beta,p\!-\!1}_{n,\mylabel{F5}{4}{4}}
 ,\\
  \eta_C\ZC{-\alpha,p,\beta,0}_{n,\mylabel{f1}{2}{2}}&= -2 \alpha \beta
  \delta_{n,p\!-\!1} \ZC{\alpha,1,\beta,1}_{0,\mylabel{B3}{4}{1}}
 ,\\
  \eta_C\ZB{-\alpha,p,\beta,0}_{n,\mylabel{f1}{2}{2}}&= 2 \alpha \beta
  \delta_{n,p-2} \ZB{\alpha,1,\beta,1}_{0,\mylabel{B3}{4}{1}}
 ,\\
  \eta_E\ZC{-\alpha,p\!-\!1,\beta,p\!-\!1}_{n,\mylabel{f5}{2}{4}}&= 2
  \ZB{-\alpha,p,\beta,0}_{n-1,\mylabel{F1}{4}{2}}-2\ZC{-\alpha,p\!-\!1,\beta,p\!-\!1}_{n-1,\mylabel{F5}{4}{4}}
 ,\\
  \eta_B\ZC{-\alpha,p\!-\!1,\beta,p\!-\!1}_{n,\mylabel{f5}{2}{4}}&= 2
  \alpha \qint{n} \ZB{-\alpha,p\!-\!1,\beta,p\!-\!1}_{n,\mylabel{F5}{4}{4}}.
\end{align*}
Together with the $c^{\kappa,\iota}$ specified on the links, this
defines a $\algU$ module, which is then seen to be maximal and
indecomposable.

\subsection{Completeness} 
We have constructed projective covers $\repQ_i\to\repZ_i$ for all
simple $\algU$ modules.  That these are all projective modules can
also be verified by calculating the sum $\sum_i\dim\repQ_i \cdot
\dim\repZ_i$:
\begin{multline}
  4\biggl( \underbrace{(p-1)4p\cdot 4p}_{\bref{proj:steinberg}} {}+{}
  \underbrace{(p-2)\sum_{s=1}^{p-1} 8p\cdot 4s}_{\bref{proj:typical}}
  {}+{} \underbrace{(12p)\cdot(2p-1)}_{\bref{proj:po}} {}+{}
  \underbrace{(12p)\cdot(2p-1)}_{\bref{proj:pp}}
  \\
  {}+{} \underbrace{\sum_{s=2}^{p-1}16p\cdot(2s-1)}_{\bref{proj:so}}
  {}+{} \underbrace{\sum_{s=1}^{p-2}16p\cdot(2s+1)}_{\bref{proj:ss}}
  {}+{} \underbrace{(24p)\cdot 1}_{\bref{proj:10}}\biggr) = 64 p^4 =
  \dim\algU
\end{multline}
(the overall $4$ is for the values taken by $\alpha$ and $\beta$).

\subsection{}
As an immediate corollary of the structure of projective
$\algU$-modules, we obtain the minimal polynomials for the Casimir
elements in~\bref{sec:casimirs}, simply by finding the eigenvalues on
modules from each linkage class and taking the Jordan-cell size into
account (which is $3$ for the atypical linkage class, $2$ for each of
the typical linkage classes, and $1$ for each of the Steinberg
classes; hence the multiplicities of the minimal-polynomial roots
in~\bref{sec:casimirs}).

As a less trivial (although generally straightforward) corollary, with
find the $\algU$ center.

\section{\textbf{The $\balgU$ center}}\label{sec:center}

\begin{Thm}
  The algebra $\algU$ has a $(5p^2 - p + 4)$-dimensional center.  The
  center $\Zcenter$ decomposes into the direct sum
  \begin{equation*}
    \Zcenter=\bigoplus_{j=1}^{p^2+p-1}\oC\cdot \mathbf{e}_j \
    {}\oplus \bigoplus_{j=1}^{4p^2-2p+5}\oC\cdot \mathbf{w}_j
  \end{equation*}
  of linear subspaces generated by primitive central idempotents
  $\mathbf{e}_j$ and central nilpotents $\mathbf{w}_j$.  The block
  decomposition of the center as an associative algebra is
  \begin{equation*}
    \Zcenter = Z_{\mathrm{at}}\oplus\bigoplus_{j=1}^{(p-1)(p-2)}Z_{\mathrm{t}}^j
    \oplus\bigoplus_{j=1}^{4(p-1)}Z_{\mathrm{st}}^j
  \end{equation*}
  where $Z_{\mathrm{at}}$ corresponds to atypical $Z_{\mathrm{t}}^j$
  to typical and $Z_{\mathrm{st}}^j$ to Steinberg linkage class.
\end{Thm}

This theorem is an application of the construction of projective
$\algU$ modules.  We calculate the $\algU$ center as the center of the
basic algebra (an approach also taken for $\Ures$ in~\cite{[Arike]}).

\subsection{The basic algebra of $\algU$}\label{basic-alg}
The basic algebra is the algebra of endomorphisms of the direct sum of
projective modules taken with multiplicity 1 each.  Basic algebra
generators can be chosen as primitive idempotents and nilpotents:
(i)~each primitive idempotent $\mathbf{e}_{\repQ}$ is the projector on
a single projective module $\repQ$, and (ii)~each nilpotent
$\mathbf{w}_{\repQ,n}$ is a morphism $\repQ\to\repQ'$ defined uniquely
by the condition that it sends the top subquotient of $\repQ$ into an
isomorphic subquotient \textit{on level two} in a projective module
$\repQ'$, is an isomorphism of these subquotients, and acts by zero on
all projective modules other than~$\repQ$.  Hence,
$n=1,\dots,N_{\repQ}$, where $N_{\repQ}$ is the number of level-2
subquotients in the linkage class that are isomorphic to the top
subquotient of $\repQ$; it is in fact equal to the number of level-2
subquotients in $\repQ$.

We describe this in more detail, invoking the structure of various
projective modules:
\begin{enumerate}
\item Each of the $4(p-1)$ Steinberg modules (simple projective
  modules) contributes only an idempotent to the basic algebra
  generators.  These idempotents are central.

\item ``Typical'' projective modules contribute $4(p-1)(p-2)$
  idempotents (projectors on each of typical projective modules) and
  $8(p-1)(p-2)$ nilpotents: for each typical projective module, these
  are two maps (distinguished by $\pm$ in front of $\beta$)
  \begin{equation*}
    \repQ^{\alpha,\beta}_{s,r} \to
    \repQ^{-\alpha,\pm\beta}_{p-s,p+r-s},\quad r\neq 0,s
  \end{equation*}
  sending the top subquotient into a level-2 subquotient.

\item\label{item3-basic} For the atypical projective modules, there
  are $4(2p-1)$ idempotents.  As regards maps to level two in
  projective modules, we see from the graphs in
  Figs.~\ref{fig:proj-po}--\ref{fig:proj-10} that the relevant number
  $N=N_{\repQ}$ of such maps is as follows for the five species of
  atypical projective modules:
  \begin{equation}
    \label{BA-gens}
    \begin{alignedat}{6}
      &\repQ^{\alpha,\beta}_{p,0}\ (\zpo):&\quad& N = 3,
      &\qquad&\repQ^{\alpha,\beta}_{p-1,p-1}\ (\zpp):\quad&& N = 3,\\
      &\repQ^{\alpha,\beta}_{s,0}\ (\zso):&& N = 4,
      &&\repQ^{\alpha,\beta}_{s,s}\ (\zss):&& N = 4,
      &\qquad&\repQ^{\alpha,\beta}_{1,0}\ (\zoz):&\quad& N = 6.
    \end{alignedat}
  \end{equation}
  There are respectively $4$, $4$, $4(p-2)$, $4(p-2)$, and $4$
  projective modules of each species, which gives the total of $
  4\cdot 3 + 4\cdot 3 + 4 (p - 2)\cdot 4 + 4 (p - 2)\cdot 4 + 4\cdot 6
  = 32 p - 16$ ``atypical'' generators $\Omega_j$ of the basic
  algebra.
\end{enumerate}

\subsection{The dimension of the center}
We now list the generators of the center.  These are primitive central
idempotents (which are enumerated immediately) and central nilpotents
(finding which requires more work).  Their total number
(see~\bref{dim-summary} below) gives the dimension of the center of
$\algU$.

\subsubsection{Central idempotents and nilpotents}\label{central-e-w}
Each linkage class (see the list in~\bref{sec:link-cl}) yields a
primitive central idempotent, the projection onto that linkage class.
In addition, the linkage classes except the Steinberg ones yield
several central nilpotents each.  The items describing them are listed
below in the order of (rapidly) increasing complexity.

\begin{enumerate}
\item Each of the $4(p-1)$ Steinberg linkage classes produces a single
  central idempotent $\est^{\alpha,\beta}_r$ ($r=1,\dots,p-1$).

\item Each of the $(p-1)(p-2)$ ``typical'' linkage classes yields one
  central idempotent $\etyp^{\alpha}_{s,r}$ and four central
  nilpotents $\wtip^{\alpha,\beta}_{s,r}(\gamma)$, where $1\leq
  r<s\leq p-1$ and $\alpha,\beta,\gamma=\pm$ (note that $\alpha$, $s$,
  and $r$ enumerate a linkage class, while $\beta$ and $\gamma$ range
  over nilpotents inside a linkage class).  These
  $\wtip^{\alpha,\beta}_{s,r}(\gamma)$ are the maps
  \begin{equation*}
    \wtip^{\alpha,\beta}_{s,r}(+): \repQ^{\alpha,\beta}_{s,r} \to \repQ^{\alpha,\beta}_{s,r},\quad r\neq0,s
  \end{equation*}
  and
  \begin{equation*}
    \wtip^{\alpha,\beta}_{s,r}(-): \repQ^{-\alpha,\beta}_{p-s,p+r-s}
    \to \repQ^{-\alpha,\beta}_{p-s,p+r-s},\quad r\neq0,s,
  \end{equation*}
  sending the top subquotient to the bottom subquotient in the same
  projective module.

\item The ``atypical'' linkage class yields one central idempotent and
  central nilpotents that are of two groups: some follow immediately
  (item~a) and the derivation of others is somewhat more involved
  (item~b).
  \begin{enumerate}
  \item There are $4(2p-1)$ central nilpotents
    $\watypb^{\alpha,\beta}_{a}$, $a=1,\dots2p-1$, $\alpha,\beta=\pm$,
    one for each atypical projective module
    ($\repQ^{\alpha,\beta}_{s,0}$ with $1\leq s\leq p$ and
    $\repQ^{\alpha,\beta}_{s,s}$ with $1\leq s\leq p-1$). Each
    $\watypb^{\alpha,\beta}_{a}$ maps the top subquotient into the
    isomorphic bottom subquotient in the same projective module (and
    is zero on all other projective modules).  It then follows that
    $\watypb^{\alpha,\beta}_{a}\watypb^{\alpha',\beta'}_{a'}=0$ for
    all values of the indices.

  \item In addition, there are central nilpotents $\watypm_{b}$ given
    by linear combinations of noncentral idempotents $\Watypm_{B}$
    defined as follows.  Each $\Watypm_{B}$ is a map from the top
    subquotient in one of the atypical projective modules into an
    isomorphic subquotient at level three in the same projective
    module.  We can therefore write $\Watypm_{B}=W_{\Omega,m}$, where
    $\repQ$ is a projective module and $m$ labels its level-$3$
    subquotients that are isomorphic to its top subquotient.
    
    The graphs in Figs.~\ref{fig:proj-po}--\ref{fig:proj-10} readily
    show that the number $M$ of such subquotients and hence the number
    of $\Watypm_{B}$ yielded by each species of projective modules is
    as follows:
    \begin{alignat*}{4}
      \qquad\qquad&\repQ^{\alpha,\beta}_{p,0} (\zpo): &&M = 4,
      &\qquad&(\repQ^{\alpha,\beta}_{s,0})_{2\leq s\leq p-1}
      (\zso):\quad &&M = 4\cdot 2(p - 2),
      \\
      &\repQ^{\alpha,\beta}_{p-1,p-1} (\zpp): &\quad&M = 4,
      &&(\repQ^{\alpha,\beta}_{s,s})_{1\leq s\leq p-2} (\zss):
      &&M = 4\cdot 2 (p - 2),\\
      &\repQ^{\alpha,\beta}_{1,0} (\zoz): &&M = 4\cdot 4
    \end{alignat*}
    (a factor of $4$ in each case is of course due to the four values
    taken by $(\alpha,\beta)$).  This gives a $(16p-8)$-dimensional
    space of noncentral ``level-3'' nilpotents.

    Taking their linear combinations
    \begin{equation}\label{small-w}
      w = \sum x_{\Omega,m} \Watypm_{\Omega,m}
      = \sum_{B=1}^{16p-8}x_B \Watypm_{B},
    \end{equation}
    we require the commutativity with the basic algebra generators.
    This $w$ already commutes with all idempotents (which act either
    as identity or by zero) and, evidently, with the all ``typical''
    nilpotents.  It remains to require that it commute with the basic
    algebra generators described in item~3 in~\bref{basic-alg}
    (page~\pageref{item3-basic}):
    \begin{equation}\label{the-eqs}
      w\,\Omega_j - \Omega_j\,w = 0,\quad j=1,\dots, 32 p - 16,
    \end{equation}

    Solving Eqs.~\eqref{the-eqs} for the $x_{\Omega,m}$, we find a
    $(2p+1)$-dimensional subspace of central ``level-3'' nilpotents.
  \end{enumerate}
\end{enumerate}
\begin{thm}\label{thm:2p+1}
  There are exactly $2p+1$ linearly independent solutions
  $\watypm_{b}$ of Eqs.~\eqref{the-eqs}.
\end{thm}
We prove this in~\bref{center-proof} by deriving the explicit form of
the equations and solving them.

\subsubsection{}\label{dim-summary}
As a corollary, we calculate the dimension of the center of $\algU$ as
follows, with ``idem.'' and ``nilp.'' referring to the idempotents and
nilpotents described in items 1--3 above:
\begin{alignat*}{2}
  \dim\Zcenter &= \underbrace{4(p-1)}_{\text{idem., item~1}} &&{}+{}
  \underbrace{(p-1)(p-2)}_{\text{idem., item~2}} {}+{}
  \underbrace{1}_{\text{idem., item~3}}
  \\*
  &&&{}+{} \underbrace{4(p-1)(p-2)}_{\text{nilp., item~2}} {}+{}
  \underbrace{10p-3}_{\text{nilp., item~3}}
  \\*
  &=5p^2-p+4.
\end{alignat*}

\begin{rem}
  It follows immediately that the algebra of the $\watypm_{b}$
  (solutions of~\eqref{the-eqs}) is
  \begin{equation*}
    \watypm_{b}\watypm_{c}=\sum_{j,\alpha,\beta}
    f^{\alpha,\beta,j}_{b,c}\watypb^{\alpha,\beta}_{j}
  \end{equation*}
  where all nonzero constants $f^{\alpha,\beta,j}_{b,c}$ can be chosen
  equal to $1$ by rescaling the $\watypb^{\alpha,\beta}_{j}$.
\end{rem}

\section{\textbf{Conclusions}}
The full structure of projective modules is a powerful tool in
studying an associative algebra.  Various consequences of the
construction in this paper are to be worked out with regard to the
properties relevant for the LCFT counterpart of our
$\overline{\mathscr{U}}_{\q} s\ell(2|1)$.  The results can have a
bearing on various facets of LCFT{} models, in the range from ``spin
chains,'' potentially allowing physical applications (see
\cite{[GV],[GSaT],[GJRSV]} and the references therein), to no less
exciting ``categorial studies''
(see~\cite{[FSS-1106],[FSS-1207],[FSS-1302]} and the references
therein).  A spin chain that suggests itself in relation to our
$\overline{\mathscr{U}}_{\q} s\ell(2|1)$ is the one composed of
alternating ``fundamental'' $(s=1,r=1)$ and ``antifundamental''
$(s=2,r=0)$ $3$-dimensional representations.  A relevant tool in its
study would be Lusztig's divided-power version of the algebra.  As
regards the relation to ``continuous'' LCFT, it is sensitive to the
Hopf algebra $H$ that defines the corresponding category of
Yetter--Drinfeld modules~\cite{[nich-sl2-1]}.  As we see, $H$ comes
with its universal $R$-matrix, which also is to play a role in CFT.
Of primary interest is the modular group action on (a subalgebra in)
the $\overline{\mathscr{U}}_{\q} s\ell(2|1)$ center, which is also
linked to the study of the corresponding LCFT torus amplitudes and
their modular and other properties.  Tensor products of
$\overline{\mathscr{U}}_{\q} s\ell(2|1)$ modules conjecturally
correspond to fusion on the LCFT side.\enlargethispage{\baselineskip}

We thank B.~Feigin, A.~Gainutdinov, and A.~Kiselev for the useful
discussions.  This paper was supported in part by the RFBR grant
13-01-00386.

\appendix
\section{}
\subsection{Simple $\Ures$
  modules~\cite{[FGST]}.}\label{app:SL2-irrep}
The $\Ures$ Hopf algebra is defined in Eqs.~\eqref{Ures}.  Its simple
modules can be labeled as $\repX^{\alpha}_{s}$, where $\alpha=\pm$ and
$s=1,\dots,p$ (with $\dim\repX^{\alpha}_{s}=s$).  A basis in
$\repX^{\alpha}_{s}$ is denoted by $\ket{\alpha,s}_n$, $n=0,\dots
s-1$, with $\ket{\alpha,s}_0$ being a highest-weight vector:
\begin{gather*}
  E \ket{\alpha,s}_0=0,\quad
  K\ket{\alpha,s}_0=\alpha\q^{s-1}\ket{\alpha,s}_0.
\end{gather*}
The $\Ures$ action on all the $\ket{\alpha,s}_n$, in a ``naturally
isomorphic'' notation, is given in~\cite{[FGST]}.

\subsection{Projective $\Uresk$ modules}\label{app:Uk-proj}
The projective $\Uresk$-module $\repP^{\alpha,\beta}_{s,r}$ that
covers $\repX^{\alpha,\beta}_{s,r}$ has a basis
$\ket{\alpha,s,\beta,r}^{\mathsf{b}}_n$,
$\ket{\alpha,s,\beta,r}^{\mathsf{a}}_n$ with $n=0,\dots,s-1$ and
$\ket{\alpha,s,\beta,r}^{\mathsf{x}}_m$,
$\ket{\alpha,s,\beta,r}^{\mathsf{y}}_m$ with $m=0,\dots,p-s-1$. The
action of $\Ures$ in this basis is given in \cite{[FGST]}, and the
action of $k$ is
\begin{alignat*}{2}
  k\ket{\alpha,s,\beta,r}^{\mathsf{a}}_n
  &=\beta\q^{-r+n}\ket{\alpha,s,\beta,r}^{\mathsf{a}}_n, \quad&
  k\ket{\alpha,s,\beta,r}^{\mathsf{b}}_n
  &=\beta\q^{-r+n}\ket{\alpha,s,\beta,r}^{\mathsf{b}}_n,
  \\
  k\ket{\alpha,s,\beta,r}^{\mathsf{x}}_n
  &=\beta\q^{s-r+n-p}\ket{\alpha,s,\beta,r}^{\mathsf{x}}_n, \quad&
  k\ket{\alpha,s,\beta,r}^{\mathsf{y}}_n
  &=\beta\q^{s-r+n}\ket{\alpha,s,\beta,r}^{\mathsf{y}}_n.
\end{alignat*}

\subsection{The rest of the $\algU$ action on simple
  modules}\label{action-on-simple}
For completeness, we give the formulas that fully define the $\algU$
action on simple modules in~\bref{simple-list} (the action of
fermionic generators on basis vectors was already described there).

For the atypical $\repZ^{\alpha,\beta}_{s,0}$ modules
in~\eqref{Zs0-decomp}, we have
\begin{alignat*}{2}
  K\ZC{\alpha,s,\beta,0}_n&=\alpha\q^{s-2n-1}\ZC{\alpha,s,\beta,0}_n,&
  k\ZC{\alpha,s,\beta,0}_n&=\beta\q^{n}\ZC{\alpha,s,\beta,0}_n,\\
  F \ZC{\alpha,s,\beta,0}_n&=\ZC{\alpha,s,\beta,0}_{n+1},&
  E \ZC{\alpha,s,\beta,0}_n&=\alpha \qint{n}\qint{s - n}\ZC{\alpha,s,\beta,0}_{n-1},\\
  K\ZB{\alpha,s,\beta,0}_m&=\alpha\q^{s-2m-2}\ZB{\alpha,s,\beta,0}_m,&
  k\ZB{\alpha,s,\beta,0}_m&=-\beta\q^{m+1}\ZB{\alpha,s,\beta,0}_m,\\
  F \ZB{\alpha,s,\beta,0}_m&=\ZB{\alpha,s,\beta,0}_{m+1},& E
  \ZB{\alpha,s,\beta,0}_m&=\alpha \qint{m}\qint{s\!-\!m\!-\!1}
  \ZB{\alpha,s,\beta,0}_{m-1},
\end{alignat*}
where we set $\ZC{\alpha,s,\beta,0}_n=0$ for $n<0$ and $n>s-1$, and
$\ZB{\alpha,s,\beta,0}_m=0$ for $m<0$ and $m>s-2$.  For $s=1$, each
module $\repZ^{\alpha,\beta}_{1,0}$ with $\alpha,\beta=\pm$ is
$1$-dimensional.  The basis consists of a single vector
$\ZC{\alpha,1,\beta,0}_0$ such that
\begin{equation*}
  K\ZC{\alpha,1,\beta,0}_0=\alpha\ZC{\alpha,1,\beta,0}_0,
  \qquad
  k\ZC{\alpha,1,\beta,0}_0=\beta\ZC{\alpha,1,\beta,0}_0,
\end{equation*}
with the other generators acting trivially.  For convenience in what
follows (in the cases where $\repZ^{\alpha,\beta}_{1,0}$ occur along
with higher-dimensional modules), we set $\ZB{\alpha,1,\beta,0}_m=0$
for all~$m$.

For the atypical $\repZ^{\alpha,\beta}_{s,s}$ modules
in~\eqref{Zss-decomp}, we have
\begin{alignat*}{2}
  K\ZC{\alpha,s,\beta,s}_n&=\alpha\q^{s-2n-1}\ZC{\alpha,s,\beta,s}_n,&
  k\ZC{\alpha,s,\beta,s}_n&=\beta\q^{-s+n}\ZC{\alpha,s,\beta,s}_n,\\
  F \ZC{\alpha,s,\beta,s}_n&=\ZC{\alpha,s,\beta,s}_{n+1},&
  E \ZC{\alpha,s,\beta,s}_n&=\alpha \qint{n}\qint{s - n}\ZC{\alpha,s,\beta,s}_{n-1},\\
  K\ZB{\alpha,s,\beta,s}_m&=\alpha\q^{s-2m}\ZB{\alpha,s,\beta,s}_m,&
  k\ZB{\alpha,s,\beta,s}_m&=-\beta\q^{-s+m}\ZB{\alpha,s,\beta,s}_m,\\
  F \ZB{\alpha,s,\beta,s}_m&=\ZB{\alpha,s,\beta,s}_{m+1},& E
  \ZB{\alpha,s,\beta,s}_m&=\alpha \qint{m}\qint{s -
    m+1}\ZB{\alpha,s,\beta,s}_{m-1},
\end{alignat*}
where we set $\ZC{\alpha,s,\beta,s}_n=0$ for $n<0$ and $n>s-1$, and
$\ZB{\alpha,s,\beta,s}_m=0$ for $m<0$ and $m>s$.

For the typical modules $\repZ^{\alpha,\beta}_{s,r}$
in~\eqref{bulk-decomp},
\begin{alignat*}{2}
  K \ZC{\alpha, s, \beta, r}_{j} &= \alpha \q^{s - 2 j - 1}
  \ZC{\alpha, s, \beta, r}_{j}, & k \ZC{\alpha, s, \beta, r}_{j} &=
  \beta \q^{-r + j} \ZC{\alpha, s, \beta, r}_{j},
  \\
  F \ZC{\alpha, s, \beta, r}_{j} &= \ZC{\alpha, s, \beta, r}_{j + 1},
  & E \ZC{\alpha, s, \beta, r}_{j} &= \alpha \qint{j} \qint{s - j}
  \ZC{\alpha, s, \beta, r}_{j - 1},
  \\
  K \ZU{\alpha, s, \beta, r}_{m} &= \alpha \q^{s - 2 m} \ZU{\alpha, s,
    \beta, r}_{m}, & k \ZU{\alpha, s, \beta, r}_{m} &= -\beta \q^{-r +
    m} \ZU{\alpha, s, \beta, r}_{m},
  \\
  F \ZU{\alpha, s, \beta, r}_{m} &= \ZU{\alpha, s, \beta, r}_{m + 1},
  & E \ZU{\alpha, s, \beta, r}_{m} &= \alpha \qint{m} \qint{s - m + 1}
  \ZU{\alpha, s, \beta, r}_{m - 1},
  \\
  K \ZD{\alpha, s, \beta, r}_{n} &= \alpha \q^{s - 2 n - 2}
  \ZD{\alpha, s, \beta, r}_{n}, & k \ZD{\alpha, s, \beta, r}_{n} &=
  -\beta \q^{-r + n + 1} \ZD{\alpha, s, \beta, r}_{n},
  \\
  F \ZD{\alpha, s, \beta, r}_{n} &= \ZD{\alpha, s, \beta, r}_{n + 1},
  & E \ZD{\alpha, s, \beta, r}_{n} &= \alpha \qint{n} \qint{s - n - 1}
  \ZD{\alpha, s, \beta, r}_{n - 1},
  \\
  K \ZB{\alpha, s, \beta, r}_{j} &= \alpha \q^{s - 2 j - 1}
  \ZB{\alpha, s, \beta, r}_{j}, & k \ZB{\alpha, s, \beta, r}_{j} &=
  \beta \q^{-r + j + 1} \ZB{\alpha, s, \beta, r}_{j},
  \\
  F \ZB{\alpha, s, \beta, r}_{j} &= \ZB{\alpha, s, \beta, r}_{j + 1},
  & E \ZB{\alpha, s, \beta, r}_{j} &= \alpha \qint{j} \qint{s - j}
  \ZB{\alpha, s, \beta, r}_{j - 1},
\end{alignat*}
where we as usual assume that the vectors outside the ranges specified
in~\eqref{four-sub} are equal to zero.  If $s=1$, the decomposition
degenerates to $\repZ^{\alpha,\beta}_{1,r}=\repX^{\alpha,\beta}_{1,r}
\oplus\repX^{\alpha,-\beta}_{2,r}\oplus\repX^{\alpha,\beta}_{1,r-1}$,
with a basis $\ZC{\alpha,1,\beta,r}_0$,
$\bigl(\ZU{\alpha,1,\beta,r}_m\bigr)_{m=0,1}$,
$\ZB{\alpha,1,\beta,r}_0$, and with the $\algU$ action easily
deducible from the general case above.

For the Steinberg modules $\repZ^{\alpha,\beta}_{p,r}$
in~\eqref{stein-decomp},
\begin{alignat*}{2}
  K \ZC{\alpha, p, \beta, r}_{n} &= \alpha \q^{-2 n - 1} \ZC{\alpha,
    p, \beta, r}_{n}, &\quad k \ZC{\alpha, p, \beta, r}_{n} &= \beta
  \q^{n - r} \ZC{\alpha, p, \beta, r}_{n},
  \\
  F \ZC{\alpha, p, \beta, r}_{n} &= \ZC{\alpha, p, \beta, r}_{n + 1},
  & E \ZC{\alpha, p, \beta, r}_{n} &= -\alpha \qint{n}^2 \ZC{\alpha,
    p, \beta, r}_{n - 1},
  \\
  K \ZU{\alpha, p, \beta, r}_{n} &= \alpha \q^{-2 n} \ZU{\alpha, p,
    \beta, r}_{n}, & k \ZU{\alpha, p, \beta, r}_{n} &= -\beta \q^{n -
    r} \ZU{\alpha, p, \beta, r}_{n},
  \\
  F \ZU{\alpha, p, \beta, r}_{n} &= \ZU{\alpha, p, \beta, r}_{n + 1},
  & E \ZU{\alpha, p, \beta, r}_{n} &= -\alpha \qint{n} \qint{n - 1}
  \ZU{\alpha, p, \beta, r}_{n - 1},
  \\
  K \ZD{\alpha, p, \beta, r}_{n} &= \alpha \q^{-2 n - 2} \ZD{\alpha,
    p, \beta, r}_{n}, & k \ZD{\alpha, p, \beta, r}_{n} &= -\beta \q^{n
    + 1 - r} \ZD{\alpha, p, \beta, r}_{n},
  \\
  F \ZD{\alpha, p, \beta, r}_{n} &= \ZD{\alpha, p, \beta, r}_{n + 1},
  & E \ZD{\alpha, p, \beta, r}_{n} &= \alpha \ZU{\alpha, p, \beta,
    r}_{n}
  \\*
  &&&\quad{}-\alpha \qint{n} \qint{n + 1}\ZD{\alpha, p, \beta, r}_{n -
    1},
  \\
  K \ZB{\alpha, p, \beta, r}_{n} &= \alpha \q^{-2 n - 1} \ZB{\alpha,
    p, \beta, r}_{n}, & k \ZB{\alpha, p, \beta, r}_{n} &= \beta \q^{n
    + 1 - r} \ZB{\alpha, p, \beta, r}_{n},
  \\
  F \ZB{\alpha, p, \beta, r}_{n} &= \ZB{\alpha, p, \beta, r}_{n + 1},
  & E \ZB{\alpha, p, \beta, r}_{n} &= -\alpha \qint{n}^2 \ZB{\alpha,
    p, \beta, r}_{n - 1}.
\end{alignat*}

\section{Proofs and calculation details}
\subsection{The universal $R$-matrix of $\algU$}\label{app:R}
We calculate the universal $R$-matrix for $\algU$ by relating $\algU$
to the Drinfeld double of~$U_{\leq}$ (see
decomposition~\eqref{UU-decomp}).

\subsubsection{}\label{sec:duality} The generators $E$ and $C$ can be
viewed as functionals on the subalgebra $U_{\leq}$ such that
\begin{equation*}
  \eval{E}{F K^i k^j} := -\ffrac{1}{\q - \q^{-1}},
  \qquad
  \eval{C}{B K^i k^j} := \ffrac{1}{\q - \q^{-1}},
\end{equation*}
and all other evaluations of $E$ and $C$ on the PBW basis elements in
$U_{\leq}$ vanish.  We let $\tilde{U}_{>}$ temporarily denote the
algebra generated by the functionals $E$ and $C$ with the product of
any two functionals $\beta$ and $\gamma$ defined standardly as
\begin{equation*}
  \eval{\beta\gamma}{x}=
  \eval{\beta}{x'}\eval{\gamma}{x''},
\end{equation*}
where $\Delta x=x'\tensor x''$ is the coproduct on $U_{\leq}$.  It
then readily follows that the relations $ E E C - \qint{2} E C E + C E
E = 0$ and $C C = 0$ hold in $\tilde{U}_{>}$.  Hence, polynomials in $ E
$ and $ C $ can be expressed linearly in terms of $E^{i}$, $E^{i + 1}
C$, $E^{i} C E$, $E^{i} C E C$, $i\geq 0$.  By induction on the power
of the relevant generator, we then find that these tentative PBW basis
elements have the following nonzero evaluations on the PBW basis
elements in $U_{\leq}$:
\begin{align*}
  \eval{E^{n}}{F^{n} K^i k^{j}} &=(-1)^{n} (\q - \q^{-1})^{-n}
  \q^{\half n (n - 1)} \qfac{n},
  \\
  \eval{E^{n + 1} C}{F^{n} B F K^i k^{j}} &= (-1)^{n + 1} (\q -
  \q^{-1})^{-n - 2} \q^{\half(n + 2) (n - 1)} \qfac{n + 1},
  \\
  \eval{E^{n} C}{F^{n} B K^i k^{j}} &= (-1)^{n} (\q - \q^{-1})^{-n -
    1} \q^{\half n (n - 1)} \qfac{n},
  \\
  \eval{E^{n} C E}{F^{n} B F K^i k^{j}} &= (-1)^{n + 1} (\q -
  \q^{-1})^{-n - 2} \q^{\half n (n - 1)}\bigl(1 + \q^{n -
    1}\qint{n}\bigr) \qfac{n},
  \\
  \eval{E^{n} C E}{F^{n + 1}B K^i k^{j}} &= (-1)^{n + 1} (\q -
  \q^{-1})^{-n - 2} \q^{\half(n + 2) (n - 1)} \qfac{n + 1},
  \\
  \eval{E^{n} C E C}{F^{n} B F B K^i k^{j}} &= (-1)^{n + 1} (\q -
  \q^{-1})^{-n - 2} \q^{\half(n + 1) (n - 2)} \qfac{n}
\end{align*}
(and all other evaluations vanish).  It then follows that $E^p=0$ in
$\tilde{U}_{>}$, and $\tilde{U}_{>}$, with its PBW basis modeled on that
in $U_{>}$, is isomorphic to $U_{>}$ as an algebra.

To diagonalize the ``nondiagonal'' part of the pairing $U_{>}\tensor
U_{\leq}\to\oC$, we define
\begin{align*}
  X_n &= \qint{n} E^{n - 1} C E - \bigl(\q^{2 - n} +
  \qint{n-1}\bigr)E^{n} C,
  \\
  Y_n &= E^{n - 1} C E - \q^{-1} E^{n} C.
\end{align*}
Then
\begin{alignat*}{2}
  \eval{X_n}{F^{n}B K^i k^j} & =(-1)^{n + 1} \ffrac{\q^{\half(n - 1)
      (n - 2)}\qfac{n}}{(\q - \q^{-1})^{n}}, \kern-60pt&
  \eval{X_n}{F^{n - 1}BF K^i k^j} ={}0&,
  \\
  \eval{Y_n}{F^{n}B K^i k^j} &=0, & \eval{Y_n}{F^{n - 1}BF K^i k^j} =
  (-1)^n \ffrac{\q^{\half n (n - 3)}\qfac{n - 1}}{(\q -
    \q^{-1})^{n}}&.
\end{alignat*}

\subsubsection{}\label{app:DD-rel} We next extend $\tilde{U}_{>}$ by
generators $L,\ell\in H^*$, also functionals on $U_{\leq}$\,, that we
require to commute with the generators of $\Uminus$ exactly as $K$ and
$k$ do:
\begin{alignat*}{2}
  F L &= \q^2 L F,&\qquad B L &= \q^{-1} L B,
  \\
  F \ell &= \q^{-1} \ell F,&\qquad B \ell &= -\ell B,
\end{alignat*}
where for $a\in U_{\leq}$ and a functional $\beta$, we evaluate the
product using the Drinfeld-double formula
\begin{equation}\label{D-d}
  a\beta = (a'\leftact \beta\rightact S^{-1}(a'''))a''
\end{equation}
(where $\Delta a = a'\tensor a''$ is the coproduct on $U_{\leq}$).  It
then follows that
\begin{equation*}
  \eval{L}{K^m k^n}=\q^{n-2m},\qquad
  \eval{\ell}{K^m k^n}=(-1)^n \q^m
\end{equation*}
(with the other evaluations vanishing).  Using~\eqref{D-d}, we then
establish further ``cross-commutator'' relations:
\begin{equation*}
  F E - E F = -\ffrac{L - K^{-1}}{\q - \q^{-1}},
  \qquad
  B C - C B = -\ffrac{\ell - k^{-1}}{\q - \q^{-1}},
\end{equation*}
as well as $F C - C F = 0$ and $B E - E B = 0$ (and $KL = L K$, etc.).

This essentially (modulo easily reconstructible details) shows that
$\algU$ is the quotient of the Drinfeld double of $U_{\leq}$ by the
Hopf ideal generated by the relations $L = K$ and $\ell = k$.  The
(inverse) universal $R$-matrix is then inherited from the Drinfeld
double in the form
\begin{multline*}
  R^{-1}= \sum_{a=0}^{p - 1} (-1)^{a} \ffrac{\q^{-\half a (a - 1)}(\q
    - \q^{-1})^{a}}{\qfac{a}} \Bigl( E^{a}\tensor F^{a} - \q^{a - 1}
  X_{a}\tensor F^{a} B
  \\
  \qquad{}- \q (\q - \q^{-1}) Y_{a + 1}\tensor F^{a} B F - \q (\q -
  \q^{-1})^2 E^{a} C E C\tensor F^{a} B F B \Bigr)\rho^{-1},
\end{multline*}
where
\begin{equation*}
  \rho^{-1}=\fffrac{1}{(2 p)^2}\sum_{i=0}^{2 p - 1}\sum_{j=0}^{2 p - 1}
  \sum_{m=0}^{2 p - 1}\sum_{n=0}^{2 p - 1}
  (-1)^{j n} \q^{2 i m - i n - j m} K^{i} k^{j} \tensor
  K^{m} k^{n}.
\end{equation*}
From here,
the universal $R$-matrix $R=(S^{-1}\tensor\id)R^{-1}$ (where $S$ is
the antipode of $\algU$) follows in the form
$R =\rho\bar{R}$, with $\rho$ in~\eqref{the-rho} and
\begin{align*}\label{the-R-1}
  \bar{R} &= \sum_{a=0}^{p - 1} \ffrac{\q^{\half a (a - 1)}(\q -
    \q^{-1})^{a}}{\qfac{a}} \Bigl(E^{a}\tensor F^{a}
  \\
  &\qquad{} - \q^{-1} \bigl(\q^{2 - a} \qint{a} E^{a - 1} C E - (1 +
  \q^{2 - a} \qint{a - 1}) E^{a} C\bigr)\tensor F^{a} B
  \\
  &\qquad{}+ \q^{-2} (\q^2 - 1) \bigl(E^{a} C E - \q E^{a + 1}
  C\bigr)\tensor F^{a} B F - \q^{-3} (\q^2 - 1)^2 E^{a} C E C \tensor
  F^{a}B F B\Bigr),
\end{align*}
which can be readily rewritten in the factored form as
in~\eqref{barR-factored}.

\subsection{Proof of~\bref{prop:M-matrix}}\label{M-calc}
The M-matrix (see~\bref{sec:M}) can obviously be written as
$M=(\rho_{21}\bar{R}_{21}\rho)\bar{R}$, where we calculate
$\rho_{21}\bar{R}_{21}\rho$ from the definition of $\rho$ and
$\bar{R}$:
\begin{align*}
  \rho_{21}\bar{R}_{21}\rho&= \ffrac{1}{(2 p)^4} \sum_{i=0}^{2 p -
  1}\sum_{j=0}^{2 p - 1} \sum_{m=0}^{2 p - 1}\sum_{n=0}^{2 p - 1}
  \sum_{i'=0}^{2 p - 1}\sum_{j'=0}^{2 p - 1} \sum_{m'=0}^{2 p -
  1}\sum_{n'=0}^{2 p - 1} (-1)^{j n + j' n'}
  \\
  &\quad\times \q^{-2 i m + j m + i n} \q^{-2 i' m' + j' m' + i'
  n'} (K^{m} k^{n}\tensor K^{i} k^{j}) \bar{R}_{21} (K^{i'}
  k^{j'}\tensor K^{m'} k^{n'}).
\end{align*}
Here, we next move all $K$ and $k$ factors to the right of
$\bar{R}_{21}$; after simple changes of summation variables, $i'
\to i' - m$ and $j' \to j' - n$, we then see that the summations
over $m$ and $n$ can give a nonzero result only for terms with
particular values of $i$ and $j$, e.g., $i = m' - a$ and $j = n'$
or $i = m' - 1 - a$ and $j = n' -1$, etc.\ (depending on the term
taken in the expression for $\bar{R}_{21}$).  This reduces the
$m'$ and $n'$ summations to the form
\begin{equation*}
  \sum_{m'=0}^{2p-1}\sum_{n'=0}^{2p-1}
  q^{x m' + y n'}K^{2m'}k^{2n'},
\end{equation*}
where $x$ and $y$ are integer linear combinations of other
summation indices.  Next, splitting the range $[0,2p-1]$ of $m'$
into $[0,p-1]\cup[p,2p-1]$ and shifting $m'$ by $p$ in the second
half of the range, we obtain that the entire sum is proportional
to $1+(-1)^x$, which acts as a selection rule for the parity of
one of the remaining summation indices involved in~$x$.  The same
is repeated for the $n'$ sum, yielding the factor $1+(-1)^y$ and
another selection rule.  The result of these straightforward
manipulations is
\begin{align*}
  \rho_{21}\bar{R}_{21}\rho&= \!\sum_{a=0}^{p - 1}\!
  \ffrac{\q^{-\frac{1}{2}a - \frac{3}{2} a^2}(\q -
  \q^{-1})^a}{\qfac{a}}
  \\
  &\quad{}\times\!  \biggl(\!\!  F^{a}\tensor E^{a} + \q^{2 a - 1}
  F^{a}B\tensor \bigl(\q^{2 - a} \qint{a} E^{a - 1}C E - (1\! +
  \!\q^{2 - a} \qint{a\!-\!1}) E^{a}C\bigr) (k\tensor k^{-1})
  \\
  &\qquad{}- \q^{-2 a - 2} (\q^2 - 1) F^{a}B F\tensor \bigl(E^{a}C E
  - \q E^{a + 1}C\bigr) (K k\tensor K^{-1} k^{-1})
  \\
  &\qquad{} - \q^{-1} (\q^2 - 1)^2 (F^{a}B F B\tensor E^{a}C E C) (K
  k^2\tensor K^{-1} k^{-2}) \!\!\biggr)
  \\
  &\quad{}\times \ffrac{1}{p^2}\sum_{i=0}^{p - 1}\sum_{j=0}^{p -
  1} \sum_{i'=0}^{p - 1}\sum_{j'=0}^{p - 1} \q^{2 i' j + 2 i j' -
  4 i i'} \Bigl(K^{2 i + a} k^{2 j}\tensor K^{2 i' - a} k^{2
  j'}\Bigr).
\end{align*}
This is immediately verified to be equal to $\bar{M}\bar{\rho}$
in~\bref{prop:M-matrix}.

\subsection{Coincidence of two Drinfeld maps}\label{app:MM}
Let $A$ be a Hopf algebra and $\Phi$ an invertible normalized
two-cocycle, i.e., an invertible element $\Phi\in A\tensor A$ such
that
\begin{equation*}
  \Phi_1' F_1 \tensor \Phi_1'' F_2 \tensor \Phi_2
  = \Phi_1 \tensor \Phi_2' F_1 \tensor \Phi_2'' F_2
  \quad\text{and}\quad
  \varepsilon(\Phi_1)\Phi_2=\Phi_1\varepsilon(\Phi_2)=1,
\end{equation*}
where $\Phi=\Phi_1\tensor\Phi_2=F_1\tensor F_2$ (and $\varepsilon$ is
the counit).  This standardly defines a new Hopf algebra
structure---the one with the same product and counit, and with the
coproduct and antipode given by
\begin{equation*}
  \tDelta(x)=\Phi^{-1}\Delta(x)\Phi,
  \qquad
  \widetilde{S}(x)=U^{-1} S(x) U
  \qquad \forall x\in A,
\end{equation*}
where
\begin{equation}\label{eq:U}
  U=S(\Phi_1)\Phi_2.
\end{equation}
We also note that $\widetilde{S}^2(x)=\xi^{-1}S^2(x)\xi$, where
\begin{equation}\label{xi}
  \xi=S(U^{-1}) U.
\end{equation}

\begin{thm}\label{thm:MM}
  Let $A$ be a quasitriangular Hopf algebra and $\Phi$ an invertible
  normalized 2-cocycle.  Then diagram~\eqref{MM-diag} is commutative,
  i.e.,
  \begin{equation*}
    \widetilde{\beta}(\widetilde{M}_1)\widetilde{M}_2
    =\beta(M_1)M_2
  \end{equation*}
  for any $\beta\in\Ch$.
\end{thm}

\subsubsection*{Proof}
First, we have a linear space isomorphism $\Ch\to\widetilde{\Ch}$
given by $\beta\mapsto(\beta\rightact\xi)$.  Indeed, if $\beta\in\Ch$,
which amounts to the condition that $\beta(xy)=\beta(S^2(y)x)$ for all
$x,y\in A$, then the functional $\widetilde{\beta}:x\mapsto \beta(\xi
x)$ is invariant under the ``tilded'' coadjoint action, i.e.,
$\widetilde{\beta}(xy)=\widetilde{\beta}(\widetilde{S}^2(y)x)$ for all
$x,y\in A$.

We next note two simple consequences of the cocycle condition:
\begin{align}\label{eq:cocycle1}
  \varphi_1\tensor\Phi_1\varphi_2\tensor\Phi_2
  &=\varphi_1\Phi_1'\tensor \varphi_2'\Phi_1''\tensor
  \varphi_2''\Phi_2,
  \\
  \label{eq:cocycle2}
  F_1\tensor F_2 \varphi_1\tensor \varphi_2 &=\varphi_1' F_1\tensor
  \varphi_1'' F_2'\tensor \varphi_2 F_2'',
\end{align}
where $\Phi^{-1}=\varphi_1\tensor\varphi_2$.
From~\eqref{eq:cocycle1}, applying the antipode and multiplying, we
find the identity $\varphi_1 S(\Phi_1\varphi_2)\Phi_2=
1$, whence, for $U$ in~\eqref{eq:U}, it follows that $U^{-1}=\varphi_1
S(\varphi_2)$.  Hence, $\xi=S(U^{-1}) U=S^2(f_2)S(f_1)S(F_1)F_2$, and
we calculate
\begin{align*}
  \widetilde{\beta}(\widetilde{M}_1)\widetilde{M}_2
  &=\beta\bigl(\xi\varphi_1 M_1 \Phi_1\bigr)\varphi_2 M_2 \Phi_2
  \\[-2pt]
  &=\beta\bigl(S(f_1)S(F_1)F_2\varphi_1 M_1 \Phi_1 f_2\bigr) \varphi_2
  M_2 \Phi_2
  \\[-2pt]
  \intertext{where we next note that
    $S(F_1)F_2\varphi_1\tensor\varphi_2 = S(F_1)F_2'\tensor F_2''$ as
    a simple consequence of~\eqref{eq:cocycle2}, and we can therefore
    continue}%
  &=\beta\bigl(S(f_1)S(F_1)F_2' M_1 \Phi_1 f_2\bigr)F_2'' M_2 \Phi_2
  \\[-2pt]
  &=\beta\bigl(S(f_1)S(F_1)M_1 F_2' \Phi_1 f_2\bigr)M_2 F_2'' \Phi_2
  \\[-2pt]
  \intertext{(by the property $M\Delta(x)=\Delta(x)M$ of the
    M-matrix), which after directly applying the cocycle condition
    becomes} &=\beta\bigl(S(f_1)S(\Phi_1' F_1)M_1 \Phi_1'' F_2
  f_2\bigr)M_2 \Phi_2
  \\[-2pt]
  &=\beta\bigl(S(\Phi_1')M_1 \Phi_1''\bigr)M_2 \Phi_2
  \\[-2pt]
  &=\beta\bigl(M_1 \Phi_1'' S^{-1}(\Phi_1')\bigr)M_2 \Phi_2.
\end{align*}
This is the same as $\beta\bigl(M_1\bigr)M_2$.

\subsection{Proof of~\bref{thm:2p+1}}\label{center-proof}
We solve Eqs.~\eqref{the-eqs}. \ We recall that $\Omega_j$ are the $32
p - 16$ nilpotent basic algebra generators defined in
item~\ref{item3-basic} in~\bref{basic-alg}, and the unknowns
$x_{\Omega,m}$ (see~\eqref{small-w}) are associated with the $16p-8$
``level-3'' nilpotent basic algebra elements $\Watypm_{\Omega,m}$
defined in item~3b in~\bref{central-e-w}.

For a projective module $\repQ=\repQ^{\alpha,\beta}_{r,s}$, we let
$\Watypm_{\repQ^{\alpha,\beta}_{r,s},m}$ be denoted as
$\Watypm^{\alpha,\beta}_{r,s}(m)$.  The corresponding unknowns,
accordingly, are then written as $x^{\alpha,\beta}_{r,s}(m)$, and we
moreover drop the uninformative ``$x$'' and distinguish the variables
pertaining to the five species of projective modules
in~\bref{proj:po}--\bref{proj:10} by an individual letter each:
\begin{equation}
  x^{\alpha,\beta}_{r,s}(m)=
  \begin{cases}
    \cpo^{\alpha,\beta}_p(m),& s=p,\ r=0\quad(\zpo),
    \\
    \cpp^{\alpha,\beta}_{p-1}(m),& s=r=p-1\quad(\zpp),
    \\
    \cso^{\alpha,\beta}_s(m),&2\leq s\leq p-1,\ r=0\quad(\zso),
    \\
    \css^{\alpha,\beta}_s(m),&1\leq s=r\leq p-2\quad(\zss),
    \\
    \coz^{\alpha,\beta}_1(m),& s=1,\ r=0\quad(\zoz).
  \end{cases}
\end{equation}
We recall that the argument $m$ labels level-$3$ subquotients in a
given projective module that are isomorphic to the top subquotient.
It is convenient, for uniformity, to let $m$ take not consecutive
values (e.g., from $1$ to $4$ for the projective module
$\repQ^{\alpha,\beta}_{1,0}$ in Fig.~\ref{fig:proj-10}) but the values
that the relevant subquotients are already assigned in the graphs
(which are $2$, $4$, $5$, and $7$ in Fig.~\ref{fig:proj-10}).  The
argument in $\Watypm^{\alpha,\beta}_{r,s}(m)$ is then of course
understood in the same way.

\begin{figure}[tbh]\footnotesize
  \centering
  \begin{equation*}
    \mbox{}\kern-20pt\xymatrix@R=40pt@C=5pt{
      &&&&\mysmallmatrix{\alpha&\beta\\s&r}\ar[dllll]\ar[dll]\ar[d]|{c_1}\ar[dr]|{c_2}\ar@{{.}{.}{>}}@/^10pt/^{1}[drrrrrrr]&&&&&&&&\mysmallmatrix{\alpha'&\beta'\\s'&r'}\ar[dl]\ar[dr]\ar[drrr]&&
      \\
      \mbox{}&\dots&\mbox{}&\dots&\mysmallmatrix{\alpha_1&\beta_1\\s_1&r_1}\ar@{{.}{.}{>}}@/_50pt/^(.3){c'_1/c_1}[drrrrrrr]&\mysmallmatrix{\alpha_2&\beta_2\\s_2&r_2}\ar@{{.}{.}{>}}@/^5pt/^{c'_2/c_2}[drrrrr]
      &&&&&&\mysmallmatrix{\alpha&\beta\\s&r}\ar[dl]|{c'_2}\ar[d]|{c'_1}\ar[drr]\ar[drrrr]&\dots&\mbox{}&\dots&\mbox{}
      \\
      &&&&&&&&&&\mysmallmatrix{\alpha_2&\beta_2\\s_2&r_2}&\mysmallmatrix{\alpha_1&\beta_1\\s_1&r_1}&\dots&\mbox{}&\dots&
    }
  \end{equation*}
  \caption[Mapping projective modules]{\small Mapping projective
    modules by an element of the basic algebra.  The top subquotient
    of the projective module on the left is sent into an isomorphic
    level-$2$ subquotient in the projective module on the right.
    Isomorphic descendants are then mapped into one another with the
    coefficients given by ratios of the weights associated with edges
    of the graphs.  Those children of
    $\mysmallmatrix{\alpha&\beta\\s&r}$ on the left that have no
    isomorphic subquotients among the children of
    $\mysmallmatrix{\alpha&\beta\\s&r}$ on the right are mapped to
    zero.  The procedure continues similarly to the lower-lying
    levels.}
  \label{fig:mod-map}
\end{figure}
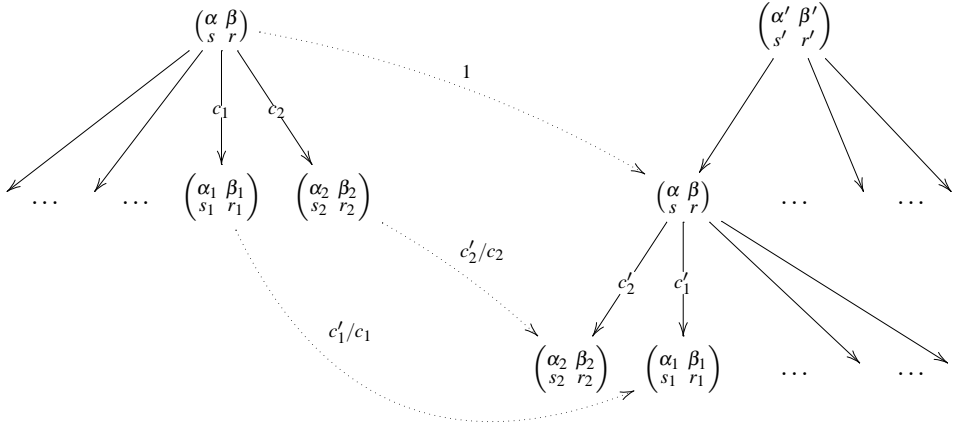

Each $\Watypm^{\alpha,\beta}_{r,s}(m)$, as well as each $\Omega_j$
in~\eqref{the-eqs}, can be regarded as a linear operator on the vector
space with basis consisting of all simple subquotients of all
projective modules, and is therefore completely determined by the
coefficients with which it sends each subquotient into (linear
combinations of) others.  By definition, the top subquotient is mapped
(into a level-$3$ subquotient by each $\Watypm$ and a level-$2$
subquotient by each $\Omega_j$) with the coefficient~$1$, while all
other coefficients are obtained from the graphs in
Figs.~\ref{fig:proj-po}--\ref{fig:proj-10}, simply from the condition
that the basic algebra elements be $\algU$ intertwiners.  This is
illustrated in Fig.~\ref{fig:mod-map}.

The commutativity in~\eqref{the-eqs} then becomes the commutativity
condition for the corresponding matrices.  We illustrate this with
equations involving the six basic algebra generators
$\Omega_{\repQ^{\alpha,\beta}_{1,0},n}$, $1\leq n\leq 6$, sending the
top one-dimensional subquotient of $\repQ^{\alpha,\beta}_{1,0}$ into
isomorphic level-2 subquotients (see~\bref{basic-alg},
item~\ref{item3-basic}).  The isomorphic subquotients occur in the
projective modules $\repQ^{-\alpha,\beta}_{p,0}$,
$\repQ^{-\alpha,-\beta}_{p,0}$, $\repQ^{-\alpha,\beta}_{p-1,p-1}$,
$\repQ^{-\alpha,-\beta}_{p-1,p-1}$, $\repQ^{\alpha,-\beta}_{2,0}$, and
$\repQ^{\alpha,\beta}_{1,1}$; we select the first one in this
(arbitrary) ordering.  Because the corresponding basic algebra
generator $\Omega_{\repQ^{\alpha,\beta}_{1,0},1}$ acts by zero on all
projective modules except $\repQ^{\alpha,\beta}_{1,0}$, commutator
equation~\eqref{the-eqs} takes the form
\begin{multline}\label{op-eq}
  \Omega_{\repQ^{\alpha,\beta}_{1,0},1}
  \Bigl(\!\coz^{\alpha,\beta}_1\vmylabel{s5}{3}{2}\Watypm^{\alpha,\beta}_{1,0}\vmylabel{s5}{3}{2}
  +\coz^{\alpha,\beta}_1\vmylabel{mt4}{3}{4}\Watypm^{\alpha,\beta}_{1,0}\vmylabel{mt4}{3}{4}
  +\coz^{\alpha,\beta}_1\vmylabel{s6}{3}{5}\Watypm^{\alpha,\beta}_{1,0}\vmylabel{s6}{3}{5}+\\
  +\coz^{\alpha,\beta}_1\vmylabel{s7}{3}{7}\Watypm^{\alpha,\beta}_{1,0}\vmylabel{s7}{3}{7}\!\Bigr)
  -\cpo^{-\alpha,\beta}_{p}\vmylabel{Cc}{3}{0}\Watypm^{-\alpha,\beta}_{p,0}\vmylabel{Cc}{3}{0}\Omega_{\repQ^{\alpha,\beta}_{1,0},1}=0.
\end{multline}

This equation for the unknowns
$\coz^{\alpha,\beta}_1\vmylabel{s5}{3}{2}$,
$\coz^{\alpha,\beta}_1\vmylabel{mt4}{3}{4}$,
$\coz^{\alpha,\beta}_1\vmylabel{s6}{3}{5}$,
$\coz^{\alpha,\beta}_1\vmylabel{s7}{3}{7}$, and
$\cpo^{-\alpha,\beta}_{p}\vmylabel{Cc}{3}{0}$ is still an ``operator''
equation in the sense that it is written in terms of maps.  A
``scalar'' equation follows by applying~\eqref{op-eq} to the top
subquotient of $\repQ^{\alpha,\beta}_{1,0}$.  Simple analysis as in
Fig.~\ref{fig:mod-map} readily shows that, with
$\Omega_{\repQ^{\alpha,\beta}_{1,0},1}$ sending the top subquotient of
$\repQ^{\alpha,\beta}_{1,0}$ as
$\myxy{\oziAIBZ}\to\myxy{\poiieemAIBZm}$, some other relevant
subquotients of $\repQ^{\alpha,\beta}_{1,0}$ are mapped by
$\Omega_{\repQ^{\alpha,\beta}_{1,0},1}$ as
\begin{gather*}
  \myxy{\oziiisVAIBZ}\to-\half\myxy{\poivEEmAIBZm},\quad
  \myxy{\oziiimtIVAIBZ}\to\half \myxy{\poivEEmAIBZm},\quad
  \myxy{\oziiisVIAIBZ}\to0,\quad \myxy{\oziiisVIIAIBZ}\to0.
\end{gather*}
Also, $\Watypm^{-\alpha,\beta}_{p,0}\vmylabel{Cc}{3}{0}$ maps as
\begin{equation*}
  \myxy{\poiieemAIBZm}\to\alpha \myxy{\poivEEmAIBZm}.
\end{equation*}
This gives the equation
\begin{equation*}
  \alpha\cpo^{-\alpha,\beta}_{p}\vmylabel{Cc}{3}{0}
  +\ffrac{\coz^{\alpha,\beta}_1\vmylabel{mt4}{3}{4}}{2}
  -\ffrac{\coz^{\alpha,\beta}_1\vmylabel{s5}{3}{2}}{2}=0.
\end{equation*}

The full list of equations that follow from commuting with
$\Omega_{\repQ^{\alpha,+}_{1,0},n}$ is
\begin{alignat*}{2}
  \alpha\cpo^{-\alpha,+}_{p}\vmylabel{Cc}{3}{0}+\ffrac{\coz^{\alpha,+}_1\vmylabel{mt4}{3}{4}}{2}-\ffrac{\coz^{\alpha,+}_1\vmylabel{s5}{3}{2}}{2}&=0,
  &
  \alpha\cpo^{-\alpha,-}_{p}\vmylabel{Cc}{3}{0}-\ffrac{\coz^{\alpha,+}_1\vmylabel{mt4}{3}{4}}{2}-\ffrac{\coz^{\alpha,+}_1\vmylabel{s5}{3}{2}}{2}&=0,
  \\
  \cpp^{-\alpha,+}_{p-1}\vmylabel{Efb}{3}{1}-\ffrac{\coz^{\alpha,+}_1\vmylabel{mt4}{3}{4}}{2}+\ffrac{\coz^{\alpha,+}_1\vmylabel{s6}{3}{5}}{2}&=0,
  \quad&
  \cpp^{-\alpha,-}_{p-1}\vmylabel{Efb}{3}{1}+\ffrac{\coz^{\alpha,+}_1\vmylabel{mt4}{3}{4}}{2}+\ffrac{\coz^{\alpha,+}_1\vmylabel{s6}{3}{5}}{2}&=0,
  \\
  -\ffrac{\cso^{\alpha,-}_2\vmylabel{Bt}{3}{0}}{\qint{2}^3}-\ffrac{1}{2}\coz^{\alpha,+}_1\vmylabel{s7}{3}{7}-\alpha\ffrac{\cso^{\alpha,-}_2\vmylabel{Tb}{3}{5}}{\qint{2}^4}&=0,
  \\
  &&\kern-200pt \css^{\alpha,+}_1\vmylabel{Bc}{3}{0}-\ffrac{\alpha}{2}
  \qint{2}\coz^{\alpha,+}_1\vmylabel{s5}{3}{2}+\ffrac{1}{2}\qint{2}\coz^{\alpha,+}_1\vmylabel{s6}{3}{5}+\ffrac{1}{2}\qint{2}\coz^{\alpha,+}_1\vmylabel{s7}{3}{7}&=0
\end{alignat*}
(note the four $\alpha$ occurring as coefficients).

Similarly, the $4(p-2)$ equations that follow from commuting with
$\Omega_{\repQ^{+,+}_{s,s},n}$, $1\leq s\leq p-2$, are as follows:
\begin{alignat*}{2}
  \ffrac{\css^{+,+}_1\vmylabel{Bc}{3}{0}}{\qint{2}}+\ffrac{\coz^{+,+}_1\vmylabel{s5}{3}{2}}{2}+\ffrac{\coz^{+,+}_1\vmylabel{s6}{3}{5}}{2}+\ffrac{\coz^{+,+}_1\vmylabel{s7}{3}{7}}{2}&=0,
  \\
  \css^{+,+}_s\vmylabel{Cef}{3}{5}+\ffrac{\cso^{-,-}_{p-s}\vmylabel{Tb}{3}{5}}{\qint{s}^2}&=0,
  \quad &&1\leq s\leq p-2,
  \\
  \css^{+,+}_s\vmylabel{Cef}{3}{5}-\ffrac{\cso^{-,+}_{p-s}\vmylabel{Tb}{3}{5}}{\qint{s}^2}&=0,
  \quad &&1\leq s\leq p-2,
  \\
  \ffrac{\css^{+,+}_{s+1}\vmylabel{Bc}{3}{0}}{\qint{s+2}}-\ffrac{\css^{+,-}_s\vmylabel{Bc}{3}{0}}{\qint{s}}+\ffrac{\css^{+,-}_s\vmylabel{Cef}{3}{5}}{\qint{s}^2}&=0,
  \quad &&1\leq s\leq p-3,
  \\
  -\ffrac{\css^{+,-}_{s+1}\vmylabel{Bc}{3}{0}}{\qint{s+2}}+\ffrac{\css^{+,+}_s\vmylabel{Bc}{3}{0}}{\qint{s}}+\ffrac{\css^{+,+}_s\vmylabel{Cef}{3}{5}}{\qint{s}^2}&=0,
  \quad &&1\leq s\leq p-3, \\
  \ffrac{\css^{+,+}_{p-2}\vmylabel{Bc}{3}{0}}{\qint{2}}+\ffrac{\css^{+,+}_{p-2}\vmylabel{Cef}{3}{5}}{\qint{2}^2}+\cpp^{+,-}_{p-1}\vmylabel{Efb}{3}{1}&=0
\end{alignat*}
(the projective modules with a level-2 subquotient isomorphic to
$\repZ^{+,+}_{s,s}$ (the top subquotient of $\repQ^{+,+}_{s,s}$) are
$\repQ^{-,-}_{p-1,0}$, $\repQ^{+,+}_{1,0}$, $\repQ^{-,+}_{p-1,0}$, and
$\repQ^{+,-}_{2,2}$ for $s=1$ and $\repQ^{-,-}_{p-s,0}$,
$\repQ^{-,+}_{p-s,0}$, $\repQ^{+,-}_{s-1,s-1}$, and
$\repQ^{+,-}_{s+1,s+1}$ for $2\leq s\leq p-2$).

It is impossible to write the entire system of equations here because
of its length. Most of the equations have two or three terms, but the
system consists of numerous ``blocks'' in accordance with the
parameterization $\Omega_{\repQ^{\alpha,\beta}_{r,s},n}$ of basic
algebra generators.  Part of the system can be written ``uniformly,''
with the equations labeled by $s$ and having the same functional form
for any $s$.  But there are also ``boundary effects'': the two
length-$(p-2)$ series of projective modules ($\zss$ and $\zso$) are
followed at the end by modules of a somewhat reduced structure ($\zpp$
and $\zpo$), and are also ``joined'' by the $\zoz$ projective module
with an ``enhanced'' structure.  Both these effects are well seen in
the above formulas ($\coz_1$ in the first and $\cpp_{p-1}$ in the last
equation).

But it is possible to give the full solution of the system, which can
be written relatively compactly.  The comparative complexity or
simplicity of the explicit solution---and indeed of the procedure of
solving---depends rather strongly on the choice of free variables in
terms of which the others are to be expressed.  In choosing the free
variables, we were guided by the desire to avoid final formulas with
the number of terms growing with $p$; this turned out to be possible,
and was actually a factor underlying the success in solving the system
explicitly.  Numerous variations of our choice are of course possible.
A drawback of the specific choice that we make is that the formulas
become slightly sensitive to the parity of~$p$; we therefore write the
solution explicitly only for odd $p$.  Specifically, the $2p+1$ free
variables are chosen as
\begin{align*}
  &\coz^{+,+}_1\vmylabel{s5}{3}{2},\quad
  \coz^{+,+}_1\vmylabel{mt4}{3}{4},\quad
  \coz^{-,-}_1\vmylabel{mt4}{3}{4},\quad
  \coz^{+,+}_1\vmylabel{s6}{3}{5},\quad
  \cso^{-,-}_2\vmylabel{Bt}{3}{0},\quad
  \\
  &\cso^{+,-}_{2i+1}\vmylabel{Bt}{3}{0},\quad i =
  1,\dots,\fffrac{p-1}{2}-1, \qquad
  \cso^{+,+}_{2i}\vmylabel{Bt}{3}{0},\quad i =
  1,\dots,\fffrac{p-1}{2},
  \\
  &\css^{+,-}_{2i-1}\vmylabel{Bc}{3}{0},\quad i =
  1,\dots,\fffrac{p-1}{2}-1, \qquad
  \css^{+,+}_{2i}\vmylabel{Bc}{3}{0},\quad i =
  1,\dots,\fffrac{p-1}{2}-1,
  \\
  &\cpp^{+,+}_{p-1}\vmylabel{Efb}{3}{1}.
\end{align*}

The other $14 p - 9$ variables are expressed in terms of these as
follows.  First, there are $22$ lower-$s$ relations, occurring because
of the special structure of the $\repQ^{\alpha,\beta}_{1,0}$
projective module:
\begin{align*}
  \coz^{-,+}_1\vmylabel{mt4}{3}{4}&=
  -\coz^{-,-}_1\vmylabel{mt4}{3}{4},
  \\
  \coz^{+,-}_1\vmylabel{mt4}{3}{4}&=
  -\coz^{+,+}_1\vmylabel{mt4}{3}{4},
  \\
  \coz^{-,\beta}_1\vmylabel{s5}{3}{2}&=\coz^{-,-}_1\vmylabel{mt4}{3}{4}
  - \ffrac{2 \cso^{+,+}_{p - 1}\vmylabel{Bt}{3}{0}}{\qint{2}},
  \\
  \coz^{+,-}_1\vmylabel{s5}{3}{2}&=\coz^{+,+}_1\vmylabel{s5}{3}{2},
  \\
  \coz^{-,\beta}_1\vmylabel{s6}{3}{5}&=
  -\coz^{-,-}_1\vmylabel{mt4}{3}{4} - 2\cpp^{+,+}_{p -
    1}\vmylabel{Efb}{3}{1},
  \\
  \coz^{+,-}_1\vmylabel{s6}{3}{5}&=\coz^{+,+}_1\vmylabel{s6}{3}{5},
  \\
  \coz^{+,+}_1\vmylabel{s7}{3}{7}&= -\coz^{+,+}_1\vmylabel{s5}{3}{2} -
  \coz^{+,+}_1\vmylabel{s6}{3}{5} -
  \ffrac{2\css^{+,-}_1\vmylabel{Bc}{3}{0}}{\qint{2}} +
  2\coz^{-,-}_1\vmylabel{mt4}{3}{4},
  \\
  \coz^{+,-}_1\vmylabel{s7}{3}{7}&= -\coz^{+,+}_1\vmylabel{s5}{3}{2} -
  \coz^{+,+}_1\vmylabel{s6}{3}{5} -
  \ffrac{2\css^{+,-}_1\vmylabel{Bc}{3}{0}}{\qint{2}},
  \\
  \coz^{-,\beta}_1\vmylabel{s7}{3}{7}&= -\beta
  \coz^{+,+}_1\vmylabel{mt4}{3}{4} + \coz^{+,+}_1\vmylabel{s5}{3}{2} +
  2\cpp^{+,+}_{p - 1}\vmylabel{Efb}{3}{1} +
  \ffrac{2\css^{+,-}_1\vmylabel{Bc}{3}{0}}{\qint{2}},
  \\
  \css^{-,\beta}_1\vmylabel{Bc}{3}{0}&=\ffrac{\beta}{2}
  \qint{2}\coz^{+,+}_1\vmylabel{mt4}{3}{4} -
  \half\qint{2}\coz^{+,+}_1\vmylabel{s5}{3}{2} + \cso^{+,+}_{p -
    1}\vmylabel{Bt}{3}{0} - \css^{+,-}_1\vmylabel{Bc}{3}{0},
  \\
  \cso^{-,\beta}_2\vmylabel{Tb}{3}{5}&=\ffrac{\beta}{2} \qint{2}^4
  \coz^{+,+}_1\vmylabel{mt4}{3}{4} - \ffrac{\beta}{2}\qint{2}^4
  \coz^{+,+}_1\vmylabel{s5}{3}{2} - \beta\qint{2}^4
  \cpp^{+,+}_{p - 1}\vmylabel{Efb}{3}{1} \\
  &\quad{}-\beta\qint{2}^3 \css^{+,-}_1\vmylabel{Bc}{3}{0} -
  \beta\qint{2} \cso^{-,-}_2\vmylabel{Bt}{3}{0},
  \\
  \cso^{+,\beta}_2\vmylabel{Tb}{3}{5}&=
  -\ffrac{\beta}{2}\qint{2}^4\coz^{+,+}_1\vmylabel{s5}{3}{2} -
  \ffrac{\beta}{2}\qint{2}^4 \coz^{+,+}_1\vmylabel{s6}{3}{5} -
  \beta\qint{2}^3 \css^{+,-}_1\vmylabel{Bc}{3}{0} + \beta\qint{2}
  \cso^{+,+}_2\vmylabel{Bt}{3}{0},
  \\
  \cso^{-,+}_2\vmylabel{Bt}{3}{0}&=\cso^{-,-}_2\vmylabel{Bt}{3}{0} -
  \qint{2}^3 \coz^{+,+}_1\vmylabel{mt4}{3}{4},
  \\
  \cso^{-,-}_3\vmylabel{Bt}{3}{0}&= -\half \qint{2}\qint{3}^3
  \coz^{+,+}_1\vmylabel{mt4}{3}{4} - \half \qint{2}\qint{3}^3
  \coz^{+,+}_1\vmylabel{s5}{3}{2} - \qint{3}^3
  \css^{+,-}_1\vmylabel{Bc}{3}{0} + \qint{3}^2 \css^{+,+}_{p -
    3}\vmylabel{Bc}{3}{0},
  \\
  \cso^{-,\beta}_3\vmylabel{Tb}{3}{5}&= -\ffrac{\beta}{2}
  \qint{2}\qint{3}^4 \coz^{+,+}_1\vmylabel{mt4}{3}{4} +
  \ffrac{\beta}{2} \qint{2}\qint{3}^4
  \coz^{+,+}_1\vmylabel{s5}{3}{2} + \beta\qint{3}^4 \css^{+,-}_1\vmylabel{Bc}{3}{0}\\
  &\quad{}- \beta\qint{3}^3 \css^{+,+}_{p - 3}\vmylabel{Bc}{3}{0} +
  \beta\ffrac{\qint{3}^4 \cso^{-,-}_2\vmylabel{Bt}{3}{0}}{\qint{2}^2}.
\end{align*}
Here and hereafter, $\beta=\pm$. \ Next, there are $14(p-1)-31$
``serial'' relations for ``generic'' values of $s$, which in our
solution are split into even and odd ones:
\begin{alignat*}{2}
  \cso^{+,\beta}_{2 i}\vmylabel{Tb}{3}{5}&=\beta\qint{2 i}
  \cso^{+,+}_{2 i}\vmylabel{Bt}{3}{0} - \beta\ffrac{\qint{2 i}^4
    \cso^{+,-}_{2 i - 1}\vmylabel{Bt}{3}{0}}{\qint{2 i - 2} \qint{2 i
      - 1}^2}, && 2\leq i\leq \ffrac{p - 1}{2},
  \\
  \cso^{+,\beta}_{2 i + 1}\vmylabel{Tb}{3}{5}&=\beta\qint{2 i +
    1}\cso^{+,-}_{2 i + 1}\vmylabel{Bt}{3}{0} - \beta\ffrac{\qint{2 i
      + 1}^4\cso^{+,+}_{2 i}\vmylabel{Bt}{3}{0}}{\qint{2 i - 1}\qint{2
      i}^2}, && 1\leq i\leq \ffrac{p - 1}{2} - 1,
  \\
  \cso^{-,\beta}_{2 i}\vmylabel{Tb}{3}{5}&= -\beta\qint{2 i}^3
  \css^{+,-}_{p - 2 i}\vmylabel{Bc}{3}{0} + \beta\ffrac{\qint{2 i}^4
    \css^{+,+}_{p + 1 - 2 i}\vmylabel{Bc}{3}{0}}{\qint{2 i - 2}}, &&
  2\leq i\leq \ffrac{p - 1}{2},
  \\
  \cso^{-,\beta}_{2 i + 1}\vmylabel{Tb}{3}{5}&= -\beta\qint{2 i + 1}^3
  \css^{+,+}_{p - 1 - 2 i}\vmylabel{Bc}{3}{0} + \beta\ffrac{\qint{2 i
      + 1}^4\css^{+,-}_{p - 2 i}\vmylabel{Bc}{3}{0}}{\qint{2 i - 1}},
  && 2\leq i\leq \ffrac{p - 1}{2}-1,
  \\
  \cso^{+,-}_{2 i}\vmylabel{Bt}{3}{0}&=\cso^{+,+}_{2
    i}\vmylabel{Bt}{3}{0} - \qint{2 i - 1}\qint{2
    i}^3\coz^{-,-}_1\vmylabel{mt4}{3}{4}, && 1\leq i\leq \ffrac{p -
    1}{2},
  \\
  \cso^{+,+}_{2 i + 1}\vmylabel{Bt}{3}{0}&=\cso^{+,-}_{2 i +
    1}\vmylabel{Bt}{3}{0} - \qint{2 i}\qint{2 i + 1}^3
  \coz^{-,-}_1\vmylabel{mt4}{3}{4},&& 1\leq i\leq \ffrac{p - 1}{2} -
  1,
  \\
  \cso^{-,\beta}_{2 i}\vmylabel{Bt}{3}{0}&= -\ffrac{\beta}{2}\qint{2 i
    - 1} \qint{2 i}^3 \coz^{+,+}_1\vmylabel{mt4}{3}{4} - \half \qint{2
    i - 1} \qint{2 i}^3
  \coz^{+,+}_1\vmylabel{s5}{3}{2} \\
  &\quad{}- \ffrac{\qint{2 i - 1}\qint{2 i}^3
    \css^{+,-}_1\vmylabel{Bc}{3}{0}}{\qint{2}} + \qint{2 i}^2
  \css^{+,-}_{p - 2 i}\vmylabel{Bc}{3}{0}, && 2\leq i\leq \ffrac{p -
    1}{2},
  \\
  \cso^{-,\beta}_{2 i + 1}\vmylabel{Bt}{3}{0}&=\ffrac{\beta}{2}\qint{2 i}\qint{2 i + 1}^3 \coz^{+,+}_1\vmylabel{mt4}{3}{4} - \half \qint{2 i}\qint{2 i + 1}^3 \coz^{+,+}_1\vmylabel{s5}{3}{2} \\
  &\quad{} - \ffrac{\qint{2 i}\qint{2 i + 1}^3
    \css^{+,-}_1\vmylabel{Bc}{3}{0}}{\qint{2}} + \qint{2 i + 1}^2
  \css^{+,+}_{p - 1 - 2 i}\vmylabel{Bc}{3}{0}, && 1\leq i\leq \ffrac{p
    - 1}{2} - 1,
  \\
  \css^{-,\beta}_{2 i}\vmylabel{Bc}{3}{0}&= -\ffrac{\beta}{2}\qint{2 i} \qint{2 i + 1}\coz^{+,+}_1\vmylabel{mt4}{3}{4} - \half\qint{2 i} \qint{2 i + 1}\coz^{+,+}_1\vmylabel{s5}{3}{2} \kern-100pt\\
  &\quad{}+ \ffrac{\cso^{+,-}_{p - 2 i}\vmylabel{Bt}{3}{0}}{\qint{2
      i}^2} - \ffrac{\qint{2 i} \qint{2 i +
      1}\css^{+,-}_1\vmylabel{Bc}{3}{0}}{\qint{2}}, && 1\leq i\leq
  \ffrac{p - 1}{2} - 1,
  \\
  \css^{-,\beta}_{2 i + 1}\vmylabel{Bc}{3}{0}&=\ffrac{\beta}{2}\qint{2 i + 1} \qint{2 i + 2}\coz^{+,+}_1\vmylabel{mt4}{3}{4} - \half\qint{2 i + 1} \qint{2 i + 2}\coz^{+,+}_1\vmylabel{s5}{3}{2} \kern-100pt\\
  &\quad{}+ \ffrac{\cso^{+,+}_{p - 1 - 2
      i}\vmylabel{Bt}{3}{0}}{\qint{2 i + 1}^2} - \ffrac{\qint{2 i +
      1}\qint{2 i + 2}\css^{+,-}_1\vmylabel{Bc}{3}{0}}{\qint{2}}, &&
  1\leq i\leq \ffrac{p - 1}{2} - 1,
  \\
  \css^{+,-}_{2 i}\vmylabel{Bc}{3}{0}&=\css^{+,+}_{2
    i}\vmylabel{Bc}{3}{0} - \qint{2 i} \qint{2 i +
    1}\coz^{-,-}_1\vmylabel{mt4}{3}{4}, && 1\leq i\leq \ffrac{p -
    1}{2} - 1,
  \\
  \css^{+,+}_{2 i + 1}\vmylabel{Bc}{3}{0}&=\css^{+,-}_{2 i +
    1}\vmylabel{Bc}{3}{0} - \qint{2 i + 1}\qint{2 i +
    2}\coz^{-,-}_1\vmylabel{mt4}{3}{4}, && 0\leq i\leq \ffrac{p -
    1}{2} - 2,
  \\
  \css^{-,\beta}_{2 i}\vmylabel{Cef}{3}{5}&=\beta\ffrac{\cso^{+,-}_{p
      - 2 i}\vmylabel{Bt}{3}{0}}{\qint{2 i}} - \beta\ffrac{\qint{2
      i}^2 \cso^{+,+}_{p - 1 - 2 i}\vmylabel{Bt}{3}{0}}{\qint{2 i +
      1}^2\qint{2 i + 2}}, && 1\leq i\leq \ffrac{p - 1}{2} - 1,
  \\
  \css^{-,\beta}_{2 i +
    1}\vmylabel{Cef}{3}{5}&=\beta\ffrac{\cso^{+,+}_{p - 1 - 2
      i}\vmylabel{Bt}{3}{0}}{\qint{2 i + 1}} - \beta\ffrac{\qint{2 i +
      1}^2\cso^{+,-}_{p - 2 - 2 i}\vmylabel{Bt}{3}{0}}{\qint{2 i +
      2}^2\qint{2 i + 3}}, && 0\leq i\leq \ffrac{p - 1}{2} - 2.
  \\
  \css^{+,\beta}_{2 i}\vmylabel{Cef}{3}{5}&=\beta\ffrac{\qint{2
      i}^2\css^{+,-}_{2 i + 1}\vmylabel{Bc}{3}{0}}{\qint{2 i + 2}} -
  \beta\qint{2 i} \css^{+,+}_{2 i}\vmylabel{Bc}{3}{0}, && 1\leq i\leq
  \ffrac{p - 1}{2} - 2,
  \\
  \css^{+,\beta}_{2 i + 1}\vmylabel{Cef}{3}{5}&=\beta\ffrac{\qint{2 i
      + 1}^2\css^{+,+}_{2 i + 2}\vmylabel{Bc}{3}{0}}{\qint{2 i + 3}} -
  \beta\qint{2 i + 1} \css^{+,-}_{2 i + 1}\vmylabel{Bc}{3}{0}, &&
  0\leq i\leq \ffrac{p - 1}{2} - 2,
\end{alignat*}
And finally, the $14$ high-$s$ relations, whose form is largely
determined by the somewhat special structure of
$\repQ^{\alpha,\beta}_{p,0}$ and $\repQ^{\alpha,\beta}_{p-1,p-1}$, are
\begin{align*}
  \css^{+,+}_{p - 3}\vmylabel{Cef}{3}{5}&= -\half \qint{2}\qint{3}^2 \coz^{+,+}_1\vmylabel{mt4}{3}{4} + \half \qint{2}\qint{3}^2\coz^{+,+}_1\vmylabel{s5}{3}{2} + \ffrac{\qint{3}^2 \cso^{-,-}_2\vmylabel{Bt}{3}{0}}{\qint{2}^2}\\
  &\quad{}+ \qint{3}^2 \css^{+,-}_1\vmylabel{Bc}{3}{0} - \qint{3}
  \css^{+,+}_{p - 3}\vmylabel{Bc}{3}{0},
  \\
  \css^{+,\beta}_{p - 2}\vmylabel{Cef}{3}{5}&=\ffrac{\beta}{2} \qint{2}^2 \coz^{+,+}_1\vmylabel{mt4}{3}{4} - \ffrac{\beta}{2}\qint{2}^2 \coz^{+,+}_1\vmylabel{s5}{3}{2} - \beta\qint{2}^2 \cpp^{+,+}_{p - 1}\vmylabel{Efb}{3}{1}\\
  &\quad{}- \beta\qint{2} \css^{+,-}_1\vmylabel{Bc}{3}{0} -
  \beta\ffrac{\cso^{-,-}_2\vmylabel{Bt}{3}{0}}{\qint{2}},
  \\
  \css^{-,\beta}_{p - 2}\vmylabel{Cef}{3}{5}&= -\ffrac{\beta}{2}
  \qint{2}^2 \coz^{+,+}_1\vmylabel{s5}{3}{2} -
  \ffrac{\beta}{2}\qint{2}^2 \coz^{+,+}_1\vmylabel{s6}{3}{5} -
  \beta\qint{2} \css^{+,-}_1\vmylabel{Bc}{3}{0} +
  \beta\ffrac{\cso^{+,+}_2\vmylabel{Bt}{3}{0}}{\qint{2}},
  \\
  \css^{+,+}_{p - 2}\vmylabel{Bc}{3}{0}&= -\half
  \qint{2}\coz^{+,+}_1\vmylabel{mt4}{3}{4} +
  \half\qint{2}\coz^{+,+}_1\vmylabel{s5}{3}{2} +
  \ffrac{\cso^{-,-}_2\vmylabel{Bt}{3}{0}}{\qint{2}^2} +
  \css^{+,-}_1\vmylabel{Bc}{3}{0} -
  \qint{2}\coz^{-,-}_1\vmylabel{mt4}{3}{4},
  \\
  \css^{+,-}_{p - 2}\vmylabel{Bc}{3}{0}&= -\half
  \qint{2}\coz^{+,+}_1\vmylabel{mt4}{3}{4} +
  \half\qint{2}\coz^{+,+}_1\vmylabel{s5}{3}{2} +
  \ffrac{\cso^{-,-}_2\vmylabel{Bt}{3}{0}}{\qint{2}^2} +
  \css^{+,-}_1\vmylabel{Bc}{3}{0},
  \\
  \cpp^{-,\beta}_{p - 1}\vmylabel{Efb}{3}{1}&=\ffrac{\beta}{2}
  \coz^{+,+}_1\vmylabel{mt4}{3}{4} - \half
  \coz^{+,+}_1\vmylabel{s6}{3}{5},
  \\
  \cpp^{+,-}_{p -
    1}\vmylabel{Efb}{3}{1}&=\coz^{-,-}_1\vmylabel{mt4}{3}{4} +
  \cpp^{+,+}_{p - 1}\vmylabel{Efb}{3}{1},
  \\
  \cpo^{+,+}_p\vmylabel{Cc}{3}{0}&=\ffrac{\cso^{+,+}_{p -
      1}\vmylabel{Bt}{3}{0}}{\qint{2}} -
  \coz^{-,-}_1\vmylabel{mt4}{3}{4},
  \\
  \cpo^{+,-}_p\vmylabel{Cc}{3}{0}&=\ffrac{\cso^{+,+}_{p -
      1}\vmylabel{Bt}{3}{0}}{\qint{2}}
  \\
  \cpo^{-,\beta}_p\vmylabel{Cc}{3}{0}&=\half
  \coz^{+,+}_1\vmylabel{s5}{3}{2} - \ffrac{\beta}{2}
  \coz^{+,+}_1\vmylabel{mt4}{3}{4}.
\end{align*}
This completes the list of formulas expressing $14 p - 9$ variables in
terms of $2p+1$ independent variables and thus solving linear
system~\eqref{the-eqs}.

\parindent0pt


\begin{thebibliography}{99}
\bibitem{[Gurarie]} V.~Gurarie, \textit{Logarithmic operators in
    conformal field theory,} Nucl.\ Phys.\ B410 (1993) 535
  [hep-th$/$\linebreak[0]9303160].

\bibitem{[Sa]} H.~Saleur, \textit{Polymers and percolation in
    two-dimensions and twisted $N=2$ supersymmetry}, Nucl.\ Phys.
  B382 (1992) 486--531 [hep-th$/$\linebreak[0]9111007].

\bibitem{[Kausch]}H.G.\;Kausch, \textit{Extended conformal algebras
    generated by a multiplet of primary fields}, Phys.\ Lett. B\;259
  (1991) 448.

\bibitem{[GK+]}M.R.\;Gaberdiel and H.G.\;Kausch,
  \textit{Indecomposable fusion products}, Nucl.\ Phys.\ B477 (1996)
  293--318 [hep-th$/$\linebreak[0]9604026]; \textit{A rational
    logarithmic conformal field theory}, Phys.\ Lett. B\;386 (1996)
  131--137 [hep-th$/$\linebreak[0]9606050]; \textit{A local
    logarithmic conformal field theory}, Nucl.\ Phys.\ B538 (1999)
  631--658 [hep-th$/$\linebreak[0]9807091].

\bibitem{[FHST]}J.\;Fuchs, S.\;Hwang, A.M.\;Semikhatov, and
  I.Yu.\;Tipunin, \textit{Nonsemisimple fusion algebras and the
    Verlinde formula}, Commun.\ Math.\ Phys.\ 247 (2004) 713--742
  [hep-th$/$\linebreak[0]0306274].

\bibitem{[STbr]}A.M.\;Semikhatov and I.Yu.\;Tipunin, \textit{The
    Nichols algebra of screenings}, Commun. Contemp. Math. 14 (2012)
  1250029, arXiv:1101.5810.

\bibitem{[Nich]}W.\;D.\;Nichols, \textit{Bialgebras of type one},
  Commun. Algebra 6 (1978) 1521--1552.

\bibitem{[Wor]}S.L.\;Woronowicz, \textit{Differential calculus on
    compact matrix pseudogroups \textup{(}quantum groups\textup{)}},
  Commun. Math. Phys. 122 (1989) 125--170.

\bibitem{[Lu-intro]}G.\;Lusztig, \textsl{Introduction to Quantum
    Groups}.  Birkh\"auser, 1993.

\bibitem{[Rosso-inv]}M.\;Rosso, \textit{Quantum groups and quantum
    shuffles}, Invent. math. 133 (1998) 399--416.

\bibitem{[AG]}N.\;Andruskiewitsch and M.\;Gra\~na, \textit{Braided
    Hopf algebras over non-abelian finite groups}, Bol. Acad. Nacional
  de Ciencias (Cordoba) 63 (1999) 45--78 [arXiv:math$/$9802074
  [math.QA]].

\bibitem{[AS-pointed]}N.\;Andruskiewitsch and H.-J.\;Schneider,
  \textit{Pointed Hopf algebras},
  in: \textsl{New directions in Hopf algebras}, MSRI Publications 43,
  pages 1--68.  Cambridge University Press, 2002.

\bibitem{[AS-onthe]}N.\;Andruskiewitsch and H.-J.\;Schneider,
  \textit{On the classification of finite-dimensional pointed Hopf
    algebras}, Ann. Math. 171 (2010) 375--417 [arXiv:math$/$0502157
  [math.QA]].

\bibitem{[Heck-Weyl]}I.\;Heckenberger, \textit{The Weyl groupoid of a
    Nichols algebra of diagonal type}, Invent. Math. 164 (2006)
  175--188.

\bibitem{[Heck-class]}I.\;Heckenberger, \textit{Classification of
    arithmetic root systems}, Adv. Math. 220 (2009) 59--124
  [math.QA$/$\linebreak[0]0605795].

\bibitem{[AHS]} N.\;Andruskiewitsch, I.\;Heckenberger, and
  H.-J.\;Schneider, \textit{The Nichols algebra of a semisimple
    Yet\-ter--Drinfeld module}, Amer. J. Math. 132 (2010) 1493--1547.

\bibitem{[ARS]}N.\;Andruskiewitsch, D\;Radford, and H.-J.\;Schneider,
  \textit{Complete reducibility theorems for modules over pointed Hopf
    algebras}, J.\ Algebra, 324 (2010) 2932--2970 [arXiv:1001.3977].

\bibitem{[Ag-0804-standard]}I.E.\;Angiono, \textit{On Nichols algebras
    with standard braiding}, Algebra \& Number Theory 3 (2009) 35--106
  [arXiv:\linebreak[0]0804.\linebreak[0]0816].

\bibitem{[Ag-1008-presentation]}I.E.\;Angiono, \textit{A presentation
    by generators and relations of Nichols algebras of diagonal type
    and convex orders on root systems}, arXiv:1008.4144,
  J. Europ. Math. Soc.

\bibitem{[KL]} D.~Kazhdan and G.~Lusztig, \textit{Tensor structures
    arising from affine Lie algebras,} I, J. Amer. Math. Soc. 6 (1993)
  905--947; II, J. Amer. Math. Soc. 6 (1993) 949--1011; III, J.\
  Amer.\ Math.\ Soc.\ 7 (1994) 335--381; IV, J.\ Amer.\ Math.\ Soc.\ 7
  (1994) 383--453.

\bibitem{[Fink]}M.~Finkelberg, \textit{An equivalence of fusion
    categories}, Geometric and Functional Analysis (GAFA) 6 (1996)
  249--267.

\bibitem{[FGST2]}B.L.\;Feigin, A.M.\;Gainutdinov, A.M.\;Semikhatov,
  and I.Yu.\;Tipunin, \textit{Kazhdan--Lusztig correspondence for the
    representation category of the triplet $W$-algebra in logarithmic
    CFT}, Theor. Math. Phys. 148 (2006) 1210--1235
  [arXiv:\linebreak[0]math$/$0512621 [math.QA]].

\bibitem{[NT]}K.\;Nagatomo and A.\;Tsuchiya, \textit{The triplet
    vertex operator algebra $W(p)$ and the restricted quantum group at
    root of unity}, arXiv:\linebreak[0]0902.4607 [math.QA].

\bibitem{[AM-lattice]}D.\;Adamovi\'c and A.\;Milas, \textit{Lattice
    construction of logarithmic modules for certain vertex algebras},
  Selecta Math. New Ser., 15 (2009) 535--561, arXiv:0902.3417.

\bibitem{[TW1]}A.\;Tsuchiya and S.\;Wood, \textit{The tensor structure
    on the representation category of the $W_p$ triplet algebra},
  arXiv:\linebreak[0]1201.0419.

\bibitem{[b-fusion]}A.M.\;Semikhatov, \textit{Fusion in the entwined
    category of Yetter--Drinfeld modules of a rank-$1$ Nichols
    algebra}, Theor. Math. Phys. 173(1) (2012) 1329--1358.

\bibitem{[FGST]}B.L.\;Feigin, A.M.\;Gainutdinov, A.M.\;Semikhatov, and
  I.Yu.\;Tipunin, \textit{Modular group representations and fusion in
    logarithmic conformal field theories and in the quantum group
    center}, Commun. Math. Phys. 265 (2006) 47--93
  [arXiv:\linebreak[0]hep-th$/$\linebreak[0]0504093].

\bibitem{[FGST3]}B.L.\;Feigin, A.M.\;Gainutdinov, A.M.\;Semikhatov,
  and I.Yu.\;Tipunin, \textit{Logarithmic extensions of minimal
    models: characters and modular transformations}, Nucl. Phys. B757
  (2006) 303--343 [arXiv:\linebreak[0]hep-th$/$\linebreak[0]0606196].

\bibitem{[FGST-q]} B.L.~Feigin, A.M.~Gainutdinov, A.M.~Semikhatov, and
  I.Yu.~Tipunin, \textit{Kazhdan--Lusztig-dual quantum group for
    logarithmic extensions of Virasoro minimal models}, J.\ Math.\
  Phys.\ 48 (2007) 032303 [math.QA$/$\linebreak[0]0606506].

\bibitem{[GRW]} M.R.\;Gaberdiel, I.\;Runkel, and S.\;Wood,
  \textit{Fusion rules and boundary conditions in the $c=0$ triplet
    model}, arXiv:\linebreak[0]0905.0916; \textit{A modular invariant
    bulk theory for the $c=0$ triplet model}, J. Phys. A44 (2011)
  015204, arXiv:\linebreak[0]1008.0082.

\bibitem{[RGW]} I.\;Runkel, M.R.\;Gaberdiel, and S.\;Wood,
  \textit{Logarithmic bulk and boundary conformal field theory and the
    full centre construction}, arXiv:1201.6273.

\bibitem{[c-charge]}A.M.\;Semikhatov, \textit{Virasoro central charges
    for Nichols algebras}, arXiv:1109.1767 [math.QA], in: ``Conformal
  field theories and tensor categories,'' Mathematical Lectures from
  Peking University Bai, C.; Fuchs, J.; Huang, Y.-Z.; Kong, L.;
  Runkel, I.; Schweigert, C. (Eds.)  2014, IX 67--92.

\bibitem{[nich-sl2-1]}A.M.\;Semikhatov and I.Yu.\;Tipunin,
  \textit{Logarithmic $\hSL{2}$ CFT{} models from Nichols
    algebras. 1}, J. Phys. A: Math. Theor. 46 (2013) 494011
  [arXiv:1301.2235].

\bibitem{[BFGT]} P.V.\;Bushlanov, B.L.\;Feigin, A.M.\;Gainutdinov, and
  I.Yu.\;Tipunin \textit{Lusztig limit of quantum sl(2) at root of
    unity and fusion of (1,p) Virasoro logarithmic minimal models},
  Nucl. Phys. B818 (2009) 179--195 [arXiv:0901.1602];

\bibitem{[BGT]} P.V.\;Bushlanov, A.M.\;Gainutdinov, and
  I.Yu.\;Tipunin, \textit{Kazhdan-Lusztig equivalence and fusion of
    Kac modules in Virasoro logarithmic models}, arXiv:1102.0271.

\bibitem{[AY]} I.~Angiono and H.~Yamane, \textit{The R-matrix of
    quantum doubles of Nichols algebras of diagonal type},
  arXiv:1304.5752 [math.QA].

\bibitem{[HY]}I.\;Heckenberger and H.\;Yamane, \textit{Drinfel'd
    doubles and Shapovalov determinants}, Revista de la Uni\'on
  Matem\'atica Argentina 51 no.~2 (2010) 107--146 [arXiv:0810.1621
  [math.QA]].

\bibitem{[HS-double]} I.\;Heckenberger and H.-J.\;Schneider,
  \textit{Yetter--Drinfeld modules over bosonizations of dually paired
    Hopf algebras}, arXiv:\linebreak[0]1111.4673.

\bibitem{[Radford-bos]}D.\;Radford, \textit{Hopf algebras with a
    projection}, J. Algebra 92 (1985) 322--347.

\bibitem{[Mj-braided]} S.~Majid, \textit{Braided groups and algebraic
    quantum field theories}, Lett. Math. Phys. 22 (1991) 167--176.

\bibitem{[Mj-trans-rank]} S.~Majid, \textit{Transmutation theory and
    rank for quantum braided groups}, Math. Proc. Camb. Phil. Soc. 113
  (1993) 45--70.

\bibitem{[Besp-TMF+]}Yu.N.\;Bespalov, \textit{Crossed modules, quantum
    braided groups, and ribbon structures}, Theor.\ Math.\ Phys. 103
  (1995) 621--637; \textit{Crossed modules and quantum groups in
    braided categories}, arXiv:q-alg/9510013.

\bibitem{[Majid-bos]}S.\;Majid, \textit{Crossed products by braided
    groups and bosonization}, J. Algebra 163 (1994) 165--190.

\bibitem{[Mj]} S.~Majid, \textit{Doubles of quasitriangular Hopf
    algebras}, Comm. Alg. 19 (1991) 3061--3073.

\bibitem{[McL]} S.~Maclane, \textsl{Homology}, Springer Verlag, 1963.

\bibitem{[Arike]}Y.~Arike, \textit{Symmetric linear functions of the
    restricted quantum group $\bar{U}_qsl_2(\mathbb{C})$},
  arXiv:0706.1113 [math.QA]; \textit{A construction of symmetric
    linear functions of the restricted quantum group $\overline{U}_q
    (sl_2)$}, arXiv:0807.0052 [math.QA].

\bibitem{[GV]}A.M.\;Gainutdinov and R.\;Vasseur, \textit{Lattice
    fusion rules and logarithmic operator product expansions},
  Nucl. Phys. B 868, 223--270 (2013) [arXiv:1203.6289].

\bibitem{[GSaT]} A.M.\;Gainutdinov, H.\;Saleur, and I.Yu.\;Tipunin,
  \textit{Lattice W-algebras and logarithmic CFTs},
  arXiv:\linebreak[0]1212.\linebreak[0]1378.

\bibitem{[GJRSV]} A.M.\;Gainutdinov, J.L.\;Jacobsen, N.\;Read,
  H.\;Saleur, and R.~Vasseur, \textit{Logarithmic conformal field
    theory: a lattice approach}, J. Phys. A: Math. Theor. 46 (2013)
  494012 [arXiv:1303.2082].

\bibitem{[FSS-1106]} J.\;Fuchs, C.\;Schweigert, and C.\;Stigner,
  \textit{Modular invariant Frobenius algebras from ribbon Hopf
    algebra automorphisms}, Journal of Algebra 363 (2012) 29--72
  [arXiv:1106.0210 [math.QA]].

\bibitem{[FSS-1207]} J.\;Fuchs, C.\;Schweigert, and C.\;Stigner,
  \textit{Higher genus mapping class group invariants from
    factorizable Hopf algebras}, ZMP-HH/12-13, Hamburger Beitr. zur
  Math. 447 [arXiv:1207.6863 [math.QA]].

\bibitem{[FSS-1302]} J.\;Fuchs, C.\;Schweigert, and C.\;Stigner,
  \textit{From non-semisimple Hopf algebras to correlation functions
    for logarithmic Conformal Field Theory}, Hamburger Beitr. zur
  Mathematik 468, ZMP-HH/13-2 [arXiv:1302.4683 [hep-th]].

\end{thebibliography}
\end{document}